\documentclass[12pt,twoside]{article}

\usepackage{amssymb}
\usepackage{amsmath}
\usepackage{a4wide}
\usepackage{epic,eepic}
\usepackage{bbm}
\usepackage{amsthm}
\usepackage[all]{xy}

\newcommand{\Hom}{\operatorname{Hom}}
\newcommand{\lo}{\left(}
\newcommand{\po}{\right)}
\newcommand{\lk}{\left\{}
\newcommand{\pk}{\right\}}
\newcommand{\ld}{\left<}
\newcommand{\pd}{\right>}
\newcommand{\lkw}{\left[}
\newcommand{\pkw}{\right]}
\newcommand{\ulmod}{\operatorname{\underline{mod}}}
\newcommand{\ulhom}{\operatorname{\underline{Hom}}}
\newcommand{\End}{\operatorname{End}}
\newcommand{\gldim}{\operatorname{gldim}}

\newcommand{\rdim}{\operatorname{repdim}}

\newcommand{\md}{\operatorname{mod}}
\newcommand{\soc}{\operatorname{soc}}
\newcommand{\rad}{\operatorname{rad}}

\newcommand{\add}{\operatorname{add}}

\newcommand{\dualnosc}{\operatorname{D}}
\newcommand{\slad}{\operatorname{Tr}}
\newcommand{\id}{\operatorname{id}}
\renewcommand{\tilde}{\widetilde}

\makeatletter

\renewcommand{\section}{\@startsection{section}{0}{0mm}%
  {\baselineskip}{3ex}{\bf}}

\newcommand{\circlearrowl}{
  \put(0,0){\oval(20,20)[l]}
  \bezier{50}(0,-10)(4.14,-10)(7.1,-7.1)
  \bezier{50}(0,10)(4.14,10)(7.1,7.1)
  \bezier{50}(7,7)(8.48,5.7)(9.24,3.83)
  \put(9.7,-3.5){\vector(1,2){0}}
}
\newcommand{\circlearrowr}{
  \put(0,0){\oval(20,20)[r]}
  \bezier{50}(0,-10)(-4.14,-10)(-7.1,-7.1)
  \bezier{50}(0,10)(-4.14,10)(-7.1,7.1)
  \bezier{50}(-7,-7)(-8.48,-5.7)(-9.24,-3.83)
  \put(-8.7,3.3){\vector(-1,-2){0}}
}

\renewenvironment{proof}{
  \trivlist
  \item[\hskip\labelsep{\noindent\textbf{Proof.}}]
}{
  \qed
  \endtrivlist
}

\makeatother

\makeatletter
\newtheoremstyle{tw}{}{}{\itshape}{\parindent}{\scshape}{.}{ }{}
\newtheoremstyle{tw1}{}{}{\itshape}{\parindent}{\scshape}{.}{ }{}
\makeatother

\newtheorem{tw}{Theorem}[section]
\newtheorem{lemma}[tw]{Lemma}
\newtheorem{Notheorem}[tw]{\hspace{-.33em}}

\newtheorem{corollary}[tw]{Corollary}

\newtheorem{proposition}[tw]{Proposition}

\newtheorem{exa}[tw]{Example}

\newtheorem{tw1}{Theorem}

\newenvironment{arablist}{
  \begin{list}{\textup{(\arabic{enumi})}}{
    \usecounter{enumi}
    \setlength{\topsep}{0pt}
    \setlength{\parsep}{0pt}
    \setlength{\itemsep}{0pt}
    \setlength{\leftmargin}{0pt}
    \settowidth{\labelwidth}{(m)}
    \setlength{\itemindent}{\parindent}
    \addtolength{\itemindent}{\labelwidth}
    \addtolength{\itemindent}{\labelsep}
  }
}{
  \end{list}
}

\newenvironment{alphlist}{
  \begin{list}{\textup{(\alph{enumi})}}{
    \usecounter{enumi}
    \setlength{\topsep}{0pt}
    \setlength{\parsep}{0pt}
    \setlength{\itemsep}{0pt}
    \setlength{\leftmargin}{0pt}
    \settowidth{\labelwidth}{(m)}
    \setlength{\itemindent}{\parindent}
    \addtolength{\itemindent}{\labelwidth}
    \addtolength{\itemindent}{\labelsep}
  }
}{
  \end{list}
}

\newenvironment{romalist}{
  \begin{list}{\textup{(\roman{enumi})}}{
    \usecounter{enumi}
    \setlength{\topsep}{0pt}
    \setlength{\parsep}{0pt}
    \setlength{\itemsep}{0pt}
    \setlength{\leftmargin}{0pt}
    \settowidth{\labelwidth}{(vi)}
    \setlength{\itemindent}{\parindent}
    \addtolength{\itemindent}{\labelwidth}
    \addtolength{\itemindent}{\labelsep}
  }
}{
  \end{list}
}

\makeatletter

\def\ps@paper{
  \let\@mkboth\@gobbletwo
  \def\@evenhead{%
     \parbox{\textwidth}{
       \upshape\footnotesize\thepage\hfill%
       \textit{R. Bocian, T. Holm and A. Skowro\'{n}ski}%
       \hfill\rule[-3pt]{0pt}{15pt}}%
     }
  \def\@oddhead{%
    \parbox{\textwidth}{%
    \upshape\footnotesize\hfill%
    \textit{Derived equivalence classification of
one-parametric selfinjective algebras}%
    \hfill\thepage\rule[-3pt]{0pt}{15pt}}%
  }
  \def\@oddfoot{}
  \let\@evenfoot\@oddfoot
}

\def\ps@first{
  \let\@mkboth\@gobbletwo
  \def\@evenhead{}%
  \def\@oddhead{}%
  \def\@oddfoot{\hfill\footnotesize[\thepage]\hfill}%
  \let\@evenfoot\@oddfoot%
}

\makeatother

\begin{document}

\thispagestyle{empty}

\author{Rafa{\l} Bocian$^{\mathrm{a}}$,
    Thorsten Holm$^{\mathrm{b}}$ and
    Andrzej Skowro\'{n}ski$^{\mathrm{a},*}$\\
\footnotesize\it
$^{\mathrm{a}}$%
Faculty of Mathematics and Computer Science, Nicolaus Copernicus
University, Chopina 12/18,\\
\footnotesize\it 87-100 Toru\'{n}, Poland \\
\footnotesize\it
$^{\mathrm{b}}$%
Department of Pure Mathematics, University of Leeds,
Leeds LS2 9JT, United Kingdom
}
\title{\textrm{Derived equivalence classification of
one-parametric selfinjective algebras}}

\date{}

\maketitle

\renewcommand{\thefootnote}{}
\noindent\footnote{\hspace{-6.5mm}$^*$Corresponding author.\\
 \textit{E-mail addresses:} rafalb@mat.uni.torun.pl (R. Bocian),
              tholm@maths.leeds.ac.uk (T.~Holm),
              skowron@mat.uni.torun.pl (A. Skowro\'{n}ski).}

\hrule

\begin{center}
{\em Dedicated to Claus Michael Ringel on the occasion of his
sixtieth birthday}
\end{center}

\hrule

\bigskip

 \begin{footnotesize}
\noindent\textbf{Abstract}
\bigskip

In continuation of our paper \cite{bib:BoHoSk} we complete the
description of the derived equivalence normal forms of all
one-parametric selfinjective algebras over algebraically closed
fields which admit simply connected Galois coverings. As a
consequence, a description of the stable equivalence normal forms
of these algebras is obtained.

\bigskip


\end{footnotesize}

\bigskip

\hrule

\bigskip\bigskip

\setcounter{section}{-1}

\section{\hspace{-2.5ex}. Introduction and the main result.}

Throughout the paper $K$ will denote a fixed algebraically closed
field. By an algebra we mean a finite dimensional $K$-algebra
(associative, with an identity), which we shall assume (without
loss of generality) to be basic and connected. An algebra $A$ can
be written as a bound quiver algebra $A\cong KQ/I$, where
$Q=Q_{A}$ is the Gabriel quiver of $A$ and $I$ is an admissible
ideal in the path algebra $KQ$ of $Q$. An algebra $A$ is called
selfinjective if the projective $A$-modules are injective.

From Drozd's Tame and Wild Theorem \cite{bib:Dr} the class of
algebras may be divided into two disjoint classes. One class
consists of the tame algebras for which the indecomposable modules
occur, in each dimension $d$, in a finite number of discrete and a
finite number of one-parametric families. The second class is
formed by the wild algebras whose representation theory comprises
the representation theories of all finite dimensional
$K$-algebras. Accordingly we may realistically hope to classify
the indecomposable finite dimensional modules only for the tame
algebras. A special class of tame algebras is formed by the
algebras of finite representation type having only finitely many
isomorphism classes of indecomposable finite dimensional modules.
The representation theory of algebras of finite representation
type is presently well understood, and in particular the Morita
equivalence classes of all selfinjective algebras of finite
representation type are classified (see
\cite{bib:BrLaRi},\cite{bib:HuWa}, \cite{bib:Rie2},
\cite{bib:Rie3}, \cite{bib:Was}). The representation theory of
arbitrary tame algebras is still only emerging.

For an algebra $A$, we denote by $\md A$ the category of finite
dimensional left $A$-modules and by $\ulmod A$ its stable
category. Recall that the objects of $\ulmod A$ are the objects of
$\md A$ without nonzero projective direct summands, and for any
two objects $M$ and $N$ in $\ulmod A$ the $K$-module
$\ulhom_{A}\lo M,N\po$ of morphisms from $M$ to $N$ is the
quotient $\Hom_{A}\lo M,N\po/P\lo M,N\po$, where $P\lo M, N\po$ is
the submodule of $\Hom_{A}\lo M,N\po$ consisting of all
$A$-homomorphisms which factorize through a projective $A$-module.
Moreover, we denote by $D^{b}\lo\md A\po$ the \textit{derived
category} of bounded complexes of modules from $\md A$
\cite{bib:Hap2}. Two algebras $A$ and $B$ are said to be
\textit{stably equivalent} if the stable module categories $\ulmod
A$ and $\ulmod B$ are equivalent. Finally, two algebras $A$ and
$B$ are said to be \textit{derived equivalent} if the derived
categories $D^{b}\lo\md A\po$ and $D^{b}\lo\md B\po$ are
equivalent as triangulated categories. Since Happel's work
\cite{bib:Hap1} interpreting tilting theory in terms of
equivalences of derived categories, the machinery of derived
categories has been of interest to representation theorists. In
\cite{bib:Ric1} J. Rickard proved his celebrated criterion: two
algebras $A$ and $B$ are derived equivalent if and only if $B$ is
the endomorphism algebra of a tilting complex over $A$. Since a
lot of interesting properties are preserved by derived
equivalences of algebras, it is for many purposes important to
classify classes of algebras up to derived equivalence, instead of
Morita equivalence. For instance, for selfinjective algebras the
representation type is an invariant of the derived category. In
fact, derived equivalent selfinjective algebras are stably
equivalent \cite{bib:Ric2}, and for any algebra a stable
equivalence preserves the representation type \cite{bib:Kra},
\cite{bib:KrZw}. Finally, we note that derived equivalent algebras
have the same number of pairwise nonisomorphic simple modules. The
derived equivalence classification has been established for some
distinguished classes of tame selfinjective algebras (see
\cite{bib:Asa}, \cite{bib:BiHoSk}, \cite{bib:BoHoSk},
\cite{bib:Hol}, \cite{bib:Kau}, \cite{bib:Mem}, \cite{bib:Ric2}
for some work in this direction).

We are concerned with the problem of the derived equivalence
classification of one-parametric selfinjective algebras. Recall
that an algebra $A$ of infinite representation type is called
\textit{one-parametric} if there exists a $A$-$K\lkw
x\pkw$-bimodule $M$ which is finitely generated and free as right
$K\lkw x\pkw$-module and, for any dimension $d$, all but a finite
number of isomorphism classes of indecomposable (left) $A$-modules
of dimensional $d$ are of the form $M \otimes K\lkw x\pkw/\lo
x-\lambda\po^{m}$ for some $\lambda \in K$ and some $m\geqslant
1$. We also mention that the class of one-parametric algebras
coincides with the class of algebras having exactly one generic
module \cite{bib:CB2}. By general theory, the class of
one-parametric selfinjective algebras splits into two classes: the
\textit{standard algebras}, having simply connected Galois
coverings, and the remaining \textit{nonstandard algebras}. It is
expected that the nonstandard one-parametric (even the
representation-infinite domestic) selfinjective algebras occur
only in characteristic $2$ and are geometric and socle
deformations of standard one-parametric selfinjective algebras.

By general theory \cite{bib:Sko}, \cite{bib:LeSk}, the class of
standard one-parametric selfinjective algebras coincides with the
class of selfinjective algebras of the form $\widehat{B}/G$, where
$\widehat{B}$ is the repetitive algebra of a tilted algebra $B$ of
Euclidean type $\widetilde{\mathbbm{A}}_{m}$ or
$\widetilde{\mathbbm{D}}_{n}$, and $G$ is an infinite cyclic group
generated by an automorphism $\varphi$ of $\widehat{B}$ whose
square $\varphi^{2}$ is a twist $\rho\nu_{\widehat{B}}$ of the
Nakayama automorphism $\nu_{\widehat{B}}$ by a rigid automorphism
$\rho$ of $\widehat{B}$. We also note that the stable
Auslander-Reiten quiver of a standard one-parametric selfinjective
algebra consists of one Euclidean component $\mathbbm{Z}\Delta$ of
type $\Delta\in\lk \widetilde{\mathbbm{A}}_{m},
\widetilde{\mathbbm{D}}_{n}\pk$ and a $\mathbbm{P}_{1}\lo
K\po$-family of stable tubes. In \cite{bib:BoSk1},
\cite{bib:BoSk2}, \cite{bib:BoSk3} we classified the Morita
equivalence classes of all standard one-parametric selfinjective
algebras by algebras arising from Brauer graphs. Moreover, the
derived (respectively, stable) equivalence normal forms of all
standard one-parametric weakly symmetric algebras has been
established in our paper \cite{bib:BoHoSk}. Recall that a
selfinjective algebra $A$ is called weakly symmetric if the socle
$\soc P$ of any indecomposable projective $A$-module $P$ is
isomorphic to its top $P/\rad P$.

The main aim of the paper is to establish the derived
(respectively, stable) equivalence normal forms of all standard
one-parametric algebras which are not weakly symmetric, using the
Morita equivalence classification of this class of algebras
established in \cite{bib:BoSk3}. For $\theta\in K\backslash\lk
0\pk$, we associate to a Brauer graph $T_{p,q,k}$ of the form
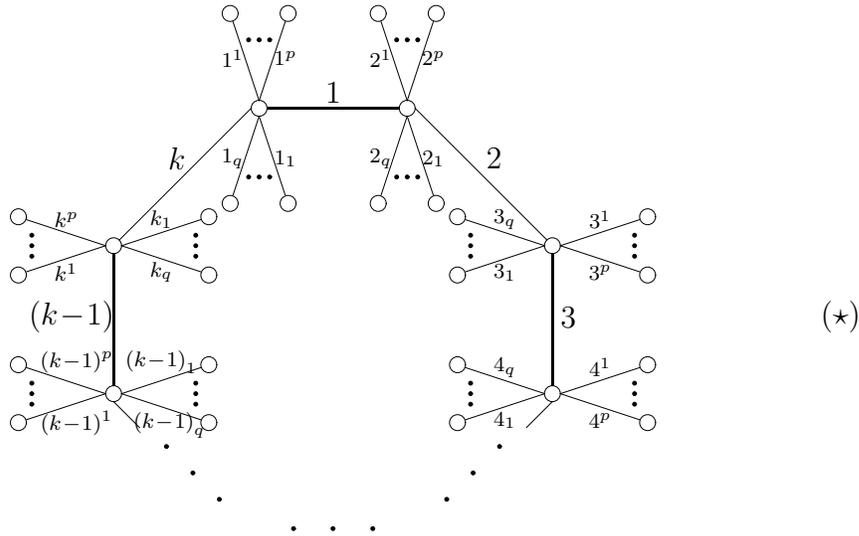
\begin{center}
\unitlength=1.000000pt
\begin{picture}(396.00,220.00)(-128.00,-115.00)
\put(81.00,59.00){\circle{6.00}} \put(-30.00,7.00){\circle{6.00}}
\put(136.00,7.00){\circle{6.00}}
\put(-30.00,-49.00){\circle{6.00}}
\put(136.00,-49.00){\circle{6.00}}
\put(25.00,59.00){\circle{6.00}}
\put(78.00,59.00){\line(-1,0){50.00}}
\put(84.00,59.00){\line(1,-1){50.00}}
\put(22.00,59.00){\line(-1,-1){50.00}}
\put(-30.00,-52.00){\line(1,-1){10.00}}
\put(136.00,-52.00){\line(-1,-1){10.00}}
\multiput(-10,-69)(10,-10){3}{\circle*{1.5}}
\multiput(116,-69)(-10,-10){3}{\circle*{1.5}}
\multiput(38,-100)(15,0){3}{\circle*{1.5}}
\put(-30.00,4.00){\line(0,-1){50.00}}
\put(136.00,4.00){\line(0,-1){50.00}}

\put(81,62){\line(-1,3){10}} \put(81,62){\line(1,3){10}}
\put(70,95){\circle{6.00}} \put(92,95){\circle{6.00}}
\multiput(77,85)(4,0){3}{\circle*{1.5}}
\put(81,56){\line(-1,-3){10}} \put(81,56){\line(1,-3){10}}
\put(70,23){\circle{6.00}} \put(92,23){\circle{6.00}}
\multiput(77,33)(4,0){3}{\circle*{1.5}}

\put(25,62){\line(-1,3){10}} \put(25,62){\line(1,3){10}}
\put(14,95){\circle{6.00}} \put(36,95){\circle{6.00}}
\multiput(21,85)(4,0){3}{\circle*{1.5}}
\put(25,56){\line(-1,-3){10}} \put(25,56){\line(1,-3){10}}
\put(14,23){\circle{6.00}} \put(36,23){\circle{6.00}}
\multiput(21,33)(4,0){3}{\circle*{1.5}}

\put(-33,7){\line(-3,1){30}} \put(-33,7){\line(-3,-1){30}}
\put(-66,18){\circle{6.00}} \put(-66,-4){\circle{6.00}}
\multiput(-61,3)(0,4){3}{\circle*{1.5}}
\put(-27,7){\line(3,1){30}} \put(-27,7){\line(3,-1){30}}
\put(6,18){\circle{6.00}} \put(6,-4){\circle{6.00}}
\multiput(1,3)(0,4){3}{\circle*{1.5}}

\put(-33,-49){\line(-3,1){30}} \put(-33,-49){\line(-3,-1){30}}
\put(-66,-38){\circle{6.00}} \put(-66,-60){\circle{6.00}}
\multiput(-61,-53)(0,4){3}{\circle*{1.5}}
\put(-27,-49){\line(3,1){30}} \put(-27,-49){\line(3,-1){30}}
\put(6,-38){\circle{6.00}} \put(6,-60){\circle{6.00}}
\multiput(1,-53)(0,4){3}{\circle*{1.5}}

\put(133,7){\line(-3,1){30}} \put(133,7){\line(-3,-1){30}}
\put(172,18){\circle{6.00}} \put(172,-4){\circle{6.00}}
\multiput(167,3)(0,4){3}{\circle*{1.5}}
\put(139,7){\line(3,1){30}} \put(139,7){\line(3,-1){30}}
\put(100,18){\circle{6.00}} \put(100,-4){\circle{6.00}}
\multiput(105,3)(0,4){3}{\circle*{1.5}}

\put(133,-49){\line(-3,1){30}} \put(133,-49){\line(-3,-1){30}}
\put(172,-38){\circle{6.00}} \put(172,-60){\circle{6.00}}
\multiput(167,-53)(0,4){3}{\circle*{1.5}}
\put(139,-49){\line(3,1){30}} \put(139,-49){\line(3,-1){30}}
\put(100,-38){\circle{6.00}} \put(100,-60){\circle{6.00}}
\multiput(105,-53)(0,4){3}{\circle*{1.5}}

\put(15.00,78.00){\makebox(0,0){\scriptsize $1^{1}$}}
\put(35.00,78.00){\makebox(0,0){\scriptsize $1^{p}$}}
\put(15.00,40.00){\makebox(0,0){\scriptsize $1_{q}$}}
\put(35.00,40.00){\makebox(0,0){\scriptsize $1_{1}$}}

\put(71.00,78.00){\makebox(0,0){\scriptsize $2^{1}$}}
\put(91.00,78.00){\makebox(0,0){\scriptsize $2^{p}$}}
\put(71.00,40.00){\makebox(0,0){\scriptsize $2_{q}$}}
\put(91.00,40.00){\makebox(0,0){\scriptsize $2_{1}$}}

\put(154.00,17.00){\makebox(0,0){\scriptsize $3^{1}$}}
\put(154.00,-3.00){\makebox(0,0){\scriptsize $3^{p}$}}
\put(118.00,17.00){\makebox(0,0){\scriptsize $3_{q}$}}
\put(118.00,-3.00){\makebox(0,0){\scriptsize $3_{1}$}}

\put(154.00,-39.00){\makebox(0,0){\scriptsize $4^{1}$}}
\put(154.00,-59.00){\makebox(0,0){\scriptsize $4^{p}$}}
\put(118.00,-39.00){\makebox(0,0){\scriptsize $4_{q}$}}
\put(118.00,-59.00){\makebox(0,0){\scriptsize $4_{1}$}}

\put(-48.00,17.00){\makebox(0,0){\scriptsize $k^{p}$}}
\put(-48.00,-3.00){\makebox(0,0){\scriptsize $k^{1}$}}
\put(-12.00,17.00){\makebox(0,0){\scriptsize $k_{1}$}}
\put(-12.00,-3.00){\makebox(0,0){\scriptsize $k_{q}$}}

\put(-44.00,-37.00){\makebox(0,0){\scriptsize $\lo
k\!-\!1\po^{p}$}} \put(-44.00,-60.00){\makebox(0,0){\scriptsize
$\lo k\!-\!1\po^{1}$}}
\put(-12.00,-37.00){\makebox(0,0){\scriptsize $\lo
k\!-\!1\po_{1}$}} \put(-9.00,-61.00){\makebox(0,0){\scriptsize
$\lo k\!-\!1\po_{q}$}}

\put(53.00,65.00){\makebox(0,0){$1$}}
\put(114.00,40.00){\makebox(0,0){$2$}}
\put(142.00,-20.00){\makebox(0,0){$3$}}
\put(-46.00,-20.00){\makebox(0,0){$\lo k\!-\!1\po$}}
\put(-6.00,40.00){\makebox(0,0){$k$}}

\put(245.00,-20,00){\makebox(0,0){$\lo\star\po$}}
\end{picture}
\end{center}
with $p,q\geqslant 0$ and $k\geqslant 2$, and a positive integer
$s$ with $1\leqslant s\leqslant k\!-\!1$, $\gcd\lo
s\!+\!2,k\po=1$, $\gcd\lo s,k\po=1$, a standard one-parametric
selfinjective algebra $\Lambda\lo p,q,k,s,\theta\po$ of Euclidean
type $\widetilde{\mathbbm{A}}_{m}$ (with $m=2\lo p+q+1\po k-1$).
Furthermore, to a Brauer tree $T_{n}$ of the form
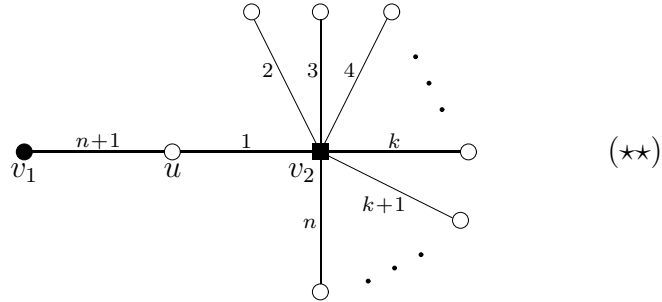
\begin{center}
\unitlength=1.000000pt
\begin{picture}(250.00,120.00)(-40.00,0.00)
\put(78.00,56.00){\rule{6.00pt}{6.00pt}}
\put(137.00,59.00){\circle{6.00}} \put(25.00,59.00){\circle{6.00}}
\put(81.00,6.00){\circle{6.00}} \put(81.00,112.00){\circle{6.00}}
\put(55.00,112.00){\circle{6.00}}
\put(134.00,33.00){\circle{6.00}}
\put(108.00,112.00){\circle{6.00}}
\put(-31.00,59.00){\circle*{6.00}}
\put(78.00,59.00){\line(-1,0){50.00}}
\put(81.00,59.00){\line(0,-1){50.00}}
\put(81.00,59.00){\line(0,1){50.00}}
\put(81.00,59.00){\line(1,2){25.30}}
\put(81.00,59.00){\line(-1,2){25.00}}
\put(81.00,59.00){\line(2,-1){50.00}}
\put(84.00,59.00){\line(1,0){50.00}}
\put(22.00,59.00){\line(-1,0){50.00}}
\multiput(117.00,95.00)(5,-10){3}{\circle*{1.5}}
\multiput(99.00,10.00)(10,5){3}{\circle*{1.5}}
\put(25.00,52.00){\makebox(0,0){\normalsize $u$}}
\put(-31.00,51.00){\makebox(0,0){\normalsize $v_{1}$}}
\put(74.00,51.00){\makebox(0,0){\normalsize $v_{2}$}}
\put(-3.00,63.00){\makebox(0,0){\scriptsize $n\!+\!1$}}
\put(53.00,63.00){\makebox(0,0){\scriptsize $1$}}
\put(109.00,63.00){\makebox(0,0){\scriptsize $k$}}
\put(77.00,32.00){\makebox(0,0){\scriptsize $n$}}
\put(105.00,39.00){\makebox(0,0){\scriptsize $k\!+\!1$}}
\put(61.00,90.00){\makebox(0,0){\scriptsize $2$}}
\put(78.00,90.00){\makebox(0,0){\scriptsize $3$}}
\put(92.00,90.00){\makebox(0,0){\scriptsize $4$}}
\put(200.00,59.00){\makebox(0,0){\normalsize $\lo\star\star\po$}}
\end{picture}
\end{center}
we associate a standard one-parametric selfinjective algebra
$\Gamma^{*}\!\lo n\po$ of Euclidean type
$\widetilde{\mathbbm{D}}_{n}$. See Section \ref{morita} below for
the precise definitions.

The following theorem is the main result of the paper.

\begin{tw1} \label{mainthm}
For a standard selfinjective algebra $A$ the following statements
are equivalent:
\begin{arablist}
\item $A$ is one-parametric but not weakly symmetric.
\item $A$ is derived equivalent to an algebra of the form
$\Lambda\lo p,q,k,s,\theta\po$ or $\Gamma^{*}\lo n\po$.
\item $A$ is stably equivalent to an algebra of the form
$\Lambda\lo p,q,k,s,\theta\po$ or $\Gamma^{*}\lo n\po$.
\end{arablist}
\end{tw1}
\noindent We note that the algebras $\Lambda\lo p,q,k,s,\theta\po$
and $\Gamma^{*}\!\lo n\po$ are never stably equivalent, because
their stable Auslander-Reiten quivers are not isomorphic (see
Propositions \ref{structard} and \ref{structar}). Moreover, by
Proposition \ref{structard}, the algebras $\Gamma^{*}\!\lo m\po$
and $\Gamma^{*}\!\lo n\po$ are stably equivalent only for $m=n$.
Finally, we refer to Proposition \ref{propiso} for a necessity
criterion for the stable equivalence of algebras of the form
$\Lambda\lo p,q,k,s,\theta\po$.

As an application of Theorem \ref{mainthm} and results of our
earlier papers \cite{bib:BoHoSk} and \cite{bib:BoHoSk2} we obtain
the following description of Auslander's representation
dimension $\rdim\lo A\po$ and Rouquier's dimension $\dim
\ulmod A$ of the triangulated category $\ulmod A$ (see Section $8$
for the relevant definitions) of arbitrary standard one-parametric
selfinjective algebras $A$.

\begin{tw1} \label{mainthm2}
Let $A$ be a standard one-parametric selfinjective algebra. Then
$\rdim\lo A\po=3$ and $\dim \ulmod A=1$.
\end{tw1}

The paper is organized as follows. In Section \ref{morita} we
recall the Morita equivalence classification of all standard
one-parametric selfinjective but not weakly symmetric algebras.
Sections \ref{tilting}, \ref{moving}, \ref{reduction} and
\ref{derived} are devoted to the proof of the implication $\lo
1\po\Longrightarrow\lo 2\po$ of Theorem \ref{mainthm} in the
Euclidean case $\widetilde{\mathbbm{A}}_{m}$, and Section
\ref{dnderived} to the proof of this implication
in the Euclidean case $\widetilde{\mathbbm{D}}_{n}$. In Section
\ref{ARseq} we describe the stable Auslander-Reiten quivers of
algebras occurring in Theorem \ref{mainthm}. Section \ref{wsalg}
is devoted to a class of weakly symmetric algebras needed in the
proof of implication $\lo 3\po\Longrightarrow\lo 1\po$ of Theorem
\ref{mainthm}. In Section \ref{proofthm1} we complete the proof of
Theorem \ref{mainthm}. In Section \ref{stableequiv} we discuss the
stable equivalence of algebras occurring in Theorem \ref{mainthm}.
The final Section \ref{secrepdim} is devoted to the proof of
Theorem \ref{mainthm2}.

For basic background on the representation theory of algebras we
refer to \cite{bib:AuReSm}, \cite{bib:Rin}, and for background
on selfinjective algebras to \cite{bib:Erd}, \cite{bib:Yam}.


\section{\hspace{-2.5ex}. One-parametric nonweakly symmetric selfinjective algebras}
\label{morita}


We describe in this section the Morita equivalence classification
of all standard one-parametric but not weakly symmetric
selfinjective algebras established in \cite{bib:BoSk3}.

A \textit{Brauer graph} $T$ is a finite connected undirected
graph, where for each vertex there is a fixed circular order on
the edges adjacent to it (see \cite{bib:Al}, \cite{bib:Kau},
\cite{bib:Mem}, \cite{bib:Rog}). In our context we assume that $T$
has at most one cycle (which may be or may not be a loop). We draw
$T$ in a plane and agree that the edges adjacent to a given vertex
are clockwise ordered. Given a Brauer graph $T$, this defines a
\textit{Brauer quiver} $Q_{T}$ as follows. The vertices of $Q_{T}$
are the edges of $T$ and there is an arrow $i\longrightarrow j$ in
$Q_{T}$ if and only if $j$ is the direct successor of $i$
in the order around some vertex (to which $i$ and $j$ are both
adjacent). We require that every vertex of $Q_{T}$ belongs to
exactly two cycles. Note that this implicitly means that, for
every end vertex of $T$, there is a loop in $Q_{T}$.

Let $T$ be a Brauer graph with exactly one cycle ${\cal R}_{k}$,
having $k\geqslant 2$ edges.
Let $v_{1},v_{2},\ldots,v_{k}$ be the vertices
of ${\cal R}_{k}$ and $e_{i}=\lk v_{i},v_{i+1}\pk,
i=1,2,\ldots,k$, where $v_{k+1}=v_{1}$, the edges of ${\cal
R}_{k}$. If $v$ is a vertex of the Brauer graph $T$ which is not a
vertex of the cycle ${\cal R}_{k}$ then by $n\!\lo v\po$ we denote
the edge incident to $v$ on the unique walk in $T$ from $v$ to
the cycle ${\cal R}_{k}$. Moreover, for $i=1,2,\ldots,k$, we
denote by $n\!\lo v_{i}\po$ the edge $e_{i}$. For a vertex $v$ of
the graph $T$, we denote by $l\lo {\cal R}_{k},v\po$ the distance
of $v$ to the cycle ${\cal R}_{k}$. Hence $l\lo {\cal
R}_{k},v\po=0$ if and only if $v$ belongs to ${\cal R}_{k}$. By
an \textit{automorphism of the Brauer graph $T$} we mean an
automorphism of the graph $T$ which preserves the fixed circular
order on the edges adjacent to any vertex. \textit{A rotation} of
the Brauer graph $T$ is an automorphism $\sigma$ of the Brauer
graph $T$ such that, for some integer $s$ with $1\leqslant s
\leqslant k-1$, we have $\sigma\!\lo v_{i}\po=v_{i+s}$ for all
$i=1,2,\ldots k$ (where $k+r=r$ for $r\geqslant 1$), and then we
set $\sigma=\sigma_{s}$.
For $k=2$, we set $\sigma_{1}\!\lo
e_{1}\po=e_{2}$ and $\sigma_{1}\!\lo e_{2}\po=e_{1}$.

Assume that $s$ is a positive integer such that $1\leqslant
s\leqslant k-1$ and $\gcd\!\lo s+2,k\po$=1. We shall define a
generalized Brauer quiver $Q\lo T,\sigma_{s}\po$, obtained from
the usual Brauer quiver $Q_{T}$ of the Brauer graph $T$ by
shifting some arrows of $Q_{T}$ using the rotation $\sigma_{s}$ of
$T$. By a $\sigma_{s}$-orbit of a vertex $v$ of $T$ we mean the
orbit of $v$ with respect to the action of the cyclic group
generated by $\sigma_{s}$ on the vertices of
$T$. We note that if two vertices $v$ and $w$ of $T$ belong to the
same $\sigma_{s}$-orbit then $l\lo {\cal R}_{k},v\po=l\lo {\cal
R}_{k},w\po$. Moreover, all $\sigma_{s}$-orbits of vertices of $T$
have the same number of elements, namely $k/d$, where $d=\gcd\lo
s,k\po$. For $m\geqslant 0$, denote by $V_{m}$ the set of all
vertices of $T$ with $l\lo {\cal R}_{k},v\po=m$. Observe that
$V_{m}$ is a disjoint union of $\frac{d|V_{m}|}{k}$
$\sigma_{s}$-orbits.

In order to define the generalized Brauer quiver $Q\lo
T,\sigma_{s}\po$, we introduce an order $p\lo T,\sigma_{s}\po$ of
the edges of the Brauer graph $T$, as the union of
$\sum_{m=0}^{\infty}\frac{d|V_{m}|}{k}$ cyclic orders $p\lo
T,\sigma_{s},v\po$ defined for the representatives $v$ of all
pairwise different $\sigma_{s}$-orbits of vertices of $T$. Let $v$
be a vertex of $T$. We define the cyclic order $p\lo
T,\sigma_{s},v\po$ by invoking the cyclic orders of edges around the
vertices $v,\sigma_{s}\lo v\po,\ldots,\sigma_{s}^{k/d-1}\lo v\po$
in the Brauer graph $T$. Let $r\in\lk 0,1,\ldots,k/d-1\pk$, let $i$
be an edge of $T$ adjacent to the vertex $\sigma_{s}^{r}\lo v\po$,
and let $j$ be the direct successor of $i$ in the cyclic order in $T$
around $\sigma_{s}^{r}\lo v\po$. If $j\neq n\lo \sigma_{s}^{r}\lo
v\po\po$, then $j$ is defined to be the direct successor of $i$ in
the cyclic order $p\lo T,\sigma_{s},v\po$. If $j= n\lo
\sigma_{s}^{r}\lo v\po\po$ then $n\lo \sigma_{s}^{r+1}\lo
v\po\po=\sigma_{s}\lo n\lo\sigma_{s}^{r}\lo v\po\po\po$ is said to
be the direct successor of $i$ in the cyclic order $p\lo
T,\sigma_{s},v\po$. Therefore, we replaced the cyclic orders
around the vertices $\sigma_{s}^{r}\lo v\po$, $0\leqslant
r\leqslant k/d-1$, by one (bigger) cyclic order $p\lo
T,\sigma_{s},v\po$. Observe also that if $e=\lk v,w\pk$ is an edge
of $T$ which is not on the cycle ${\cal R}_{k}$, or if $e$ is on the
cycle ${\cal R}_{k}$ and $d>1$, then $e$ belongs to exactly two
cyclic orders, namely $p\lo T,\sigma_{s},v\po$ and $p\lo
T,\sigma_{s},w\po$. On the other hand, if $e=\lk v,w\pk$ is an
edge of the cycle ${\cal R}_{k}$ and $d=1$, then $e$ occurs twice
in the cyclic order $p\lo T,\sigma_{s},v\po=p\lo
T,\sigma_{s},w\po$.

We define \textit{the generalized Brauer quiver} $Q\lo
T,\sigma_{s}\po$ as follows. The vertices of $Q\lo
T,\sigma_{s}\po$ are the edges of $T$ and there is an arrow
$i\longrightarrow j$ in $Q\lo T,\sigma_{s}\po$ if and only if $j$
is the direct successor of $i$ in the order $p\lo
T,\sigma_{s}\po$.

For $\theta\in K\backslash\lk 0\pk$, we define the algebra
$\Omega^{\lo 1\po}\!\lo T,\sigma_{s},\theta\po$ as the bound
quiver algebra $KQ(T,\sigma_{s})/\overline{I}^{\lo 1\po }\!\lo
T,\sigma_{s},\theta\po$, where $KQ(T,\sigma_{s})$ is the path
algebra of the quiver $Q(T,\sigma_{s})$ and $\overline{I}^{\lo
1\po }\!\lo T,\sigma_{s},\theta\po$ is the ideal in
$KQ(T,\sigma_{s})$ generated by the elements:
\begin{arablist}
\item $\alpha\beta$ where $\alpha=i_{1}\longrightarrow i_{2}$,
$\beta=i_{2}\longrightarrow i_{3}$ and $i_{1},i_{2},i_{3}$ are not
consecutive elements in the cyclic order $p\lo T,\sigma_{s}\po$.
\item $C\!\lo i,p\lo T,\sigma_{s},v\po\po-C\!\lo i,p\lo
T,\sigma_{s},w\po\po$, for $i\neq e_{1}$ and $C\!\lo e_{1},p\lo
T,\sigma_{s},v\po\po-\theta C\!\lo e_{1},p\lo
T,\sigma_{s},w\po\po$, where $i=\lk v,w\pk$ is an edge of $T$,
$C\!\lo i,p\lo T,\sigma_{s},v\po\po$ and $C\!\lo i,p\lo
T,\sigma_{s},w\po\po$ are the paths from $i$ to $\sigma_{s}\!\lo
i\po$ in the quiver $Q(T,\sigma_{s})$, corresponding to the
consecutive elements $i,\ldots,\sigma_{s}^{}\!\lo i\po$ of the
cyclic orders $p\lo T,\sigma_{s},v\po$ and $p\lo
T,\sigma_{s},w\po$, respectively.
\end{arablist}
If $T=T_{p,q,k}$ is of the form $\lo\star\po$ and $\gcd\lo
s,k\po=1$, then we denote $\Omega^{\lo 1\po}\!\lo
T,\sigma_{s},\theta\po$ by $\Lambda\lo p,q,k,s,\theta\po$.

\begin{exa} \rm \label{przyklad_1}
Let $s=2$ and $T_{2,1,3}$ be the following Brauer graph
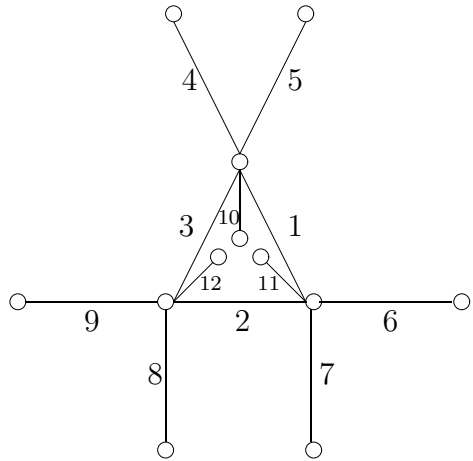
\begin{center}
\unitlength=1.000000pt
\begin{picture}(172.00,171.00)(25.00,2.00)
\put(130,143){\makebox(0,0){\normalsize $5$}}
\put(130,88){\makebox(0,0){\normalsize $1$}}
\put(90,143){\makebox(0,0){\normalsize $4$}}
\put(110,52){\makebox(0,0){\normalsize $2$}}
\put(89,88){\makebox(0,0){\normalsize $3$}}
\put(166,52){\makebox(0,0){\normalsize $6$}}
\put(142,32){\makebox(0,0){\normalsize $7$}}
\put(77,32){\makebox(0,0){\normalsize $8$}}
\put(53,52){\makebox(0,0){\normalsize $9$}}
\put(120,66){\makebox(0,0){\scriptsize $11$}}
\put(98,66){\makebox(0,0){\scriptsize $12$}}
\put(105,91){\makebox(0,0){\scriptsize $10$}}
\put(81.00,59.00){\circle{6.00}} \put(137.00,59.00){\circle{6.00}}
\put(137.00,3.00){\circle{6.00}} \put(109.00,112.00){\circle{6}}
\put(81.00,3.00){\circle{6.00}} \put(84.00,168.00){\circle{6.00}}
\put(134.00,168.00){\circle{6.00}}
\put(101.00,76.00){\circle{6.00}}
\put(117.00,76.00){\circle{6.00}}
\put(109.00,83.00){\circle{6.00}}
\put(193.00,59.00){\circle{6.00}} \put(25.00,59.00){\circle{6.00}}
\put(84.00,59.00){\line(1,2){25.00}}
\put(134.00,59.00){\line(-1,2){25.00}}
\put(81.00,56.00){\line(0,-1){50.00}}
\put(84.00,59.00){\line(1,0){50.00}}
\put(136.00,56.00){\line(0,-1){50.00}}
\put(139.00,59.00){\line(1,0){50.00}}
\put(109.00,115.00){\line(-1,2){25.00}}
\put(109.00,115.00){\line(1,2){25.00}}
\put(84.00,59.00){\line(1,1){15.00}}
\put(78.00,59.00){\line(-1,0){50.00}}
\put(109.00,109.00){\line(0,-1){23.00}}
\put(134.00,59.00){\line(-1,1){15.00}}
\end{picture}
\end{center}
Then the order $p\lo T_{2,1,3},\sigma_{2}\po$ is the union of the
four cycles:
\begin{arablist}
\item $4,5,3,12,2,8,9,2,11,1,6,7,1,10,3$,
\item $4,8,6$,
\item $5,9,7$,
\item $10,12,11$.
\end{arablist}
Then the generalized Brauer quiver $Q\lo T_{2,1,3},\sigma_{2}\po$
of the graph $T_{2,1,3}$ is of the form
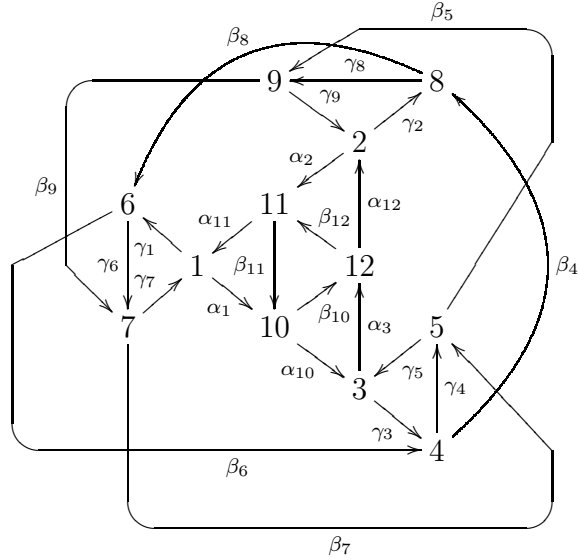
\begin{displaymath}
\xymatrix @C=1.2pc@R=0.8pc{
& & & & & & & & \\
& & & & 9 \ar[rd]^{\gamma_{9}} \ar `l^d[lll]`d^[dddlll]_{\beta_{9}}
[lldddd] & & 8 \ar[ll]_{\gamma_{8}} \ar@/_3.0pc/[lllldd]_{\beta_{8}}
& & & \\
& & & & & 2 \ar[dl]_{\alpha_{2}} \ar[ur]_{\gamma_{2}} & & & \\
& & 6 \ar[dd]_{\gamma_{6}} \ar `l^l[lld]`d^r[lldddd]
[rrrrdddd]_{\beta_{6}} & & 11 \ar[dl]_{\alpha_{11}} \ar[dd]_{\beta_{11}}
& & \\
& & & 1 \ar[dr]_{\alpha_{1}} \ar[lu]^{\gamma_{1}} & & 12
\ar[uu]_{\alpha_{12}} \ar[ul]_{\beta_{12}} & &\\
& & 7 \ar[ru]^{\gamma_{7}} \ar
`d^d[ddd]`d^r[rrrrrddd]`r^u[rrrrrr]_{\beta_{7}}`u^[rrrrrrdd] [rrrr]
& & 10 \ar[dr]_{\alpha_{10}} \ar[ur]_{\beta_{10}} &
& 5 \ar[dl]^{\gamma_{5}} \ar
`r^r[rruuu]`u^[ru]`u^l[llluuuuu]`d^[luuuuu]_{\beta_{5}}[lluuuu] & & \\
& & & & & 3 \ar[uu]_{\alpha_{3}} \ar[dr]_{\gamma_{3}} & & & & & \\
& & & & & & 4 \ar[uu]_{\gamma_{4}} \ar@/_3.5pc/[uuuuuu]_{\beta_{4}} & & \\
& & & & & & & & & & & \\
}
\end{displaymath}
and $\Lambda\lo 2,1,3,2,\theta\po$ is given by the above quiver
and the ideal $\overline{I}^{\lo 1\po }\!\lo
T_{2,1,3},\sigma_{2},\theta\po$ in $KQ\lo T_{2,1,3},\sigma_{2}\po$
generated by the elements: $\alpha_{1}\beta_{10}$,
$\alpha_{2}\beta_{11}$, $\alpha_{3}\beta_{12}$,
$\beta_{10}\alpha_{12}$, $\beta_{11}\alpha_{10}$,
$\beta_{12}\alpha_{11}$, $\gamma_{1}\beta_{6}$,
$\beta_{6}\gamma_{4}$, $\gamma_{3}\beta_{4}$,
$\beta_{4}\gamma_{8}$, $\gamma_{2}\beta_{8}$,
$\beta_{8}\gamma_{6}$, $\gamma_{4}\beta_{5}$,
$\beta_{5}\gamma_{9}$, $\gamma_{6}\beta_{7}$,
$\beta_{7}\gamma_{5}$, $\gamma_{8}\beta_{9}$,
$\beta_{9}\gamma_{7}$, $\alpha_{11}\alpha_{1}$,
$\alpha_{10}\alpha_{3}$, $\alpha_{12}\alpha_{2}$,
$\gamma_{9}\gamma_{2}$, $\gamma_{5}\gamma_{3}$,
$\gamma_{7}\gamma_{1}$,
$\alpha_{1}\alpha_{10}\gamma_{3}\gamma_{4}\gamma_{5}-
\theta\gamma_{1}\gamma_{6}\gamma_{7}\alpha_{1}\alpha_{10}$,
$\alpha_{2}\alpha_{11}\gamma_{1}\gamma_{6}\gamma_{7}-
\gamma_{2}\gamma_{8}\gamma_{9}\alpha_{2}\alpha_{11}$,
$\alpha_{3}\alpha_{12}\gamma_{2}\gamma_{8}\gamma_{9}-
\gamma_{3}\gamma_{4}\gamma_{5}\alpha_{3}\alpha_{12}$,
$\alpha_{10}\gamma_{3}\gamma_{4}\gamma_{5}\alpha_{3}-\beta_{10}$,
$\alpha_{11}\gamma_{1}\gamma_{6}\gamma_{7}\alpha_{1}-\beta_{11}$,
$\alpha_{12}\gamma_{2}\gamma_{8}\gamma_{9}\alpha_{2}-\beta_{12}$,
$\gamma_{4}\gamma_{5}\alpha_{3}\alpha_{12}\gamma_{2}-\beta_{4}$,
$\gamma_{5}\alpha_{3}\alpha_{12}\gamma_{2}\gamma_{8}-\beta_{5}$,
$\gamma_{6}\gamma_{7}\alpha_{1}\alpha_{10}\gamma_{3}-\beta_{6}$,
$\gamma_{7}\alpha_{1}\alpha_{10}\gamma_{3}\gamma_{4}-\beta_{7}$,
$\gamma_{8}\gamma_{9}\alpha_{2}\alpha_{11}\gamma_{1}-\beta_{8}$,
$\gamma_{9}\alpha_{2}\alpha_{11}\gamma_{1}\gamma_{6}-\beta_{9}$.
\end{exa}

Let $T$ be a Brauer tree. Then the simple cycles of the Brauer
quiver $Q_{T}$ may be divided into two camps, $\alpha$-camps and
$\beta$-camps, in such a way that any two cycles which intersect
nontrivially belong to different camps. We denote by $\alpha_{i}$
(respectively, $\beta_{i}$) the arrow of the $\alpha$-camp
(respectively, $\beta$-camp) of $Q_{T}$ starting at a vertex $i$,
and by $\alpha\lo i\po$ (respectively, $\beta\lo i\po$) the end
vertex of $\alpha_{i}$ (respectively, $\beta_{i}$). We also denote
by $A_{i}$ (respectively, $B_{i}$) the cycle from $i$ to $i$ going
once around the $\alpha$-cycle (respectively, $\beta$-cycle)
through $i$.

Let $T$ be a Brauer tree with two (different) distinguished
vertices $v_{1}$ and $v_{2}$ such that $v_{1}$ is the end of
exactly one edge $a$. Let the edge $b$ be the direct successor of
the edge $a$ and let $c$ be the direct predecessor of the edge $a$ in
the cyclic order of edges at the end vertex $u$ of $a$ different
from $v_{1}$. The vertices $v_{1}$, $v_{2}$ and edges $b$, $c$
determine a subtree
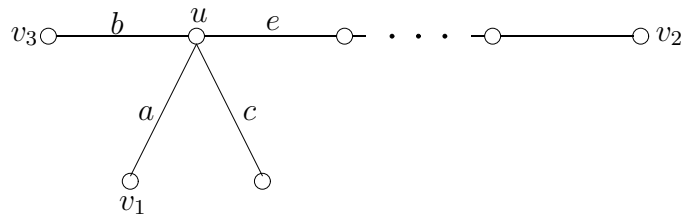
\begin{center}
\unitlength=1.000000pt
\begin{picture}(271.00,80.00)(7.00,-10.00)
\put(110,64){\makebox(0,0){\normalsize $e$}}
\put(51,64){\makebox(0,0){\normalsize $b$}}
\put(62,30){\makebox(0,0){\normalsize $a$}}
\put(101,30){\makebox(0,0){\normalsize $c$}}
\put(57,-5){\makebox(0,0){\normalsize $v_{1}$}}
\put(82,67){\makebox(0,0){\normalsize $u$}}
\put(260,59){\makebox(0,0){\normalsize $v_{2}$}}
\put(16,59){\makebox(0,0){\normalsize $v_{3}$}}
\put(81.00,59.00){\circle{6.00}} \put(137.00,59.00){\circle{6.00}}
\put(193.00,59.00){\circle{6.00}}
\put(249.00,59.00){\circle{6.00}} \put(56.00,4.00){\circle{6.00}}
\put(25.00,59.00){\circle{6.00}} \put(106.00,4.00){\circle{6.00}}
\put(78.00,59.00){\line(-1,0){50.00}}
\put(81.00,56.00){\line(-1,-2){25.00}}
\put(81.00,56.00){\line(1,-2){25.00}}
\put(84.00,59.00){\line(1,0){50.00}}
\put(140.00,59.00){\line(1,0){5.00}}
\put(185.00,59.00){\line(1,0){5.00}}
\multiput(155,59)(10,0){3}{\circle*{1.5}}
\put(196.00,59.00){\line(1,0){50.00}}
\end{picture}
\end{center}
of the Brauer tree $T$, where possibly $u=v_{2}$, $v_{2}=v_{3}$,
$b=e$, $c=e$, $b=c=e$, but every time $a\neq b$ and $a\neq c$. We
assume that the Brauer quiver $Q_{T}$ has exactly one exceptional
cycle (with multiplicity two) given by the edges of $T$ converging
at the exceptional vertex $v_{2}$. Moreover, we assume that the
cycle in $Q_{T}$ corresponding to the vertex $u$ is an
$\alpha$-cycle.

We define the algebra $\Omega^{\lo 2\po }\!\lo T,v_{1},v_{2}\po$
as the bound quiver algebra $K\overline{Q}^{\lo
2\po}_{T}/\overline{I}^{\lo 2\po}\!\lo T,v_{1},v_{2}\po$, where
$K\overline{Q}^{\lo 2\po}_{T}$ is the path algebra of the quiver
\begin{displaymath}
\overline{Q}^{\lo 2\po}_{T}=\lo\lo Q_{T}\po_{0}\!\cup\! \lk
w\pk,\lo Q_{T}\po_{1}\!\cup \!\lk \gamma_{1}:c\longrightarrow
w,\gamma_{2}:w\longrightarrow
b\pk\!\backslash\!\lk\beta_{a}:a\longrightarrow a\pk\po
\end{displaymath}
and $\overline{I}^{\lo 2\po}\!\lo T,v_{1},v_{2}\po$ is the ideal
in $K\overline{Q}^{\lo 2\po} _{T}$ generated by the elements:
\begin{arablist}
\item $\alpha_{i}\beta_{\alpha\lo i\po}$, for all vertices $i$ of
$Q_{T}$ different from $c$,
\item $\beta_{i}\alpha_{\beta\lo i\po}$, for all vertices $i$ of
$Q_{T}$ different from $a$,
\item $A_{j}-B_{j}$, if the both $\alpha$-cycle and $\beta$-cycle
through the vertex $j$ are not exceptional,
\item $A_{j}^{2}-B_{j}$, if the $\alpha$-cycle through the vertex
$j$ is exceptional but the $\beta$-cycle through $j$ is not
exceptional,
\item $A_{j}-B_{j}^{2}$, if the $\alpha$-cycle through the vertex
$j$ is not exceptional but the $\beta$-cycle through the vertex
$j$ is exceptional,
\item $\gamma_{2}\beta_{b}$,
$\beta_{\beta^{-1}\lo c\po}\gamma_{1}$,
\item $\gamma_{2}\alpha_{b}\ldots \alpha_{\alpha^{-1}\lo c\po}\gamma_{1}$,
$A_{a}$ ($\gamma_{2}\gamma_{1}$, $A_{a}$, if $b=c=e$), if the
$\alpha$-cycle through the vertex $a$ is not exceptional,
\item $\gamma_{2}A_{b}\alpha_{b}\ldots\alpha_{\alpha^{-1}\lo c\po}
\gamma_{1}$,
$A_{a}^{2}$ ($\gamma_{2}A_{b}\gamma_{1}$, $A_{a}^{2}$, if
$b=c=e$), if the $\alpha$-cycle through the vertex $a$ is
exceptional,
\item $\alpha_{c}\alpha_{a}-\gamma_{1}\gamma_{2}$.
\end{arablist}
If $\overline{Q}^{\lo 2\po}_{T}$ is of the form $\lo\star\star\po$
then we denote $\Omega^{\lo 2\po }\!\lo T,v_{1},v_{2}\po$ by
$\Gamma^{*}\lo n\po$.
\begin{exa}
Let $T_{6}$ be the following Brauer tree
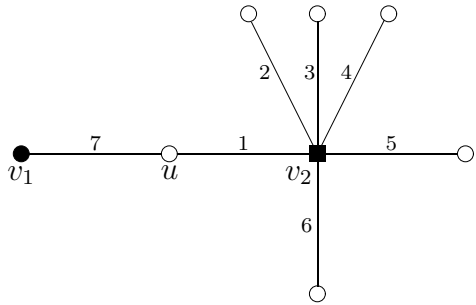
\begin{center}
\unitlength=1.000000pt
\begin{picture}(250.00,120.00)(-40.00,0.00)
\put(78.00,56.00){\rule{6.00pt}{6.00pt}}
\put(137.00,59.00){\circle{6.00}} \put(25.00,59.00){\circle{6.00}}
\put(81.00,6.00){\circle{6.00}} \put(81.00,112.00){\circle{6.00}}
\put(55.00,112.00){\circle{6.00}}
\put(108.00,112.00){\circle{6.00}}
\put(-31.00,59.00){\circle*{6.00}}
\put(78.00,59.00){\line(-1,0){50.00}}
\put(81.00,59.00){\line(0,-1){50.00}}
\put(81.00,59.00){\line(0,1){50.00}}
\put(81.00,59.00){\line(1,2){25.30}}
\put(81.00,59.00){\line(-1,2){25.00}}
\put(84.00,59.00){\line(1,0){50.00}}
\put(22.00,59.00){\line(-1,0){50.00}}
\put(25.00,52.00){\makebox(0,0){\normalsize $u$}}
\put(-31.00,51.00){\makebox(0,0){\normalsize $v_{1}$}}
\put(74.00,51.00){\makebox(0,0){\normalsize $v_{2}$}}
\put(-3.00,63.00){\makebox(0,0){\scriptsize $7$}}
\put(53.00,63.00){\makebox(0,0){\scriptsize $1$}}
\put(109.00,63.00){\makebox(0,0){\scriptsize $5$}}
\put(77.00,32.00){\makebox(0,0){\scriptsize $6$}}
\put(61.00,90.00){\makebox(0,0){\scriptsize $2$}}
\put(78.00,90.00){\makebox(0,0){\scriptsize $3$}}
\put(92.00,90.00){\makebox(0,0){\scriptsize $4$}}
\end{picture}
\end{center}
Then the Brauer quiver $\overline{Q}_{T_{6}}^{\lo 2\po}$ of the
tree $T_{6}$ is of the form
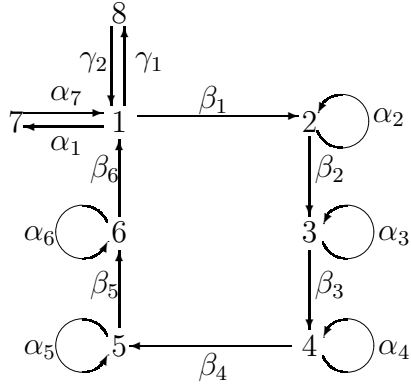
\begin{center}
\unitlength=1.0pt
\begin{picture}(166,161)(-85,-145)
 \put(48,-35){\circlearrowr}
 \put(36,-35){\makebox(0,0){\normalsize $2$}}
 \put(-6,-30){\normalsize $\beta_{1}$}
 \put(-72,-80){\normalsize $\alpha_{6}$}
 \put(-72,-123){\normalsize $\alpha_{5}$}
 \put(-6,-131){\normalsize $\beta_{4}$}
 \put(38,-100){\normalsize $\beta_{3}$}
 \put(38,-56){\normalsize $\beta_{2}$}
 \put(-47,-100){\normalsize $\beta_{5}$}
 \put(-47,-56){\normalsize $\beta_{6}$}
 \put(-30,-14){\normalsize $\gamma_{1}$}
 \put(-51,-14){\normalsize $\gamma_{2}$}
 \put(-36,-35){\makebox(0,0){\normalsize $1$}}
 \put(-36,6){\makebox(0,0){\normalsize $8$}}
 \put(-75,-35){\makebox(0,0){\normalsize $7$}}
 \put(-29,-33){\vector(1,0){61}}
 \put(-72,-32){\vector(1,0){30}}
 \put(-42,-37){\vector(-1,0){30}}
 \put(-34,-29){\vector(0,1){30}}
 \put(-39,1){\vector(0,-1){30}}
 \put(-36,-113){\vector(0,1){29}}
 \put(-36,-71){\vector(0,1){29}}
 \put(-36,-120){\makebox(0,0){\normalsize $5$}}
 \put(-36,-77){\makebox(0,0){\normalsize $6$}}
 \put(29,-120){\vector(-1,0){61}}
 \put(36,-120){\makebox(0,0){\normalsize $4$}}
 \put(50,-120){\circlearrowr}
 \put(36,-84){\vector(0,-1){30}}
 \put(36,-77){\makebox(0,0){\normalsize $3$}}
 \put(50,-77){\circlearrowr}
 \put(36,-41){\vector(0,-1){30}}

 \put(-50,-77){\circlearrowl}
 \put(-50,-120){\circlearrowl}
 \put(62,-80){\normalsize $\alpha_{3}$}
 \put(62,-123){\normalsize $\alpha_{4}$}
 \put(60,-35){\normalsize $\alpha_{2}$}
 \put(-62,-27){\normalsize $\alpha_{7}$}
 \put(-62,-46){\normalsize $\alpha_{1}$}
\end{picture}
\end{center}
and $\Gamma^{*}\lo 6\po$ is given by the above quiver and the
ideal $\overline{I}^{\lo 2\po }\!\lo T_{6},v_{1},v_{2}\po$ in
$K\overline{Q}_{T_{6}}^{2}$ generated by the elements:
$\beta_{1}\alpha_{2}$, $\alpha_{2}\beta_{2}$,
$\beta_{2}\alpha_{3}$, $\alpha_{3}\beta_{3}$,
$\beta_{3}\alpha_{4}$, $\alpha_{4}\beta_{4}$,
$\beta_{4}\alpha_{5}$, $\alpha_{5}\beta_{5}$,
$\beta_{5}\alpha_{6}$, $\alpha_{6}\beta_{6}$,
$\beta_{6}\alpha_{1}$, $\alpha_{7}\beta_{1}$,
$\gamma_{2}\beta_{1}$, $\beta_{6}\gamma_{1}$,
$\gamma_{2}\gamma_{1}$, $\alpha_{7}\alpha_{1}$,
$\lo\beta_{1}\beta_{2}\beta_{3}\beta_{4}\beta_{5}\beta_{6}
\po^{2}-\alpha_{1}\alpha_{7}$,
$\lo\beta_{2}\beta_{3}\beta_{4}\beta_{5}\beta_{6}\beta_{1}
\po^{2}-\alpha_{2}$,
$\lo\beta_{3}\beta_{4}\beta_{5}\beta_{6}\beta_{1}\beta_{2}
\po^{2}-\alpha_{3}$,
$\lo\beta_{4}\beta_{5}\beta_{6}\beta_{1}\beta_{2}\beta_{3}
\po^{2}-\alpha_{4}$,
$\lo\beta_{5}\beta_{6}\beta_{1}\beta_{2}\beta_{3}\beta_{4}
\po^{2}-\alpha_{5}$,
$\lo\beta_{6}\beta_{1}\beta_{2}\beta_{3}\beta_{4}\beta_{5}
\po^{2}-\alpha_{6}$,
$\alpha_{1}\alpha_{7}-\gamma_{1}\gamma_{2}$.
\end{exa}

The following theorem is the main result of \cite{bib:BoSk3}.
\begin{tw} \label{thm1par}
Let $A$ be a basic connected standard selfinjective algebra. Then
$A$ is one-parametric but not weakly symmetric if and only if $A$
is isomorphic to one of the algebras of the forms $\Omega^{\lo
1\po}\!\lo T,\sigma_{s},\theta\po$ or $\Omega^{\lo 2\po}\!\lo
T,v_{1},v_{2}\po$.
\end{tw}

For dealing with the algebras $\Omega^{(1)}\!\lo
T,\sigma_s,\theta\po$, the following definitions will be
convenient in the sequel. Let $i$ be an edge in $T$. We define the
\textit{level} of $i$ to be the distance of $i$ from the cycle of
$T$. In particular, the edges of level 0 are precisely the edges
on the cycle of $T$. For any edge $i$ of level $r\geqslant 2$ we
denote by $n\lo i\po$ the unique edge of level $r\! -\! 1$
adjacent to $i$. At each vertex $v_{1},\ldots,v_{k}$ on the cycle
of $T$ there are two Brauer trees attached, an inner one and an
outer one. For $u=1,\ldots,k$, we denote these Brauer trees by
$T_{u}^{inn}$ and $T_{u}^{out}$, respectively, with
$T_{u}^{inn}\cap T_{u}^{out}=\lk v_{u}\pk$. Recall that we denote
the edges of the cycle of $T$ by $e_u$, for $1\le u\le k$, where
$e_u$ has adjacent vertices $v_u$ and $v_{u+1}$.

The following result completely describes the Cartan matrices of
the algebras $\Omega^{(1)}\lo T,\sigma_{s},\theta\po$. We state it
for the convenience of the reader since it might be useful to
recall it when we are dealing with tilting complexes later. The
result follows directly from the defining relations as listed in
the definition above. Therefore the details of the proof are left
to the reader.

\begin{proposition} \label{cartan}
Given the algebra $\Omega^{\lo 1\po}\!\lo T,\sigma_{s},\theta\po$,
and edges $i$ of level $r$ and $j$ of level $t$ of the Brauer
graph $T$, the following holds for the Cartan entries $\dim\Hom\lo
P\lo i\po,P\lo j\po\po$.
\begin{arablist}
\item If $\left|r\! -\! t\right|\geqslant 2$, then $\dim\Hom\lo P\lo
i\po,P\lo j\po\po=0$.
\item Let $\left|r\! -\! t\right|=1$.
\begin{alphlist}\setlength{\itemindent}{3em}
\addtolength{\itemindent}{\labelwidth}
\item For $t=r\! -\! 1$, the following holds.
\begin{itemize}\setlength{\itemindent}{2em}
\addtolength{\itemindent}{\labelwidth}
\item Let $r\geqslant 2$. Then $\dim\Hom\lo P\lo i\po,P\lo
j\po\po=0$, unless $j=\sigma\lo n\lo i\po\po$, in which case
$\Hom\lo P\lo i\po,P\lo j\po\po$ is one-dimensional.
\item Let $r=1$ and $i\in T_{u}^{out}$ for some $1\leqslant u\leqslant
k$. Then $\dim\Hom\lo P\lo i\po,P\lo j\po\po=0$, unless
$j=\sigma\lo e_{u}\po$ or $j=\sigma\lo e_{u-1}\po$, in which case
$\Hom\lo P\lo i\po,P\lo j\po\po$ is one-dimensional.
\item Let $r=1$ and $i\in T_{u}^{inn}$ for some $1\leqslant u\leqslant
k$. Then
\begin{displaymath}
\dim\Hom(P(i),P(j)) = \lk
\begin{array}{ll} 2 & if \mbox{ } s=k\! -\! 1\mbox{ and } j=e_{u-1} \\
1 & if \mbox{ } s\neq k\! -\! 1\mbox{ and }\lo
j=e_{u-1}\mbox{ or }j=\sigma\lo e_{u}\po\po \\
0 & if \mbox{ } j\neq e_{u-1} \mbox{ and } j\neq \sigma\lo
e_{u}\po
\end{array}
\right.
\end{displaymath}
\end{itemize}
\item For $t=r\! +\! 1$, the following holds.
\begin{itemize}\setlength{\itemindent}{2em}
\addtolength{\itemindent}{\labelwidth}
\item Let $r\geqslant 1$. Then $\dim\Hom\lo P\lo i\po,P\lo
j\po\po=0$, unless $j$ is an edge adjacent to $i$, in which case
$\Hom\lo P\lo i\po,P\lo j\po\po$ is one-dimensional.
\item Let $r=0$ and $i=e_{u}$ for some $1\leqslant u\leqslant
k$. Then
\begin{displaymath}
\dim\Hom(P(i),P(j)) = \lk
\begin{array}{ll} 2 & if \mbox{ } s=k\! -\! 1\mbox{ and } j\in
T_{u}^{inn} \\
1 & if \mbox{ } \lo s\neq k\! -\! 1\mbox{ and } j\in T_{u}^{inn}\po \\
  & or \mbox{ } \lo j\in T_{u}^{out}\cup T_{u+1}^{out}\cup
T_{\sigma\lo u+1\po}^{inn}\po \\
0 & if \mbox{ } j\notin T_{u}^{out}\cup T_{u}^{inn}\cup
T_{u+1}^{out}\cup T_{\sigma\lo u+1\po}^{inn}
\end{array}
\right.
\end{displaymath}
\end{itemize}
\end{alphlist}
\item Let $\left|r\! -\! t\right|=0$.
\begin{itemize}\setlength{\itemindent}{1em}
\addtolength{\itemindent}{\labelwidth}
\item Let $r\geqslant 1$ and $i=j$. Then $\dim\Hom\lo P\lo i\po,P\lo
j\po\po=1$.
\item Let $r\geqslant 2$ and $i\in T_{u}^{out}\cup T_{u}^{inn}$
for some $1\leqslant u\leqslant k$. If $j\notin T_{u}^{out}\cup
T_{u}^{inn}\cup T_{\sigma\lo u\po}^{out}\cup T_{\sigma\lo
u\po}^{inn}$, then $\dim\Hom\lo P\lo i\po,P\lo j\po\po=0$.
Moreover, $\dim\Hom\lo P\lo i\po,P\lo j\po\po=0$, unless one of
the following cases holds, in which $\Hom\lo P\lo i\po,P\lo
j\po\po$ is one-dimensional.

(i) $j$ is a successor
of $i$ and a predecessor of $n\lo i\po$ in the cyclic order in the
graph $T$ around the common vertex of $i$ and $n\lo i\po$

(ii) $j$ is
a predecessor of $\sigma\lo i\po$ and a successor of
$n\lo\sigma\lo i\po\po=\sigma\lo n\lo i\po\po$ in the cyclic order
in the graph $T$ around the common vertex of $\sigma\lo i\po$ and
$n\lo\sigma\lo i\po\po$.

\item Let $r=1$ and $i\in T_{u}^{out}\cup T_{u}^{inn}$ for some
$1\leqslant u\leqslant k$. If $j\notin T_{u}^{out}\cup
T_{u}^{inn}\cup T_{\sigma\lo u\po}^{out}\cup T_{\sigma\lo
u\po}^{inn}$, then $\dim\Hom\lo P\lo i\po,P\lo j\po\po=0$.
Moreover, if $i\in T_{u}^{out}$ (respectively, $i\in
T_{u}^{inn}$), then $\dim\Hom\lo P\lo i\po,P\lo j\po\po=0$, unless
one of the following cases holds, in which $\Hom\lo P\lo i\po,P\lo
j\po\po$ is one-dimensional.

(i) $j\in T_{\sigma\lo
u\po}^{inn}$ (respectively, $j\in T_{u}^{out}$)

(ii) $j\in
T_{u}^{out}$ (respectively, $j\in T_{u}^{inn}$) and $j$ is
a successor of $i$ in the cyclic order in the graph $T$ around
$v_{u}$

(iii) $j\in T_{\sigma\lo u\po}^{out}$ (respectively,
$j\in T_{\sigma\lo u\po}^{inn}$) and $j$ is a predecessor of
$\sigma\lo i\po$ in the cyclic order in the graph $T$ around
$\sigma\lo v_{u}\po$.

\item Let $r=0$ and $i=e_{u}$. Then
\begin{displaymath}
\dim\Hom(P(i),P(j)) = \lk
\begin{array}{ll} 2 & if\mbox{ } s=k\! -\!1\mbox{ and }\lo
j=i\mbox{ or }j=e_{u-1}\po \\
1 & if\mbox{ } s\neq k\! -\!1 \\
  & and \mbox{ }\lo j=i\mbox{ or } j=e_{u-1}\mbox{ or }
j=\sigma\lo e_{u}\po\mbox{ or } j=\sigma\lo e_{u+1}\po\po\\
0 & if \mbox{ } j\neq i\mbox{ and } j\neq e_{u-1}\mbox{ and }
j\neq\sigma\lo e_{u}\po\mbox{ and } j\neq\sigma\lo e_{u+1}\po
\end{array}
\right.
\end{displaymath}
\end{itemize}
\item Let $i\in T_{u}^{out}$ of
level $\geqslant 2$, and let $j\in T_{w}^{inn}$ of level
$\geqslant 1$ for any $1\leqslant u,w\leqslant k$. Then
$\dim\Hom\lo P\lo i\po,P\lo j\po\po=0=\dim\Hom\lo P\lo j\po,P\lo
i\po\po$.
\item Let $i\in T_{u}^{inn}$ of level $\geqslant 2$, and let
$j\in T_{w}^{out}$ of level $\geqslant 1$ for any $1\leqslant
u,w\leqslant k$. Then $\dim\Hom\lo P\lo i\po,P\lo
j\po\po=0=\dim\Hom\lo P\lo j\po,P\lo i\po\po$.
\end{arablist}
\end{proposition}


\section{\hspace{-2.5ex}. Tilting complexes for Brauer graph algebras with a cycle}
\label{tilting}


In this section we are going to prove a first reduction step
towards the derived equivalence classification of the algebras
$\Omega^{(1)}\lo T,\sigma_s,\theta\po$. We show that, up to
derived equivalence, we can assume the Brauer trees attached to
the vertices of the cycle to be stars, i.e., all edges have level
at most 1. The proof will be based on the well-known construction
of derived equivalences for Brauer tree algebras \cite{bib:Ric2}.
But of course, the crucial new aspects are that for the algebras
$\Omega^{(1)}\lo T,\sigma_{s},\theta\po$ the Brauer graph contains
a cycle and that we have to take the automorphism $\sigma_{s}$
into account.

We shall need some notation. Let $\Omega^{(1)}\lo
T,\sigma_s,\theta\po$ be as defined above, with Brauer graph $T$
having exactly one cycle of length $k$, and rotation automorphism
$\sigma_s$. For abbreviation, we usually omit the index and just
write $\sigma$ for $\sigma_s$.

In a slight abuse of notation
we denote the successive edges on the cycle of $T$ just by
$1,2,\ldots,k$, instead of $e_1,\ldots,e_k$. Then the vertex
$v_{u}$ of $T$ lying on the cycle has adjacent edges $u\! -\! 1$
and $u$.

Recall that the edges of the Brauer graph $T$ correspond to the
simple modules of $\Omega^{(1)}\lo T,\sigma_s,\theta\po$. For an
edge $i$ we denote the corresponding projective indecomposable
module by $P(i)$.
\smallskip

The following notation will be very convenient and will be used
frequently throughout the paper. Let $i$ and $j$ be edges in $T$.
We denote by $[i,j]$ the 'interval' from $i$ to $j$ in the order
$p(T,\sigma)$ of edges of $T$, that is, the set of edges
$\{f~|~i\leqslant f\leqslant j\}$ (which might be empty).
Note that the terminology 'interval' is slightly misleading since
the set $[i,j]$ is not necessarily linearly ordered. We shall only
use this notation as a convenient way of denoting homomorphisms.
In fact, a sequence
of consecutive edges between $i$ and $j$ describes a path
in the Brauer quiver; multiplication with any such
path gives a homomorphism
$P(i)\to P(j)$ between the corresponding projective indecomposable
modules, which we also denote by $[i,j]:P(i)\to P(j)$.
Whenever we use this notation, it will always be clear from the
context which sequence of edges is meant.
As described in detail in Proposition \ref{cartan}, many of the Cartan
invariants are 1; in these cases the above notation for
homomorphisms is unambigious.


\smallskip
Our main result in this section is the following.

\begin{proposition} \label{redstep1}
Let $\Omega^{(1)}\lo T,\sigma,\theta\po$ be as defined above. Then
$\Omega^{(1)}\lo T,\sigma,\theta\po$ is derived equivalent to
$\Omega^{(1)}(\tilde{T},\sigma,\theta)$ with Brauer graph
$\tilde{T}$ of the following shape
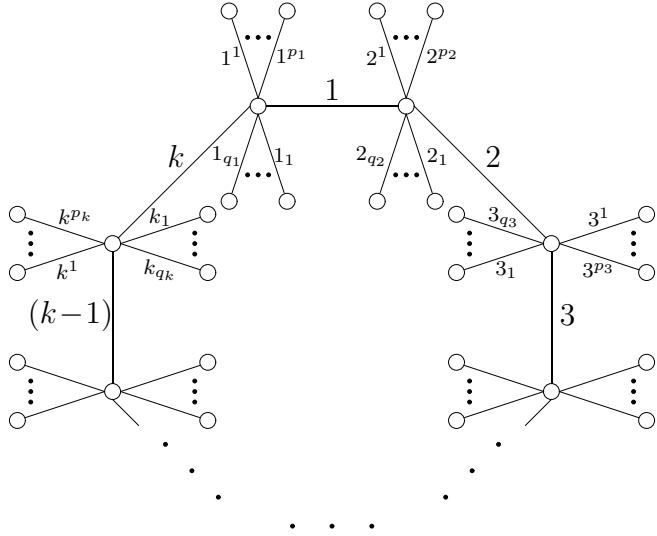
\begin{center}
\unitlength=1.000000pt
\begin{picture}(396.00,220.00)(-128.00,-115.00)
\put(81.00,59.00){\circle{6.00}} \put(-30.00,7.00){\circle{6.00}}
\put(136.00,7.00){\circle{6.00}}
\put(-30.00,-49.00){\circle{6.00}}
\put(136.00,-49.00){\circle{6.00}}
\put(25.00,59.00){\circle{6.00}}
\put(78.00,59.00){\line(-1,0){50.00}}
\put(84.00,59.00){\line(1,-1){50.00}}
\put(22.00,59.00){\line(-1,-1){50.00}}
\put(-30.00,-52.00){\line(1,-1){10.00}}
\put(136.00,-52.00){\line(-1,-1){10.00}}
\multiput(-10,-69)(10,-10){3}{\circle*{1.5}}
\multiput(116,-69)(-10,-10){3}{\circle*{1.5}}
\multiput(38,-100)(15,0){3}{\circle*{1.5}}
\put(-30.00,4.00){\line(0,-1){50.00}}
\put(136.00,4.00){\line(0,-1){50.00}}

\put(81,62){\line(-1,3){10}} \put(81,62){\line(1,3){10}}
\put(70,95){\circle{6.00}} \put(92,95){\circle{6.00}}
\multiput(77,85)(4,0){3}{\circle*{1.5}}
\put(81,56){\line(-1,-3){10}} \put(81,56){\line(1,-3){10}}
\put(70,23){\circle{6.00}} \put(92,23){\circle{6.00}}
\multiput(77,33)(4,0){3}{\circle*{1.5}}

\put(25,62){\line(-1,3){10}} \put(25,62){\line(1,3){10}}
\put(14,95){\circle{6.00}} \put(36,95){\circle{6.00}}
\multiput(21,85)(4,0){3}{\circle*{1.5}}
\put(25,56){\line(-1,-3){10}} \put(25,56){\line(1,-3){10}}
\put(14,23){\circle{6.00}} \put(36,23){\circle{6.00}}
\multiput(21,33)(4,0){3}{\circle*{1.5}}

\put(-33,7){\line(-3,1){30}} \put(-33,7){\line(-3,-1){30}}
\put(-66,18){\circle{6.00}} \put(-66,-4){\circle{6.00}}
\multiput(-61,3)(0,4){3}{\circle*{1.5}}
\put(-27,7){\line(3,1){30}} \put(-27,7){\line(3,-1){30}}
\put(6,18){\circle{6.00}} \put(6,-4){\circle{6.00}}
\multiput(1,3)(0,4){3}{\circle*{1.5}}

\put(-33,-49){\line(-3,1){30}} \put(-33,-49){\line(-3,-1){30}}
\put(-66,-38){\circle{6.00}} \put(-66,-60){\circle{6.00}}
\multiput(-61,-53)(0,4){3}{\circle*{1.5}}
\put(-27,-49){\line(3,1){30}} \put(-27,-49){\line(3,-1){30}}
\put(6,-38){\circle{6.00}} \put(6,-60){\circle{6.00}}
\multiput(1,-53)(0,4){3}{\circle*{1.5}}

\put(133,7){\line(-3,1){30}} \put(133,7){\line(-3,-1){30}}
\put(172,18){\circle{6.00}} \put(172,-4){\circle{6.00}}
\multiput(167,3)(0,4){3}{\circle*{1.5}}
\put(139,7){\line(3,1){30}} \put(139,7){\line(3,-1){30}}
\put(100,18){\circle{6.00}} \put(100,-4){\circle{6.00}}
\multiput(105,3)(0,4){3}{\circle*{1.5}}

\put(133,-49){\line(-3,1){30}} \put(133,-49){\line(-3,-1){30}}
\put(172,-38){\circle{6.00}} \put(172,-60){\circle{6.00}}
\multiput(167,-53)(0,4){3}{\circle*{1.5}}
\put(139,-49){\line(3,1){30}} \put(139,-49){\line(3,-1){30}}
\put(100,-38){\circle{6.00}} \put(100,-60){\circle{6.00}}
\multiput(105,-53)(0,4){3}{\circle*{1.5}}

\put(15.00,78.00){\makebox(0,0){\scriptsize $1^{1}$}}
\put(38.00,78.00){\makebox(0,0){\scriptsize $1^{p_1}$}}
\put(13.50,40.00){\makebox(0,0){\scriptsize $1_{q_1}$}}
\put(35.00,40.00){\makebox(0,0){\scriptsize $1_{1}$}}

\put(70.00,78.00){\makebox(0,0){\scriptsize $2^{1}$}}
\put(95.00,78.00){\makebox(0,0){\scriptsize $2^{p_2}$}}
\put(68.00,40.00){\makebox(0,0){\scriptsize $2_{q_2}$}}
\put(93.00,40.00){\makebox(0,0){\scriptsize $2_{1}$}}

\put(154.00,17.00){\makebox(0,0){\scriptsize $3^{1}$}}
\put(154.00,-3.00){\makebox(0,0){\scriptsize $3^{p_3}$}}
\put(118.00,17.00){\makebox(0,0){\scriptsize $3_{q_3}$}}
\put(119.00,-3.00){\makebox(0,0){\scriptsize $3_{1}$}}


\put(-44.00,17.00){\makebox(0,0){\scriptsize $k^{p_k}$}}
\put(-47.00,-3.00){\makebox(0,0){\scriptsize $k^{1}$}}
\put(-12.00,17.00){\makebox(0,0){\scriptsize $k_{1}$}}
\put(-12.00,-3.00){\makebox(0,0){\scriptsize $k_{q_k}$}}


\put(53.00,65.00){\makebox(0,0){$1$}}
\put(114.00,40.00){\makebox(0,0){$2$}}
\put(142.00,-20.00){\makebox(0,0){$3$}}
\put(-46.00,-20.00){\makebox(0,0){$\lo k\!-\!1\po$}}
\put(-6.00,40.00){\makebox(0,0){$k$}}
\end{picture}
\end{center}

\noindent where $p_{u}$ is the number of edges of
$T_{\sigma(u)+1}^{out}$ and
$q_{u}$ is the number of edges of $T_{u-1}^{inn}$.

In particular, in the Brauer graph $\tilde{T}$, all edges are of level
at most 1.
\end{proposition}

Note the difference between these algebras
$\Omega^{(1)}\lo\tilde{T},\sigma,\theta\po$ and the normal forms
$\Lambda(p,q,k,s,\theta)$ given in our main result Theorem
\ref{mainthm}. In this proposition we still have different numbers
of edges attached to the vertices of the cycle.
\smallskip

The rest of this section will be devoted to a proof of Proposition
\ref{redstep1}.
\smallskip

{\em Construction of the tilting complex.} For each edge $i$ of
$T$ we consider the following path from an edge on the cycle of
$T$ to $i$. If $i\in T_u^{out}$ then we take the unique shortest
path from $u-1$ to $i$. If $i\in T_u^{inn}$ we take the unique
shortest path from $u$ to $i$. Denote the edges on this path by
$i_0,i_1,\ldots,i_l=i$. By construction, there is a unique (up to
scalar multiplication) homomorphism $[i_j,i_{j+1}]:P(i_j)\to
P(i_{j+1})$ between the corresponding projective indecomposable
modules (see Proposition \ref{cartan}).

We then consider the following bounded complex of projective
$\Omega^{(1)}\lo T,\sigma,\theta\po$-modules
$$Q(i)~:~0\to P(i_0)\stackrel{[i_0,i_1]}{\longrightarrow}
P(i_1)\to \ldots \to P(i_{l-1})
\stackrel{[i_{l-1},i]}{\longrightarrow} P(i)\to 0$$
in which all maps are non-zero and where $P(i_{0})$ is in degree 0.
The entries in degree $-r$
are corresponding to edges of level $r$.

Note that the composition of any two maps is zero by the usual
Brauer tree relations which also hold in $\Omega^{(1)}\lo
T,\sigma_s,\theta\po$. For the edges on the cycle of $T$ this
construction just gives stalk complexes concentrated in degree 0.

\begin{proposition}
The complex $Q:=\bigoplus_{i\in T} Q(i)$ is a tilting complex for
$\Omega^{(1)}\lo T,\sigma,\theta\po$.
\end{proposition}
\begin{proof}
Exactly as in Rickard's construction for Brauer tree algebras
\cite{bib:Ric2} one sees that $\add(Q)$ generates the homotopy
category $K^b\lo\Omega^{(1)}\lo T,\sigma,\theta\po\po$ as a
triangulated category.

We have to show that in the homotopy category, $\Hom(Q,Q[n])=0$
for all $n\neq 0$.

By Proposition \ref{cartan}\,(1), it is clear that $\Hom(Q,Q[n])=0$
unless $n\in\{-1,0,1\}$.

We first consider a map $\alpha$ of complexes from $Q(i)$ to
$Q(j)[-1]$. Hence we have a commutative diagram of the form
$$\begin{array}{ccccccccccccccc}
0 & \to & P(i_0) & \to & P(i_1) & \to & \ldots & \to &
 P(i_r) & \to & P(i_{r+1}) & \to & P(i_{r+2}) & \to & \ldots \\
 & & & & \downarrow\mbox{\hskip-.1cm{\footnotesize $\alpha_0$}}
& & & & \downarrow & & \downarrow
\mbox{\hskip-.1cm{\footnotesize $\alpha_r$}} &
 & \downarrow & & \\
 & & 0 & \to & P(j_0) & \to & \ldots & \to &
 P(j_{r-1}) & \to & P(j_{r}) & \to & P(j_{r+1}) & \to & \ldots \\
\end{array}
$$
Let $r$ be minimal such that $\alpha_r\neq 0$. By Proposition
\ref{cartan}, the composition $P(i_r)\to P(i_{r+1}) \to P(j_{r})$
is non-zero. In fact, if the level $r\ge 1$ then this
composition is the socle map; if $r=0$ then there are two
possibilities for $j_0$, but in both
cases the composition is non-zero (recall that the map $P(i_0)\to
P(i_1)$ is given by multiplication with the {\em shortest} path from
$i_0$ to $i_1$).
This contradicts the fact that
$\alpha$ is a map of complexes, which
proves $\Hom(Q,Q[-1])=0$, as desired.
\smallskip

Secondly, we have to consider maps $\alpha$ of complexes
from $Q(i)$ to $Q(j)[1]$. We then have a commutative diagram of the form
$$\begin{array}{ccccccccccccccc}
& & 0 & \to & P(i_0) & \to & \ldots & \to &
 P(i_{r-1}) & \to & P(i_{r}) & \to & P(i_{r+1}) & \to & \ldots \\
 & & & & \downarrow\mbox{\hskip-.1cm{\footnotesize $\alpha_0$}}
 & & & & \downarrow & &
\downarrow\mbox{\hskip-.1cm{\footnotesize $\alpha_r$}}
 & & \downarrow  & & \\
0 & \to & P(j_0) & \to & P(j_1) & \to & \ldots & \to &
 P(j_r) & \to & P(j_{r+1}) & \to & P(j_{r+2}) & \to & \ldots \\
\end{array}
$$
Let $r$ be maximal with the property $\alpha_r\neq 0$. Moreover,
we choose $\alpha$ from its homotopy class such that this value of
$r$ is minimal. Assume first that $r\geqslant 1$. The existence of
a non-zero map $\alpha_r$ from level $r\geqslant 1$ to level $r+1$
implies that $j_{r+1}$ is adjacent to $i_r$ (see Proposition
\ref{cartan}\,(2b)). But this implies $j_r=i_r$, which means that
$\alpha_r$ factors through $P(j_r)$. This factorisation gives a
homotopy from $\alpha$ to a map $\beta$ for which the minimal
non-zero map occurs in a degree strictly smaller than for
$\alpha$. This contradicts the choice of $\alpha$ as being of
minimal such degree in its homotopy class.

Now assume that $r=0$. There are four possibilities for a non-zero
map from level 0 to level 1: either $j_1$ is in $T_{i_0+1}^{out}$,
in $T_{i_0}^{inn}$, in $T_{i_0}^{out}$
or in $T_{\sigma(i_0)+1}^{inn}$. In the first two cases we deduce
that $j_0=i_0$, in the third case that $j_0=i_0-1$, and
in the latter case that $j_0=\sigma(i_0)+1$.
In all cases, $\alpha_0$ factors through $P(j_0)$, implying that
$\alpha$ is homotopic to zero.

Hence, $\Hom(Q,Q[1])=0$ in the homotopy category, as desired.

This completes the proof that $Q$ is a tilting complex for
$\Omega^{(1)}\lo T,\sigma,\theta\po$.
\end{proof}
\medskip

By Rickard's criterion \cite{bib:Ric1} the endomorphism ring of
the tilting complex $Q$ for $\Omega^{(1)}\lo T,\sigma,\theta\po$
is derived equivalent to the algebra $\Omega^{(1)}\lo
T,\sigma,\theta\po$. So it remains to prove that this
endomorphism ring has the structure described in Proposition
\ref{redstep1}.

\begin{proposition} \label{redstep1b}
Let $\Omega^{(1)}\lo T,\sigma,\theta\po$ be as above, with
Brauer trees attached to the vertices of the cycle denoted by
$T_1^{out},\ldots,T_k^{out}$ and $T_1^{inn},\ldots,T_k^{inn}$,
where $T_u^{out}$ and $T_u^{inn}$ are attached to the vertex
between the edges $u-1$ and $u$ on the cycle.

Let $Q=\oplus_{i\in T} Q(i)$ be the tilting complex for
$\Omega^{(1)}\lo T,\sigma,\theta\po$ as defined above. Then the
endomorphism ring of $Q$ is isomorphic to the algebra
$\Omega^{(1)}\lo\tilde{T},\sigma,\theta\po$, where the Brauer
graph $\tilde{T}$ can be described by the following properties.
\begin{romalist}
\item The cycle of $\tilde{T}$ is the same as the cycle of $T$.
\item All edges of $\tilde{T}$ have level at most 1.
\item Denote by $E_u^{out}$ and $E_u^{inn}$
the numbers of edges of the Brauer trees $T_u^{out}$ and
$T_u^{inn}$, respectively. In $\tilde{T}$, attached to the vertex
adjacent to the edges $u-1$ and $u$, we have
$E_{\sigma(u)+1}^{out}$ outer edges and $E_{u-1}^{inn}$ inner
edges.
\end{romalist}
\end{proposition}

Note that compared with the algebra $\Omega^{(1)}\lo
T,\sigma,\theta\po$ the inner Brauer trees have been shifted by
one vertex on the cycle, and the outer Brauer trees have been
shifted by $-(s+1)$ where $\sigma=\sigma_s$. This means that the
order $p(\tilde{T},\sigma)$ on the edges of $\End(Q)$ is different
from the order of edges of the original algebra $\Omega^{(1)}\lo
T,\sigma,\theta\po$. But note that the automorphism $\sigma$ is
the same for both $\Omega^{(1)}\lo T,\sigma,\theta\po$ and the
endomorphism ring of $Q$.
\medskip

The rest of this section is devoted to a proof of Proposition
\ref{redstep1b}. The proof consists of three steps. We first
describe morphisms between direct summands of $Q$, which
correspond to the arrows of the algebra
$\Omega^{(1)}\lo\tilde{T},\sigma,\theta\po$. Secondly, we show
that these maps satisfy the defining relations of
$\Omega^{(1)}\lo\tilde{T},\sigma,\theta\po$. In the third step we
show that $\End(Q)$ is actually isomorphic to
$\Omega^{(1)}\lo\tilde{T},\sigma,\theta\po$. To this end, we need
to show that there are no further relations in $\End(Q)$.
\medskip

\noindent {\em I. The morphisms between summands of $Q$.} The
direct summands of the tilting complex $Q$ are complexes defined
completely analogous to Rickard's construction for the classical
Brauer tree algebras. Recall that for a Brauer tree algebra, the
endomorphism ring was a star; more precisely, the cyclic ordering
of the edges for the endomorphism ring can be obtained by the
famous 'walk around the Brauer tree' (see \cite{bib:Green}).
 Note that in contrast to the usual Brauer
tree setting, in our case one can not walk completely around the
Brauer tree since our algebras are not symmetric. Still, one gets
for the endomorphism ring of the complexes corresponding to any
Brauer tree $T_u^{out}$ and $T_u^{inn}$ a cyclic ordering of the
edges. In other words, in $\End(Q)$, all edges will have level at
most 1.

Let us now define explicitly
the morphisms between direct summands of $Q$
corresponding to successive edges in the order $p(\tilde{T},\sigma)$.
\smallskip

(1) {\em The maps inside the Brauer trees.}

Within any of the Brauer trees $\tilde{T}_u^{out}$ and
$\tilde{T}_u^{inn}$ of $\End(Q)$
these morphisms are given as in Rickard's
construction by a walk around the Brauer tree. Actually, for our proof
we don't need to know these maps in detail; for checking relations it
will
only be important that their composition in degree 0 is the identity
(which is immediate from Rickard's classical construction).

Next, we show that we actually have non-zero sequences of morphisms
starting from $Q(u)$ through all
complexes corresponding to vertices of $T_{u-1}^{inn}$ to
$Q(u-1)$ and further on through all complexes corresponding to
edges of $T_{\sigma(u)+1}^{out}$ and finally to $Q(\sigma(u))$.
\smallskip

(2) {\em The maps $Q(u)\to Q(T_{u-1}^{inn})$.}

For any edge $u$ on the cycle consider $[u,u-1]$, that is, the
sequence of edges $u$, consecutive edges of level 1 in $T_u^{inn}$
and $u-1$.
Then multiplication by $[u,u-1]$ gives a homomorphism $P(u)\to
P(u-1)$.

Then for any $f$ in $T_{u-1}^{inn}$ we have
a homomorphism of complexes $Q(u)\to Q(f)$ given by
multiplication with $[u,u-1]$ in degree 0.

In fact, this is a map
of complexes since $Q(f)$ starts with a map from $P(u-1)$ given
by walking around the endpoint of $u-1$ not adjacent to $u$.
Since the edges are not consecutive in the order
$p(T,\sigma)$, the composition with $[u,u-1]$ is zero.
\smallskip

(3) {\em The maps $Q(T_{u-1}^{inn})\to Q(u-1)$.}

Let $f$ be in
$T_{u-1}^{inn}$. There is an obvious map from $Q(f)$
to the stalk complex $Q(u-1)$ given by the identity in
degree 0.
\smallskip

(4) {\em The maps $Q(u-1)\to Q(T_{\sigma(u)+1}^{out})$.}

Let $q_{u}$ denote the map $P(u-1)\to P(\sigma(u))$ given by
multiplication with $[u-1,\sigma(u)]$, i.e., the sequence of edges $u$,
all edges of level 1 in $T_u^{out}$, and $\sigma(u)$.

Let $g$ be an edge in $T_{\sigma(u)+1}^{out}$. We have a
homomorphism of complexes $Q(u-1)\to Q(g)$ given by multiplication
with $q_{u}$ in degree 0.

In fact, this is a map of complexes since $Q(g)$ starts with a map
from $P(\sigma(u))$ given by an edge not consecutive to the last
edge of $q_{u}$ in the order $p(T,\sigma)$.
\smallskip

(5) {\em The maps $Q(T_{\sigma(u)+1}^{out})\to Q(\sigma(u))$.}

Finally, for any $g$ in $T_{\sigma(u)+1}^{out}$ there
is an obvious map from $Q(g)$ to the stalk complex
$Q(\sigma(u))$ given by the identity in degree 0.
\medskip

The composition of the maps just described is the map from
$Q(u)$ to $Q(\sigma(u))$ given by the socle map.
In particular, this composition is non-zero in $\End(Q)$, as desired.
\medskip

%
%
%
%
%

(6) {\em The maps around the endpoints of $\tilde{T}$.}

So far, we have described maps for all consecutive elements in the
ordering $p(\tilde{T},\sigma)$ around the vertices on the cycle.
It remains to define maps corresponding to the endpoints of the
edges of level 1 of the Brauer graph $\tilde{T}$. The edges of level 1
in $\tilde{T}$ correspond to complexes $Q(f)$ where
$f$ is any edge in $T$ not on the cycle.

We define the map $Q(f)\to Q(\sigma(f))$
to be multiplication with the socle map
$[f,\sigma(f)]$ in degree 0 and the zero map in all other degrees.

We shall show below that this map of complexes is not homotopic to
zero, that is, it is non-zero in $\End(Q)$.
\medskip

We have thus described maps of complexes corresponding to all
consecutive edges in the order $p(\tilde{T},\sigma)$. In other
words, we have given maps corresponding to all arrows in the
quiver of $\Omega^{(1)}\lo\tilde{T},\sigma,\theta\po$.
\medskip

\noindent
{\em II. Checking the defining relations.}

We have to show that the maps defined in (1)-(6) above actually
satisfy the defining relations of
$\Omega^{(1)}\lo\tilde{T},\sigma,\theta\po$ (see Section
\ref{morita} for the definition). The Brauer trees attached to
vertices on the cycle of $\tilde{T}$ are denoted by
$\tilde{T}_u^{out}=Q(T_u^{out})$ and
$\tilde{T}_u^{inn}=Q(T_u^{inn})$, respectively. Hence, attached to
the vertex in $\tilde{T}$ adjacent to $u-1$ and $u$ we have
$\tilde{T}_{\sigma(u)+1}^{out}$ and $\tilde{T}_{u-1}^{inn}$.
\smallskip

(i) {\em Products of non-consecutive edges.}

Let $f\in\tilde{T}_{u+1}^{out}$ and $g\in\tilde{T}_{\sigma u+2}^{out}$.
Then $f,u,g$ are not consecutive in the order of edges of $\tilde{T}$.
We have to show that the corresponding composition of the maps
defined above in (5) and (4) is homotopic to zero.
The composition is given in degree 0 by
$$Q(f)\stackrel{id}{\longrightarrow} Q(u)
\stackrel{[u,\sigma u+1]}{\longrightarrow} Q(g)$$
and by the zero map in all other degrees.
Write $Q(f)\,:\,0\to P(u)\to P(i_1)\to\ldots$
and $Q(g)\,:\,0\to P(\sigma(u)+1)\to P(j_1)\to\ldots$.

A suitable homotopy map from $Q(f)$ to $Q(g)$ is given by
$[i_1,\sigma\lo u\po+1]:P(i_1)\to P(\sigma\lo u\po+1)$ from degree
$-1$ to degree 0 and the zero maps in all other degrees. In fact,
in degree 0 we clearly get $[u,\sigma\lo u\po+1]$, and in degree
$-1$ we indeed get the zero map since $i_1,\sigma\lo u\po+1,j_1$
are not consecutive in the order of edges of $T$ ($i_1$ and $j_1$
are both outer edges).
\smallskip

The argument for inner edges works similarly. If
$f\in\tilde{T}_u^{inn}$ and $g\in\tilde{T}_{u-1}^{inn}$ then the
composition $Q(f)\to Q(u)\to Q(g)$ is given by $[u,u-1]$ in degree
0 and the zero maps in all other degrees. This composition is
homotopic to zero via the homotopy map $[i_1,u-1]:P(i_1)\to
P(u-1)$ from degree $-1$ to degree 0, and zero in all other
degrees. (Use that $i_1,u,j_1$ are not consecutive since $i_1$ and
$j_1$ are both inner edges.)
\smallskip

(ii) {\em Socle relations.}

Let $f$ be any edge not on the cycle. The edge in $\tilde{T}$
corresponding to the complex $Q(f)$ is of level 1. So there are
two maps $Q(f)\to Q(\sigma(f))$, given by the orders around the
two endpoints of this edge. Using the definitions in (1)-(6) above
one directly sees that both compositions are given by
the socle map $[f,\sigma(f)]$ in degree 0 and the zero maps in all other
degrees.
\smallskip


Now consider an edge $u$ on the cycle, with endpoints $v$ and $w$.
Then by the above definitions, the compositions of maps
$Q(u)\to Q(\sigma(u))$ in $\End(Q)$ around the two endpoints are
given by multiplication with $C(u,p(T,\sigma,v))$ and
$C(u,p(T,\sigma,w))$ (i.e. by the socle maps $P(u)\to
P(\sigma(u))$). By definition of the relations in
$\Omega^{(1)}(T,\sigma,\theta)$ these maps are equal for
all $u\neq 1$ and for the distinguished vertex $u=1$
they are equal up to multiplication by the scalar $\theta$. Hence
we get exactly the desired relations for
$\Omega^{(1)}(\widetilde{T},\sigma,\theta)$.
\medskip

\noindent
{\em III. No further relations.}

In the final step, we have to prove that there are no further
relations in $\End(Q)$.
To this end it suffices to show that for all edges $f$
the socle map $Q(f)\to Q(\sigma(f))$ is not homotopic to zero,
that is, non-zero in $\End(Q)$. This is clear for the edges on the
cycle, since the complexes are stalk complexes and the socle map
is just the socle map on the projective indecomposable modules.
\smallskip

Let $f$ be any edge in $\tilde{T}_{u+1}^{out}$ for some $u$.
The corresponding complex $Q(f)$ has the form
$$0\to P(u)\stackrel{[u,i_1]}{\longrightarrow} P(i_1)
\stackrel{[i_1,i_2]}{\longrightarrow} P(i_2)\longrightarrow \ldots
\longrightarrow P(i_{l-1})\stackrel{[i_{l-1},f]}{\longrightarrow}
P(f)\to 0$$ where $u,i_1,i_2,\ldots,i_{l-1},f$ is the unique
shortest walk in $T$ from the edge $u$ on the cycle to $f$. The
socle map $\alpha_f:Q(f)\to Q(\sigma \lo f\po)$ is given by the
zero map in all non-zero degrees and in degree 0 it is given by
the composition $Q(f)\stackrel{id}{\longrightarrow} Q(u)
\stackrel{[u,u-1]}{\longrightarrow} Q(u-1) \stackrel{[u,\sigma
u]}{\longrightarrow} Q(\sigma\lo f\po)$, that is, by the socle map
$[u,\sigma\lo u\po]:P(u)\to P(\sigma\lo u\po)$. We claim that this
map of complexes $Q(f)\to Q(\sigma\lo f\po)$ is not homotopic to
zero. In fact, assuming there were a homotopy map, then for degree
0 this map has to be $[i_1,\sigma\lo u\po]:P(i_1)\to P(\sigma\lo
u\po)$. Since $\alpha_f$ is 0 in degree -1, the homotopy map has
to be given by $-[i_2,\sigma\lo i_1\po]: P(i_2)\to P(\sigma \lo
i_1\po)$. Inductively, the final homotopy map $P(f)\to P(i_{l-1})$
must be given by $\pm [f,\sigma\lo i_{l-1}\po]$. But then in
degree $-l$ the composition $P(f)\to P(\sigma\lo i_{l-1}\po) \to
P(\sigma\lo f\po)$ is the non-zero socle map, a contradiction,
since $\alpha_f$ is zero in degree $-l$. Hence, $\alpha_f$ is not
homotopic to zero, as desired.
\smallskip

A completely analogous argument shows that, for any inner edge
$g$, the socle map $Q(g)\to Q(\sigma\lo g\po)$ is not homotopic to
zero. We leave the details to the reader.
\smallskip

This completes the proof of Proposition \ref{redstep1b}
$\hfill\Box$


\section{\hspace{-2.5ex}. Moving edges} \label{moving}


In the previous section we have proved that any algebra
$\Omega^{(1)}\lo T,\sigma,\theta\po$ is derived equivalent to an
algebra $\Omega^{(1)}\lo\tilde{T},\sigma,\theta\po$ where all
edges in $\tilde{T}$ are of level at most 1 (i.e., the Brauer
trees attached to the vertices on the cycle are stars).

In this section we shall prove a further reduction step
by showing that, up to derived equivalence, we can concentrate
all edges of level 1 around the vertices of a {\em single}
$\sigma$-orbit.

There is nothing to prove if the $\sigma$-action is already
transitive. So in this section we can assume that the action of
$\sigma$ on the vertices of the cycle is not transitive.

The main result of this section is the following.

\begin{proposition}\label{redstepmoving}
Let $A=\Omega^{(1)}\lo\tilde{T},\sigma,\theta\po$ be as in
Proposition \ref{redstep1}, where all edges are of level at most
1. Fix any vertex $v$ on the cycle. Then $A$ is derived equivalent
to the algebra $\Omega^{(1)}\lo\overline{T},\sigma,\theta\po$
where the Brauer graph $\overline{T}$ has the following
properties:
\begin{romalist}
\item The cycle of $\overline{T}$ is the same as the
cycle of $\widetilde{T}$.
All edges of $\overline{T}$ are of level at most 1.
\item All edges of level 1 are attached to some vertex
in the $\sigma$-orbit of the fixed vertex $v$.
\end{romalist}
\end{proposition}

This result will follow immediately from the following two lemmas.
The first lemma asserts that, up to derived equivalence, we can
move inner edges by one vertex on the cycle. The analogous second
lemma states that outer edges can be moved by $-(s+1)$ steps on
the cycle, where $\sigma=\sigma_s$. By applying these two results
about moving inner and outer edges inductively, we can clearly
move all edges of level 1 in $\tilde{T}$ to some vertex in the
$\sigma$-orbit of the fixed vertex $v$.

\begin{lemma} \label{redstep2}
Let $A=\Omega^{(1)}\lo\tilde{T},\sigma,\theta\po$ be as above with
all edges of level at most 1. Fix any edge $u$ on the cycle, and
denote the subsequent inner edges (in the order
$p(\tilde{T},\sigma)$) by $u_1,\ldots,u_{q_{u}}$, and the inner
edges subsequent to $u+1$ by $\lo u+1\po_1,\ldots,\lo
u+1\po_{q_{u+1}}$. Then $A$ is derived equivalent to an algebra
$\Omega^{(1)}\lo\overline{T},\sigma,\theta\po$, where the order
$p(\overline{T},\sigma)$ is given as follows. For any integer $j$
the order around the common vertex of $\sigma^j(u+1)$ and
$\sigma^j(u)$ is given by $\ldots,\sigma^j(u+1),\sigma^j(\lo
u+1\po_1),\ldots,\sigma^j(\lo u+1\po_{q_{u+1}}),
\sigma^j(u_1),\ldots,\sigma^j(u_{q_{u}}),\sigma^j(u),\ldots$, and
around the other endpoint of $\sigma^j(u)$ by
$\sigma^j(u),\sigma^j(u-1)$.

In other words, the inner edges
$\{\sigma^j(u_1),\ldots,\sigma^j(u_{q_{u}})~|~j\in\mathbb{Z}\}$
can be moved to the next vertex on the cycle, up to derived
equivalence (and all other edges remain unchanged).
\end{lemma}
\begin{proof}
We consider the following bounded complexes of projective
$A$-modules. For any edge $z\not\in
S:=\{\sigma^j(u_i)~|~j\in\mathbb{Z},1\leqslant i\leqslant q_{u}\}$
let $Q(z)$ be the stalk complex with
the corresponding projective indecomposable module
$P(z)$ in degree 0. For the edges in the $\sigma$-orbit of
$u_1,\ldots,u_{q_{u}}$ we define the following complexes
$$Q(\sigma^j(u_i))~:~0\to P(\sigma^j(u))\stackrel{[\sigma^j(u),
\sigma^j(u_i)]}{\longrightarrow} P(\sigma^j(u_i))\to 0$$
in degrees 0 and $-1$.
\smallskip

{\em Claim: The direct sum $Q:=\oplus_{z\in\tilde{T}} Q_{z}$ is a
tilting complex for $\Omega^{(1)}\lo\tilde{T},\sigma,\theta\po$.}
\smallskip

{\em Proof of the claim.} From the construction it is clear that
$\add(Q)$ generates the homotopy category
$K^b\lo\Omega^{(1)}\lo\tilde{T},\sigma,\theta\po\po$ as
triangulated category. It remains to show that, in the homotopy
category, $\Hom(Q,Q[n])=0$ for all $n\neq 0$. This is clear for
$|n|\geqslant 2$, since we are only dealing with two-term
complexes. From the relations of $A$ (cf. also Proposition
\ref{cartan}) we can directly deduce that all non-zero maps
$P(z)\to P(\sigma^j(u_i))$ with $z\notin S$ factor through
$P(\sigma^j(u))$. Hence, $\Hom(Q,Q[1])=0$. Similarly, one sees
that any non-zero map $P(\sigma^j(u_i))\to P(z)$ with $z\notin S$
remains non-zero when composed with $[\sigma^j(u),\sigma^j(u_i)]$.
Hence, $\Hom(Q,Q[-1])=0$, which proves the claim.
\smallskip

{\em The endomorphism ring of $Q$.} By Rickard's theorem, the
endomorphism ring of $Q$ (in the homotopy category) is derived
equivalent to $A$. So in order to complete the proof we have to
show that $\End(Q)$ has the structure described above, where the
edges in the $\sigma$-orbit of $u_1\ldots,u_{q_{u}}$ are moved by
one step to the next vertex on the cycle (and all other edges
remain unchanged). Since 'most' summands of $Q$ are stalk
complexes, it is clear that the Brauer graph for $\End(Q)$ will
only change as far as the edges in the $\sigma$-orbit of
$u_1\ldots,u_{q_{u}}$ are concerned.

We define the following maps of complexes between summands of $Q$
involving some $Q(\sigma^j(u_i))$. For simplicity we describe them
for the $Q(u_i)$; the other maps are obtained by applying
$\sigma^j$. Let $Q(\lo u+1\po_{q_{u+1}})\to Q(u_1)$ be defined by
multiplication with $[\lo u+1\po_{q_{u+1}},u]$ in degree 0. Note
that this is indeed a map of complexes since $\lo
u+1\po_{q_{u+1}},u,u_1$ are not consecutive in the ordering of
edges of $\tilde{T}$. For any $i\in\{1,\ldots,q_{u}-1\}$ we define
a map of complexes $Q(u_i)\to Q(u_{i+1})$ by the identity in
degree 0 and $[u_i,u_{i+1}]$ in degree $-1$. Finally, we define a
map $Q(u_{q_{u}})\to Q(u)$ by the identity in degree 0 and a map
$Q(u)\to Q(u-1)$ by $[u,u-1]$. We denote the edges of the Brauer
graph for $\End(Q)$ corresponding to $Q(\sigma^j(u_i))$ by
$\overline{\sigma^j(u_i)}$, and the 'old' edges corresponding to
the stalk complexes as before. All relations of
$\Omega^{(1)}\lo\overline{T},\sigma,\theta\po$ not involving these
new complexes immediately follow from the corresponding same
relations in $\Omega^{(1)}\lo\tilde{T},\sigma,\theta\po$. For the
other relations, by applying $\sigma^j$, it suffices to check the
relations involving $\overline{u_i}$. The edges
$\overline{u_{q_{u}}},u,u-1$ are claimed not to be consecutive in
the order of $\End(Q)$. In fact, the corresponding map
$Q(u_{q_{u}})\to P(u)\to P(u-1)$ is given by $[u,u-1]$ in degree
0. It is homotopic to zero via the homotopy map $[u_{q_{u}},u-1]$,
giving the desired zero relation in $\End(Q)$. Since the maps
$Q(u_i)\to Q(u_{i+1})$, for $i\in\{1,\ldots,q_{u}-1\}$, are in
degree 0 just the identity, all relations starting at any of the
edges $\lo u+1\po_j$ still hold in $\End(Q)$ with the new order of
edges. Moreover, the composition of the socle map
$Q(\sigma^{-1}(\lo u+1\po_{q_{u+1}}))\to Q(\lo u+1\po_{q_{u+1}})$
with the map $Q(\lo u+1\po_{q_{u+1}})\to Q(u_1)$ is zero, since
the latter is given by $[\lo u+1\po_{q_{u+1}},u]$.

It remains to show that the socle relations hold for the edges
in $\End(Q)$ corresponding to the complexes $Q(\sigma^j(u_i))$.
Again, for simplicity we only consider $\overline{u_i}=Q(u_i)$.
Note that the new order of edges on $\End(Q)$ is claimed to have
the following consecutive edges:
$$\overline{u_i},\ldots,\overline{u_{q_{u}}},
[u,\sigma(u+1)],\sigma(\lo u+1\po_1),\ldots,\sigma(\lo
u+1\po_{q_{u+1}}), \sigma(\overline{u_1}),
\ldots,\sigma(\overline{u_{i-1}}),\sigma(\overline{u_i}).$$ The
corresponding composition $Q(u_i)\to Q(\sigma(u_i))$ is given by
the socle map $P(u)\to P(\sigma(u))$ in degree 0 and the zero map
in degree $-1$. We claim that this composition is not homotopic to
zero. In fact, by the relations of
$\Omega^{(1)}\lo\tilde{T},\sigma,\theta\po$, this socle map is
given by the order around one of the endpoints of the edge $u$.
So in degree 0 it could only be factored by the homotopy map
$[u_i,\sigma(u)]$, but then in degree $-1$ this gives the non-zero
socle map $P(u_i)\to P(\sigma(u_i))$, which proves the claim. Note
that if we compose this socle map $Q(u_i)\to Q(\sigma(u_i))$ with
the next map $Q(\sigma(u_i))\to Q(\sigma(u_{i+1}))$, then this
composition becomes homotopic to zero via the homotopy map
$[u_i,\sigma(u)]$.

Note that for the socle relation on $P(e_1)$, where $e_1$
is the distinguished vertex on the cycle, involving the
scalar $\theta$ , the corresponding relation on $Q(e_1)$
in $\End(Q)$ involves the {\em same} scalar $\theta$.
In fact, the above defined maps on the $Q(u_i)$'s are
the identity in degree 0, so nothing changes for the
stalk complexes.

This completes the proof of the lemma
\end{proof}
\medskip

There is the following analogous result for moving outer edges.
Since the proof is very similar to the proof of Lemma \ref{redstep2}
we leave the details to the reader.

\begin{lemma} \label{redstep2b}
Let $A=\Omega^{(1)}\lo\tilde{T},\sigma,\theta\po$ be as above with
all edges of level at most 1. Fix any edge $u$ on the cycle, and
denote the outer edges succeeding $u-1$ (in the order
$p(\tilde{T},\sigma)$) by $u^1,\ldots,u^{p_u}$, and the outer
edges succeeding $u-2$ by $\lo u-1\po^1,\ldots,\lo
u-1\po^{p_{u-1}}$. Then $A$ is derived equivalent to an algebra
$\Omega^{(1)}\lo\overline{T},\sigma,\theta\po$ where the order
$p(\overline{T},\sigma)$ of edges is given as follows. For any
integer $j$ the order around the common vertex of $\sigma^j(u-2)$
and $\sigma^j(u-1)$ is given by
$$\ldots,\sigma^j(u-2),\sigma^j(\lo u-1\po^1),
\ldots,\sigma^j(\lo u-1\po^{p_{u-1}}),
\sigma^{j+1}(u^1),\ldots,\sigma^{j+1}(u^{p_u}),\sigma^{j+1}(u-1),
\ldots$$
and around the other endpoint of $\sigma^j(u-1)$ by
$\ldots,\sigma^j(u-1),\sigma^{j+1}(u),\ldots$.

In other words, all outer edges in the $\sigma$-orbit of
$u^1,\ldots,u^{p_u}$ have been shifted by $-(s+1)$ vertices on the cycle
(and all other edges remain unchanged). $\hfill\Box$
\end{lemma}


\section{\hspace{-2.5ex}. Reduction to the transitive case} \label{reduction}


In the last two sections we proved that any algebra
$\Omega^{(1)}\lo T,\sigma,\theta\po$ is derived equivalent to an
algebra $\Omega^{(1)}\lo\overline{T},\sigma,\theta\po$ where all
edges of the Brauer graph $\overline{T}$ are of level at most 1,
and all edges of level 1 are attached to vertices in a {\em
single} $\sigma$-orbit.

The main aim of this section is to prove that, up to derived
equivalence, we can even assume that the $\sigma$-action is {\em
transitive}. If the $\sigma$-action on
$\Omega^{(1)}\lo\overline{T},\sigma,\theta\po$ is not already
transitive then we will have vertices on the cycle of
$\overline{T}$ with no edges of level 1 attached. Derived
equivalences for this situation are provided by the following
result. The crucial aspect of this result is that if the
$\sigma$-action is not transitive, then one can reduce the number
of edges on the cycle.

Recall that the edge $e_1$ on the cycle plays a special role in
that the socle relation involves the scalar $\theta$. The results
in the preceding section allow us to concentrate the edges in
the $\sigma$-orbit of {\em any} vertex on the cycle. So for
our purposes we can assume that the distinguished edge $e_1:=1$
has no edges attached to it, i.e. the technical assumption in the
following result can always be guaranteed in our context.

\begin{proposition} \label{redstep3} Let
$A=\Omega^{(1)}\lo\overline{T},\sigma,\theta\po$ be as above, with
all edges of $\overline{T}$ of level at most 1. Assume that for
some edge $u$ on the cycle there are no edges of level 1 attached
to the vertices between $u-1$ and $u$ and between $u$ and $u+1$
(that is, in the ordering of edges we have that $u,u-1,\sigma(u)$
are consecutive and $u+1,u,\sigma(u+1)$ are consecutive).
Moreover, assume that $u-1$ and $u+1$ are not in the
$\sigma$-orbit of the distinguished edge.

Then $A$
is derived equivalent to an algebra $\Omega^{(1)}\lo
T^*,\sigma^*,\theta\po$ where all edges of $T^*$ are of level at
most 1, and where all edges in the $\sigma$-orbit of $u-1$ become
inner edges of level 1 and all edges in the $\sigma$-orbit of
$u+1$ become outer edges of level 1 in $T^*$.

In particular, the number of edges on the cycle of $\overline{T}$
can be reduced, up to derived equivalence. The automorphism
$\sigma^*$ is again given by $f\mapsto \sigma(f)$ for all edges
$f$. But since the number of edges decreases, the rotation
parameter changes.
\end{proposition}
\begin{proof}
We consider the following two-term complex of projective
$A$-modules. For any edge
$z\not\in\{\sigma^j(u)~|~j\in\mathbb{Z}\}$ we let $Q(z)\,:\,0\to
P(z)\to 0$ be the stalk complex in degree 0. By assumption, for
any $j\in\mathbb{Z}$, there are no edges between $\sigma^j(u)$ and
$\sigma^j(u-1)$ and between $\sigma^j(u)$ and $\sigma^j(u+1)$.
Hence, every map from $P(z)$ to $P(\sigma^j(u))$ factors through
$P(\sigma^{j-1}(u-1))$ or $P(\sigma^{j}(u+1))$. Thus we are led to
define the complexes
$$Q(\sigma^j(u))\,:\,0\to P(\sigma^{j-1}(u-1))\oplus
P(\sigma^{j}(u+1))\to P(\sigma^j(u))\to 0,$$ in degrees 0 and
$-1$. Then the direct sum $Q:=\oplus_{z\in\overline{T}} Q_{z}$ is
a tilting complex for $A$. We leave the straightforward
verification to the reader.
\smallskip

We use the same notation as before for edges in the Brauer graph
$\overline{T}$
of $A$. We denote by $\lo u-1\po_1,\ldots,\lo u-1\po_{q_{u-1}}$
the inner edges and by $\lo u-1\po^1,\ldots,\lo u-1\po^{p_{u-1}}$
the outer edges between $u-1$ and $u-2$, respectively. Moreover,
let $\lo u+2\po_1,\ldots,\lo u+2\po_{q_{u+2}}$ be the inner edges
and $\lo u+2\po^1,\ldots,\lo u+2\po^{p_{u+2}}$ the outer edges
between $u+1$ and $u+2$, respectively. Note that any of the cases
$q_{u-1}=0$, $p_{u-1}=0$, $q_{u+2}=0$ or $p_{u+2}=0$ is allowed.

We have to determine the endomorphism ring of $Q$ in the homotopy
category $K^b(A)$
(it is derived equivalent
to $A$ by Rickard's theorem).
For describing the Brauer graph of the endomorphism ring $\End(Q)$,
we denote the
'new' two-term complexes by $\tilde{\sigma}^j(u):=Q(\sigma^j(u))$.
The edges corresponding to the stalk complexes are denoted as before.

We claim that $\End(Q)$ is of the form $\Omega^{(1)}\lo
T^*,\sigma^*,\theta\po$ where the Brauer graph $T^*$ is given as
follows. Note that the Brauer graph will only change around the
new complexes $Q(\sigma^j(u))$ (since the others are stalk
complexes). It suffices to describe the structure around one of
the new complexes, say $\tilde{\sigma}(u)$.
\smallskip

{\em Claim: The ordering $p(T^*,\sigma^*)$ of edges around the one
endpoint of $\tilde{\sigma}(u)$ is given by
$$\ldots,\tilde{\sigma}(u),
u-1,\lo u-1\po_1,\ldots,\lo
u-1\po_{q_{u-1}},(\tilde{\sigma}^0(u-3)),\hskip2cm
$$
$$\hskip4cm u-2,\lo
u-1\po^1,\ldots,\lo u-1\po^{p_{u-1}},\tilde{\sigma}^2(u),\ldots$$
where $\tilde{\sigma}^0(u-3)$ occurs exactly if $u-3$ is contained
in the $\sigma$-orbit of $u$. In this case,
$\tilde{\sigma}^0(u-3)$ is on the cycle of $T^*$; otherwise, $u-2$
is on the cycle.

The order $p(T^*,\sigma^*)$
around the other endpoint of $\tilde{\sigma}(u)$
is given by
$$\ldots,\tilde{\sigma}(u),
\sigma(u+1),\sigma(\lo u+2\po^1),\ldots,\sigma(\lo
u+2\po^{p_{u+2}}),
(\tilde{\sigma}^3(u+3)),
$$
$$\hskip2cm\sigma^2(u+2),\sigma^2(\lo
u+2\po_1),\ldots, \sigma^2(\lo
u+2\po_{q_{u+2}}),\tilde{\sigma}^2(u),\ldots
$$
where
$\tilde{\sigma}^3(u+3)$ occurs precisely if $u-3$ is in the
$\sigma$-orbit of $u$. In this case, $\tilde{\sigma}^3(u+3)$ is on
the cycle of $T^*$; otherwise $\sigma^2(u+2)$ is on the cycle.

The adjacent edges on the cycle of $T^*$ are
$\ldots,\tilde{\sigma}^0(u-3),\tilde{\sigma}(u),\tilde{\sigma}^2(u+3),
\ldots$ if $u-3$ is in the $\sigma$-orbit of $u$; otherwise
the edges on the cycle are given by
$\ldots,u-2,\tilde{\sigma}(u),\sigma(u+2),\ldots$.
}

\smallskip
Note that the crucial aspect of this description is that the edges
$\sigma^j(u-1)$ and $\sigma^j(u+1)$, which were on the cycle in
$\overline{T}$, become inner and outer edges of level 1 in the
Brauer graph for $\End(Q)$. In particular, the cycle for $\End(Q)$
has fewer edges than the one for $A$.
\smallskip

{\em Proof of the claim.} We have to define maps of complexes between
summands of $Q$ corresponding to the consecutive edges as given in
the claim.

We define a map $Q(\sigma(u))\to Q(u-1)$ by the projection onto
$P(u-1)$ in degree 0.

The maps between the following stalk complexes are defined by
$[u-1,\lo u-1\po_1]:Q(u-1)\to Q\lo(u-1)_{1}\po$, and, for
$i=1,\ldots,q_{u-1}-1$ by $[\lo u-1\po_i,\lo u-1\po_{i+1}]:Q(\lo
u-1\po_i)\to Q(\lo u-1\po_{i+1}).$

Moreover, there is the map $[\lo u-1\po_{q_{u-1}},u-2]: Q(\lo
u-1\po_{q_{u-1}})\to Q(u-2)$. If $u-3$ is contained in the
$\sigma$-orbit of $u$ then this map factors through the complex
$$Q(u-3)\,:\,0\to P(\sigma^{-1}(u-4))\oplus P(u-2)\to P(u-3)\to 0$$
(indeed a map of complexes since $\lo u-1\po_{q_{u-1}},u-2,u-3$
are not consecutive edges). This factorization then makes $u-2$ an
edge of level 1 in $T^*$. Otherwise, $u-2$ is on the cycle of
$T^*$, between $Q(\lo u-1\po_{q_{u-1}})$ and $Q(\lo u-1\po^1)$.

The maps $[\lo u-1\po^i,\lo u-1\po^{i+1}]:Q(\lo u-1\po^i)\to Q(\lo
u-1\po^{i+1})$ are as before for $A$.

Finally, we have a map $Q(\lo u-1\po^{p_{u-1}})\to Q(\sigma^2(u))$
given by $([\lo u-1\po^{p_{u-1}},\sigma(u-1)],0)$ in degree 0.
Note that this indeed is a map of complexes since $\lo
u-1\po^{p_{u-1}},\sigma(u)-1,\sigma^2(u)$ are not consecutive
edges.
\smallskip

Next, we consider the edges around the other endpoint of the edge
$\tilde{\sigma}(u)=Q(\sigma(u))$.

Let the map $Q(\sigma(u))\to Q(\sigma(u+1))$ be defined by the
projection onto $P(\sigma(u+1))$ in degree 0. Moreover, we define
the following maps between stalk complexes. We set
$[\sigma(u+1),\sigma(\lo u+2\po^1)]:Q(\sigma(u+1)) \to
Q(\sigma(\lo u+2\po^1))$ and, for all $i=1,\ldots,p_{u+2}-1$, we
set $[\sigma(\lo u+2\po^i),\sigma(\lo u+2\po^{i+1})]:Q(\sigma(\lo
u+2\po^i)) \to Q(\sigma(\lo u+2\po^{i+1}))$. As in $A$, there is a
map $[\sigma(\lo u+2\po^{p_{u+2}}),\sigma^2(u+2)]:Q(\sigma(\lo
u+2\po^{p_{u+2}}) \to Q(\sigma^2(u+2))$.

If $u-3$ is in the $\sigma$-orbit of $u$, then also
$\sigma^3(u+3)$ is in the $\sigma$-orbit of $u$. Then the above
map $[\sigma(\lo u+2\po^{p_{u+2}}),\sigma^2(u+2)]$ factors through
(the first summand of) the complex
$$Q(\sigma^3(u+3))\,:\,0\to P(\sigma^2(u+2))\oplus P(\sigma^3(u+4))
\to P(\sigma^3(u+3))\to 0.$$ Note that this indeed gives a map of
complexes since $\sigma(\lo
u+2\po^{p_{u+2}}),\sigma^2(u+2),\sigma^3(u+3)$ are not consecutive
edges. In this case, $Q(\sigma^2(u+3))$ is the next edge on the
cycle adjacent to $Q(\sigma(u))$. Otherwise, $\sigma(u+2)$ is the
next edge on the cycle.

The maps between the successive inner edges come from stalk
complexes and are given exactly as in $A$, that is, we consider
$[\sigma^2(u+2),\sigma^2(\lo u+2\po_1)]$, and, for
$i=1,\ldots,q_{u+2}-1$, we have $[\sigma^2(\lo
u+2\po_i),\sigma^2(\lo u+2\po_{i+1})]$.

Finally, we need a map $Q(\sigma^2(\lo u+2\po_{q_{u+2}}))\to
Q(\sigma^2(u))$ which we define by
$$(0,[\sigma^2(\lo u+2\po_{q_{u+2}}),\sigma^2(u+1)]):
P(\sigma^2(\lo u+2\po_{q_{u+2}}))
\to P(\sigma(u-1))\oplus P(\sigma^2(u+1))$$ in degree 0. Note that
this is a map of complexes since $\sigma^2(\lo
u+2\po_{q_{u+2}}),\sigma^2(u+1),\sigma^2(u)$ are not consecutive
edges.

Clearly, all maps defined above are not homotopic to zero.

In order to complete the proof we have to check that they satisfy
the defining relations of $\Omega^{(1)}\lo
T^*,\sigma^*,\theta\po$, up to homotopy.
\smallskip

(1) {\em Non-consecutive edges.} We only have to check sequences
of edges where one of the new complexes occurs. (For the stalk complexes
it directly follows from the same relations in $A$ that non-consecutive
edges give zero relations.)

The edges $\lo u-1\po^{p_{u-1}},\tilde{\sigma}^2(u),\sigma^2(u+1)$
are not consecutive in the order $p(T^*,\sigma^*)$. The
corresponding composition of maps of complexes, as defined above,
is given by first mapping to the first summand in degree 0, then
projecting onto the second summand. Clearly, this composition is
zero. Similarly, the inner edges $\sigma(\lo
u+2\po_{q_{u+2}}),\tilde{\sigma}(u), u-1$ are not consecutive. The
corresponding composition of maps of complexes first maps into the
second summand in degree 0 and then projects onto the first
summand, which is obviously zero.
\smallskip

(2) {\em Socle relations.} From the above definitions, it is easy
to check that for all edges $f$ corresponding to stalk complexes
the socle map $Q(f)\to Q(\sigma^*(f))$ is just given by the socle
map $P(f)\to P(\sigma(f))$. Then the socle relations immediately
follow from the socle relations of $A$. (Clearly, around the
endpoint of $f$ not attached to the cycle, we also define $Q(f)\to
Q(\sigma^*(f))$ by the same socle map.) So it only remains to
prove the socle relations for the new complexes $Q(\sigma^j(u))$.
Again, it suffices to consider
$$Q(\sigma(u))\,:\,0\to P(u-1)\oplus P(\sigma(u+1))\to P(\sigma(u))
\to 0.$$
Around the one endpoint
of $\tilde{\sigma}(u)$ we obtain the map $C_1^*:Q(\sigma(u))\to
Q(\sigma^2(u))$
which is given by the socle map $P(u-1)\to P(\sigma(u-1))$
on the first summand
in degree 0. Around the other endpoint we obtain the map $C_2^*$
given by the socle map $P(\sigma(u+1))\to P(\sigma^2(u+1))$
on the second summand in degree 0.
Both maps are zero in degree $-1$.

In order to satisfy the relations of $\Omega^{(1)}\lo
T^*,\sigma^*,\theta\po$ these maps have to be homotopic
(up to the scalar $\theta$ in case of the distinguished
edge involved).
By assumption, $u-1$ and $u+1$ are not in the $\sigma$-orbit
of the distinguished edge 1. In particular, the distinguished
edge remains on the cycle of $\End(Q)$. For checking the
relation involving the scalar $\theta$, we can assume
that $\sigma(u)=1$ is the distinguished edge (otherwise, set
$\theta=1$ in the sequel).
Then the map $C_1^*-\theta C_2^*$ is homotopic to zero
via the homotopy map
$$([\sigma(u),\sigma(u-1)],-\theta [\sigma(u),\sigma^2(u+1)]):
P(\sigma(u))\to P(\sigma(u-1))\oplus P(\sigma^2(u+1)).$$
(For checking the details, one also repeatedly uses the fact that
non-consecutive edges in the order of $\overline{T}$
give zero relations.)

This completes the proof of Proposition \ref{redstep3}.
\end{proof}


\section{\hspace{-2.5ex}. Derived normal forms for the algebras $\bf \Omega^{\lo
1\po}$} \label{derived}


In this section we indicate how one can combine the reduction
steps in the previous sections to obtain a proof of the fact that
any algebra $\Omega^{(1)}\lo T,\sigma,\theta\po$ is derived
equivalent to one of the normal forms $\Lambda(p,q,k,s,\theta)$.

We can actually prove the following more precise result.

\begin{tw} \label{normalform}
Let $A=\Omega^{(1)}\lo T,\sigma_s,\theta\po$ be an algebra as
defined in Section \ref{morita}, where the cycle of $T$ has length
$k$. Denote by $o(T)$ the (total) number of outer edges in $T$,
and by $i(T)$ the (total) number of inner edges in $T$. Then $A$
is derived equivalent to the normal form
$\Lambda(p',q',k',s',\theta)$, where the parameters $p',q',k',s'$
are given as follows

$p'= \frac{\gcd(s,k)}{k}\cdot o(T) + \frac{\gcd(s,k)-1}{2}$,

$q'= \frac{\gcd(s,k)}{k}\cdot i(T) + \frac{\gcd(s,k)-1}{2}$,

$k'= \frac{k}{\gcd(s,k)}$,

$s'= \frac{s}{\gcd(s,k)}.$

\end{tw}

\noindent {\bf Remark.} It might not be obvious that the
parameters $p'$ and $q'$ as given in the theorem are actually
natural numbers. Let us recall the assumptions on the parameters
$s,k$ in the definition of $\Omega^{(1)}\lo T,\sigma_s,\theta\po$,
where $T$ has a cycle of length $k\geqslant 2$ (see Section
\ref{morita}).

The assumption for the one-parametric algebras has been that
$1\leqslant s\leqslant k-1$ and $\gcd(s+2,k)=1$. In particular,
this implies that $\gcd(s,k)$ is odd, that is, the second summand
in the formulae for $p'$ and $q'$ is an integer.

The automorphism $\sigma_s$ is just rotation of $T$ by $s$ steps
on the cycle of length $k$. Hence, the order of $\sigma_s$ (as an
automorphism) is $\frac{k}{\gcd(s,k)}$. In particular,
the $\sigma_s$-orbit of any edge (or any vertex) consists of
exactly $\frac{k}{\gcd(s,k)}$ edges (or vertices). The group
generated by $\sigma_s$ acts on the outer (and inner) edges of
$T$. Hence, the orbit length $\frac{k}{\gcd(s,k)}$ divides the
numbers $o(T)$
and $i(T)$. In particular, the formulae for $p'$ and $q'$
actually give natural numbers.

\medskip

\noindent {\bf Proof of Theorem \ref{normalform}.} Let
$A=\Omega^{(1)}\lo T,\sigma_s,\theta\po$ be an algebra as defined
in Section \ref{morita}, with a cycle of length $k$. For the proof
of the theorem we need to have a very detailed look at the
reduction steps of the previous sections. In particular, we need
to know how the parameters $s,k,o(T),i(T)$ occurring on the right
hand side are changed in each step.
\smallskip

The first reduction step Proposition \ref{redstep1} states that,
up to derived equivalence, the inner and outer Brauer trees can be
made into stars, that is, all edges have level at most 1. The
cycle of $T$ remains unchanged, as does the automorphism
$\sigma_s$. So, the parameters considered in the statement of the
theorem are not affected by this first reduction step.
\smallskip

The second reduction step Proposition \ref{redstepmoving} states
that $A$ is derived equivalent to an algebra
$\Omega^{(1)}\lo\overline{T},\sigma_s,\theta\po$ where all edges
of level 1 are attached to some vertex in the $\sigma_s$-orbit of
{\em one} fixed vertex $v$. Recall that in this second reduction
step inner (resp. outer) edges remain inner (resp. outer) edges.
In particular, the numbers $o(T)$ and $i(T)$ are not changed by
this reduction step. Note also that again the cycle and the
automorphism $\sigma_s$ are not changed.

We can determine the number of outer and inner edges attached to
each vertex in the $\sigma_s$-orbit of the fixed vertex $v$ of
$\overline{T}$. There are $o(T)$ outer edges in total; they are
equally distributed to the $\frac{k}{\gcd(s,k)}$ vertices in the
$\sigma_s$-orbit. Hence, at each such a vertex we have
$\frac{\gcd(s,k)}{k}\cdot o(T)$ outer and
$\frac{\gcd(s,k)}{k}\cdot i(T)$ inner edges attached. Note that
these terms are exactly the first summands in the formulae for
$p'$ and $q'$.
\smallskip

In the third reduction step Proposition \ref{redstep3} the length
$k$ of the cycle and the rotation parameter $s$ are changed. Recall the
situation of Proposition \ref{redstep3}. We assume that there are
two consecutive vertices on the cycle of $\overline{T}$ with no
edges of level 1 attached to them (i.e., these vertices are not in the
$\sigma_s$-orbit of the distinguished vertex $v$). The adjacent
edges were denoted $u-1,u,u+1$. Then, up to derived equivalence,
we can change the Brauer graph $\overline{T}$ to the Brauer graph
$T^*$ where all edges in the $\sigma_s$-orbit of $u-1$ become
inner edges of level 1, and all edges in the $\sigma_s$-orbit of
$u+1$ become outer edges of level 1. By applying again our second
reduction step (moving edges), if necessary, we can assume that
the new inner and outer edges are also attached to vertices in the
$\sigma_s$-orbit of the fixed vertex $v$ (recall that this second
reduction step does not affect the parameters).

We have to determine how the parameters $s$ and $k$ are changed.
The $\frac{k}{\gcd(s,k)}$ vertices in the $\sigma_s$-orbit of the
fixed vertex $v$ of $\overline{T}$ are at equal distance apart
from each other, where this distance is $\gcd(s,k)$ (the length of
the cycle divided by the length of the orbit). By Proposition
\ref{redstep3}, in each of the
intervals between two vertices of the $\sigma_s$-orbit of $v$, two
edges on the cycle become edges of level 1. Hence, the length of
the cycle of $T^*$ is given by $k_1:=k-2\cdot
\frac{k}{\gcd(s,k)}$. Concerning the rotation parameter, there are
exactly $\frac{s}{\gcd(s,k)}$ such intervals between the fixed
vertex $v$ and $\sigma_s(v)$. Hence, the automorphism $\sigma^*$
rotates $T^*$ by $s_1:=s-2\cdot\frac{s}{\gcd(s,k)}$.

We can now apply the third reduction step inductively. We set
$k_0=k$ and $s_0=s$. After $i$ iterations we obtain a Brauer graph
$T^*_{(i)}$ with a cycle of length $k_i:=k_{i-1}-2\cdot
\frac{k_{i-1}}{\gcd(s_{i-1},k_{i-1})}$ and rotation parameter
$s_i:=s_{i-1}-2\cdot \frac{s_{i-1}}{\gcd(s_{i-1},k_{i-1})}$.

By our assumptions we know that the length $\gcd(s,k)$ of the
intervals is odd (see the preceding remark). So, after
$\frac{\gcd(s,k)-1}{2}$ iterations we get a Brauer graph for which
the intervals are of length 1, that is, the action of the rotation
is transitive. The corresponding algebra is one of the normal
forms $\Lambda(p',q',k',s',\theta)$, that is, we completed the proof of
the implication $(1)\Rightarrow (2)$ in our main result Theorem
\ref{mainthm}.

It remains to determine the parameters for the normal forms.

After the second reduction step, in the Brauer graph
$\overline{T}$ we have $\frac{\gcd(s,k)}{k}\cdot o(T)$ outer edges
and $\frac{\gcd(s,k)}{k}\cdot i(T)$ inner edges attached to each
vertex in the $\sigma_s$-orbit of the fixed vertex $v$.

After $\frac{\gcd(s,k)-1}{2}$ iterations of the third reduction
step we get at each of these vertices exactly
$\frac{\gcd(s,k)-1}{2}$ additional outer and inner edges. This
proves the formulae for $p'$ and $q'$ given in Theorem
\ref{normalform}.
\smallskip

Closed formulae for the number $k'$ of edges on the cycle and the
rotation parameter $s'$ are provided by the following observation
on the numbers $k_i$ and $s_i$ which were defined inductively
above.
\smallskip

\noindent {\em Claim: For any $i=0,\ldots,\frac{\gcd(s,k)-1}{2}$,
we have $k_i= k\cdot\frac{\gcd(s,k)-2i}{\gcd(s,k)}$ and $s_i=
s\cdot\frac{\gcd(s,k)-2i}{\gcd(s,k)}$.}
\smallskip

\noindent {\em Proof of the claim.} We first prove by induction
that $\gcd(s_i,k_i)=\gcd(s,k)-2i$.

By definition, this holds for $i=0$. Let $i>0$. Then
$$\begin{array}{rcl}
\gcd(s_i,k_i) & = &
\gcd(s_{i-1}(1-\frac{2}{\gcd(s_{i-1},k_{i-1})}),
k_{i-1}(1-\frac{2}{\gcd(s_{i-1},k_{i-1})})) \\
 & = & (1-\frac{2}{\gcd(s_{i-1},k_{i-1})}) \gcd(s_{i-1},k_{i-1}) \\
 & = & \gcd(s_{i-1},k_{i-1})-2 \\
 & = & \gcd(s,k)-2i,
\end{array}
$$
by induction.

Using this formula for the $\gcd(s_{i},k_i)$, the formulae for
$k_i$ and $s_i$ now follow easily by induction. This completes the
proof of the claim.
\smallskip

Recall that we got to our normal form
$\Lambda(p',q',k',s',\theta)$ after exactly
$i=\frac{\gcd(s,k)-1}{2}$ iterations of our third reduction step.
Using the above claim, the length of the cycle is easily computed
to be $k'=k_{\frac{\gcd(s,k)-1}{2}} = \frac{k}{\gcd(s,k)}$ and the
rotation parameter is equal to $s'=s_{\frac{\gcd(s,k)-1}{2}} =
\frac{s}{\gcd(s,k)}$. This completes the proof of Theorem
\ref{normalform}. $\hfill\Box$

\medskip

The following example illustrates the above considerations.
\begin{exa}\rm
Let $T$ be the Brauer graph of the form
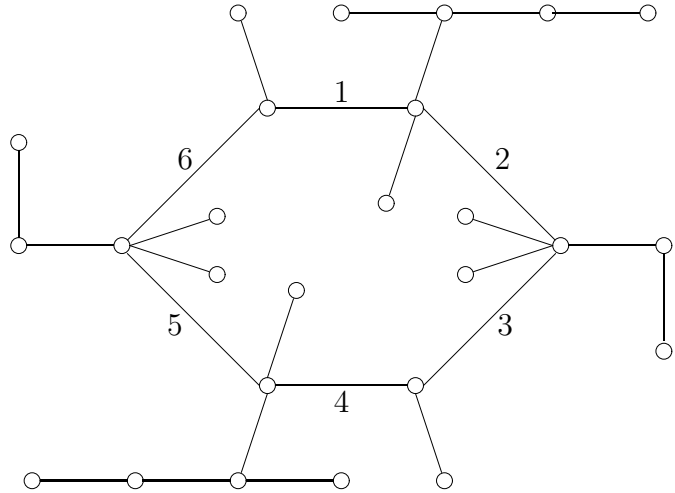
\begin{center}
\unitlength=1.000000pt
\begin{picture}(258.00,200.00)(-78.00,-95.00)
\put(81.00,59.00){\circle{6.00}} \put(-30.00,7.00){\circle{6.00}}
\put(136.00,7.00){\circle{6.00}} \put(25.00,59.00){\circle{6.00}}
\put(78.00,59.00){\line(-1,0){50.00}}
\put(84.00,59.00){\line(1,-1){50.00}}
\put(22.00,59.00){\line(-1,-1){50.00}}
\put(-28.00,4.00){\line(1,-1){50.00}}
\put(134.00,4.00){\line(-1,-1){50.00}}
\put(81.00,-46.00){\circle{6.00}}
\put(25.00,-46.00){\circle{6.00}}
\put(28.00,-46.00){\line(1,0){50.00}}

\put(95,95){\line(1,0){33}} \put(81,62){\line(1,3){10}}
\put(131,95){\circle{6.00}} \put(92,95){\circle{6.00}}
\put(133,95){\line(1,0){33}} \put(169,95){\circle{6.00}}
\put(89,95){\line(-1,0){33}} \put(53,95){\circle{6.00}}
\put(81,56){\line(-1,-3){10}} \put(70,23){\circle{6.00}}

\put(25,62){\line(-1,3){10}} \put(14,95){\circle{6.00}}

\put(-33,7){\line(-1,0){33}} \put(-69,7){\circle{6.00}}
\put(-69,10){\line(0,1){33}} \put(-69,46){\circle{6.00}}
\put(-27,7){\line(3,1){30}} \put(-27,7){\line(3,-1){30}}
\put(6,18){\circle{6.00}} \put(6,-4){\circle{6.00}}

\put(133,7){\line(-3,1){30}} \put(133,7){\line(-3,-1){30}}
\put(175,7){\circle{6.00}} \put(175,-33){\circle{6.00}}
\put(139,7){\line(1,0){33}} \put(175,4){\line(0,-1){33}}
\put(100,18){\circle{6.00}} \put(100,-4){\circle{6.00}}

\put(81,-49){\line(1,-3){10}} \put(92,-82){\circle{6.00}}

\put(11,-82){\line(-1,0){33}} \put(25,-49){\line(-1,-3){10}}
\put(-25,-82){\circle{6.00}} \put(36,-10){\circle{6.00}}
\put(-28,-82){\line(-1,0){33}} \put(-64,-82){\circle{6.00}}
\put(17,-82){\line(1,0){33}} \put(53,-82){\circle{6.00}}
\put(25,-43){\line(1,3){10}} \put(14,-82){\circle{6.00}}

\put(53.00,65.00){\makebox(0,0){$1$}}
\put(114.00,40.00){\makebox(0,0){$2$}}
\put(115.00,-23.00){\makebox(0,0){$3$}}
\put(-10.00,-23.00){\makebox(0,0){$5$}}
\put(53.00,-52.00){\makebox(0,0){$4$}}
\put(-6.00,40.00){\makebox(0,0){$6$}}
\put(185.00,-90.00){\makebox(0,0){.}}
\end{picture}
\end{center}
and $s=3$. Then, for $k=6$, we have $\gcd\lo s\! +\! 2,k\po=1$ but
$\gcd\lo s,k\po=3\neq 1$. Furthermore, using the notation of
Theorem \ref{normalform}, we have $o(T)=14$, $i(T)=6$ and thus
obtain $p^{\prime}=8$, $q^{\prime}=4$, $k^{\prime}=2$, and
$s^{\prime}=1$. Therefore, for $\theta\in K\backslash\lk 0\pk$,
the algebra $\Omega^{\lo 1\po}\!\lo T,\sigma_{3},\theta\po$ is
derived equivalent to $\Lambda\lo 8,4,2,1,\theta\po=\Omega^{\lo
1\po}\!\lo T_{8,4,2},\sigma_{1},\theta\po$, where the Brauer graph
$T_{8,4,2}$ is of the form
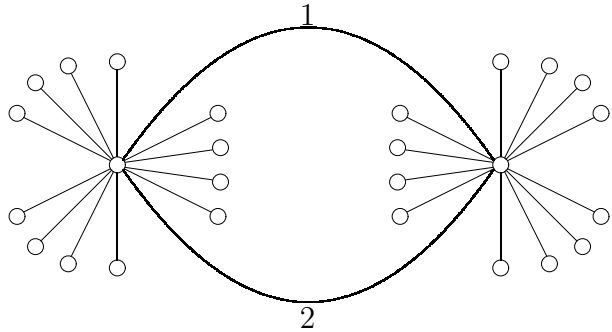
\begin{center}
\unitlength=1.000000pt
\begin{picture}(238.00,130.00)(-20.00,-15.00)
\bezier{1000}(30,50)(100,-50)(170,50)
\bezier{1000}(30,54)(100,154)(170,54)
\put(28.00,52.00){\circle{6.00}} \put(173.00,52.00){\circle{6.00}}
\put(100.00,109.00){\makebox(0,0){$1$}}
\put(100.00,-6.00){\makebox(0,0){$2$}}

\put(173,55){\line(0,1){33}} \put(173,49){\line(0,-1){33}}
\put(175,54){\line(1,1){27}} \put(175,50){\line(1,-1){27}}
\put(175.5,53.5){\line(2,1){33}} \put(175.5,50.5){\line(2,-1){33}}
\put(174.5,54.5){\line(1,2){16}} \put(174.5,49.5){\line(1,-2){16}}

\put(173.00,91.00){\circle{6.00}}
\put(173.00,13.00){\circle{6.00}}
\put(191.50,89.50){\circle{6.00}}
\put(191.50,14.50){\circle{6.00}}
\put(204.00,83.00){\circle{6.00}}
\put(204.00,21.00){\circle{6.00}}
\put(211.00,71.50){\circle{6.00}}
\put(211.00,32.50){\circle{6.00}}

\put(170.5,53.5){\line(-2,1){33}}
\put(170.5,50.5){\line(-2,-1){33}} \put(170,53){\line(-7,1){33}}
\put(170,51){\line(-7,-1){33}}

\put(135.00,71.50){\circle{6.00}}
\put(135.00,32.50){\circle{6.00}}
\put(134.00,58.50){\circle{6.00}}
\put(134.00,45.50){\circle{6.00}}

\put(28,55){\line(0,1){33}} \put(28,49){\line(0,-1){33}}
\put(26,54){\line(-1,1){27}} \put(26,50){\line(-1,-1){27}}
\put(25.5,53.5){\line(-2,1){33}} \put(25.5,50.5){\line(-2,-1){33}}
\put(26.5,54.5){\line(-1,2){16}} \put(26.5,49.5){\line(-1,-2){16}}

\put(28.00,91.00){\circle{6.00}} \put(28.00,13.00){\circle{6.00}}
\put(9.50,89.50){\circle{6.00}} \put(9.50,14.50){\circle{6.00}}
\put(-3.00,83.00){\circle{6.00}} \put(-3.00,21.00){\circle{6.00}}
\put(-10.00,71.50){\circle{6.00}}
\put(-10.00,32.50){\circle{6.00}}

\put(30.5,53.5){\line(2,1){33}} \put(30.5,50.5){\line(2,-1){33}}
\put(31,53){\line(7,1){33}} \put(31,51){\line(7,-1){33}}

\put(66.00,71.50){\circle{6.00}} \put(66.00,32.50){\circle{6.00}}
\put(67.00,58.50){\circle{6.00}} \put(67.00,45.50){\circle{6.00}}
\end{picture}
\end{center}
\end{exa}


\section{\hspace{-2.5ex}. Derived normal forms for the algebras $\bf\Omega^{\lo
2\po}$} \label{dnderived}


The aim of this section is to indicate a proof of the following
derived equivalence classification for the algebras
$\Omega^{(2)}(T,v_1,v_2)$ as defined in Section \ref{morita}.

\begin{proposition} \label{normaldn}
Any algebra $\Omega^{(2)}(T,v_1,v_2)$ is derived equivalent to the
normal form $\Gamma^*(n)$ where $n+2$ is the number of simple
modules of $\Omega^{(2)}(T,v_1,v_2)$.

In particular, two algebras in the family $\Omega^{(2)}$ are
derived equivalent if and only if they have the same number
of simple modules.
\end{proposition}

The algebras $\Omega^{(2)}(T,v_1,v_2)$ are very similar to weakly
symmetric algebras of Euclidean type studied in \cite{bib:BoHoSk}.
The difference only consists of a few relations which interchange
the socles of two projective indecomposable modules making
$\Omega^{(2)}(T,v_1,v_2)$ selfinjective but not weakly symmetric.
Apart from that, all remaining relations are the same.

In \cite[Section 4]{bib:BoHoSk} we established the derived normal
forms $\Gamma\lo n\po$ for the weakly symmetric
algebras $\Gamma^{(2)}(T,v_1,v_2)$.
It turned out that the normal forms $\Gamma\lo n\po$ are very
similar to our new normal forms $\Gamma^{*}(n)$, up to a twist in
the relations, exactly corresponding to the twist in the relations
for going from $\Gamma^{(2)}(T,v_1,v_2)$ to
$\Omega^{(2)}(T,v_1,v_2)$. In fact, the algebra $\Gamma^{*}\lo
n\po$ is given by the following quiver and relations
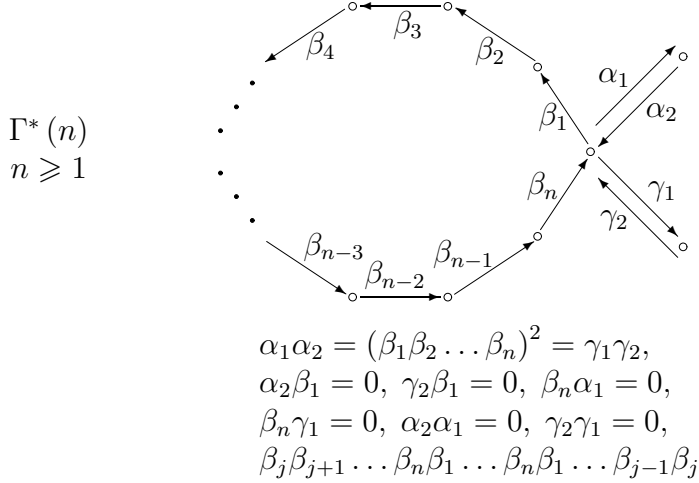
\begin{center}
\unitlength=1.0pt
\begin{picture}(441,196)(-140,-145)
 \put(99,-12){\vector(-2,3){17.8}}
 \put(79,19.5){\vector(-3,2){30}}
 \put(43,40){\vector(-1,0){30}}
 \put(7,39){\vector(-3,-2){30}}

 \put(81,-44.5){\vector(2,3){17.8}}
 \put(48.5,-69){\vector(3,2){30}}
 \put(13,-70){\vector(1,0){30}}
 \put(-22.5,-49){\vector(3,-2){30}}

 \put(80,17){\circle{3}}
 \put(46,40){\circle{3}}
 \put(10,40){\circle{3}}

 \put(80,-47){\circle{3}}
 \put(46,-70){\circle{3}}
 \put(10,-70){\circle{3}}

 \put(100,-15){\circle{3}}

 \put(102,-5){\vector(1,1){30}}
 \put(133,17){\vector(-1,-1){30}}
 \put(103,-17){\vector(1,-1){30}}
 \put(133,-55){\vector(-1,1){30}}
 \put(135,21){\circle{3}}
 \put(135,-51){\circle{3}}
 \put(109,13){\makebox(0,0){\normalsize $\alpha_{1}$}}
 \put(127,0){\makebox(0,0){\normalsize $\alpha_{2}$}}
 \put(127,-30){\makebox(0,0){\normalsize $\gamma_{1}$}}
 \put(109,-41){\makebox(0,0){\normalsize $\gamma_{2}$}}

 \multiput(-28,11)(-6,-9){3}{\circle*{1.5}}
 \multiput(-40,-23)(6,-9){3}{\circle*{1.5}}

 \put(-105,-15){\makebox(0,0){\normalsize $\begin{array}{c} \Gamma^{*}\lo n\po\\
 n \geqslant 1\end{array}$}}

 \put(86,-2){\makebox(0,0){\normalsize $\beta_{1}$}}
 \put(62,23){\makebox(0,0){\normalsize $\beta_{2}$}}
 \put(31,33){\makebox(0,0){\normalsize $\beta_{3}$}}
 \put(-2,25){\makebox(0,0){\normalsize $\beta_{4}$}}

 \put(82,-28){\makebox(0,0){\normalsize $\beta_{n}$}}
 \put(52,-54){\makebox(0,0){\normalsize $\beta_{n-1}$}}
 \put(26,-62){\makebox(0,0){\normalsize $\beta_{n-2}$}}
 \put(3,-51){\makebox(0,0){\normalsize $\beta_{n-3}$}}

 \put(100,-110){\makebox(0,0){\normalsize $\begin{array}{l}
 \alpha_{1}\alpha_{2}=\lo\beta_{1}\beta_{2}\ldots\beta_{n}\po^{2}=\gamma_{1}\gamma_{2},\\
 \alpha_{2}\beta_{1}=0,\;\gamma_{2}\beta_{1}=0,\;\beta_{n}\alpha_{1}=0,\\
 \beta_{n}\gamma_{1}=0,\;\alpha_{2}\alpha_{1}=0,\;\gamma_{2}\gamma_{1}=0,\\
 \beta_{j}\beta_{j+1}\ldots\beta_{n}\beta_{1}\ldots\beta_{n}\beta_{1}\ldots\beta_{j-1}\beta_{j}=0,\;
 2\leqslant j \leqslant n,
 \end{array}$}}
\end{picture}
\end{center}
and $\Gamma\lo n\po$ is obtained from $\Gamma^{*}\lo n\po$ by
replacing the relations $\alpha_{2}\alpha_{1}=0$ and
$\gamma_{2}\gamma_{1}=0$ by the relations $\alpha_{2}\gamma_{1}=0$
and $\gamma_{2}\alpha_{1}=0$.

Because of the similarities described above it is therefore not
surprising that the proof for the derived equivalence
classification of the weakly symmetric algebras carries over
almost verbatim to our new situation. We refrain from repeating
the proof here, and refer for details to our previous paper
\cite[Section 4]{bib:BoHoSk}.


\section{\hspace{-2.5ex}. Auslander-Reiten sequences} \label{ARseq}


In order to prove Theorem \ref{mainthm}, we need more detailed
information on the Auslander-Reiten quivers of standard
one-parametric selfinjective algebras and the action of Heller's
syzygy operator on their stable module categories.

Let $A$ be a selfinjective algebra. We denote by $\Gamma_{A}$ the
Auslander-Reiten quiver of $A$ and by $\tau_{A}$ and
$\tau_{A}^{-1}$ the Auslander-Reiten translations $\dualnosc\slad$
and $\slad\dualnosc$, respectively. We will identify an
indecomposable module $M$ from $\md A$ with the vertex $\lkw
M\pkw$ of $\Gamma_{A}$ corresponding to it. Further, we denote by
$\Gamma_{A}^{s}$ the stable Auslander-Reiten quiver of $A$,
obtained from $\Gamma_{A}$ by removing the projective modules and
the arrows attached to them. Recall also that for an
indecomposable projective $A$-module $P$ we have in $\md A$ an
Auslander-Reiten sequence of the form
\begin{displaymath}
0\longrightarrow\rad P\longrightarrow\rad P/S\oplus P
\longrightarrow P/S\longrightarrow 0,
\end{displaymath}
where $S$ is the socle of $P$ \cite[Proposition
V.5.5]{bib:AuReSm}. Finally, we denote by $\Omega_{A}$ Heller's
syzygy operator on $\md A$ which assigns to any module $M$ the
kernel of its projective cover $P\lo M\po\longrightarrow M$.

The first result of this section describes the shape of the stable
Auslander-Reiten quivers of the algebras $\Gamma^{*}\lo n\po$.

\begin{proposition} \label{structard}
$\Gamma_{\Gamma^{*}\lo n\po}^{s}$ consists of one Euclidean
component $\mathbbm{Z}\widetilde{\mathbbm{D}}_{2n+3}$, two stable
tubes of rank $2$, one stable tube of rank $2n\!+\!1$, and a
family of stable tubes of rank $1$ indexed by $K\backslash\lk
0\pk$.
\end{proposition}
\begin{proof}
It follows from \cite[Section $2$]{bib:BoSk3} and
\cite[Proposition 5.3]{bib:BoHoSk} that $\Gamma^{*}\lo n\po\cong
\widehat{C}\lo n\po/\lo\rho_{n}\psi_{n}\po$, where $C\lo n\po$ is
the bound quiver algebra $K\Delta\lo n\po/I\lo n\po$ given by the
quiver
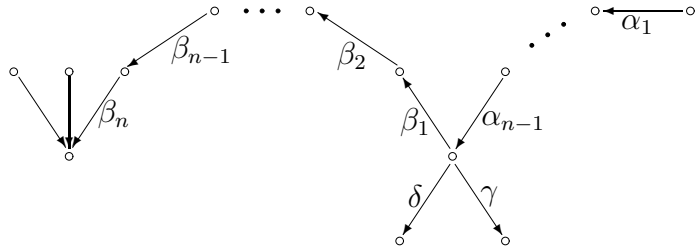
\begin{center}
\unitlength=1.0pt
\begin{picture}(260,110)(-68,-55)
 \put(99,-12){\vector(-2,3){17.8}}
 \put(119,14.5){\vector(-2,-3){17.8}}
 \put(79,19.5){\vector(-3,2){30}}
 \put(-26,15){\vector(-2,-3){17.8}}
 \put(187,40){\vector(-1,0){30}}
 \put(7,39){\vector(-3,-2){30}}
 \put(101,-18){\vector(2,-3){17.8}}
 \put(99,-18){\vector(-2,-3){17.8}}
 \put(-64,15){\vector(2,-3){17.8}}
 \put(-45,15){\vector(0,-1){27}}

 \put(80,17){\circle{3}}
 \put(120,17){\circle{3}}
 \put(46,40){\circle{3}}
 \put(154,40){\circle{3}}
 \put(10,40){\circle{3}}
 \put(190,40){\circle{3}}
 \put(-24,17){\circle{3}}

 \put(-45,-15){\circle{3}}
 \put(100,-15){\circle{3}}

 \put(80,-47){\circle{3}}
 \put(120,-47){\circle{3}}

 \put(-45,17){\circle{3}}
 \put(-66,17){\circle{3}}

 \multiput(34,40)(-6,0){3}{\circle*{1.5}}
 \multiput(130,27)(6,4){3}{\circle*{1.5}}

 \put(86,-2){\makebox(0,0){\normalsize $\beta_{1}$}}
 \put(123,-2){\makebox(0,0){\normalsize $\alpha_{n-1}$}}
 \put(62,23){\makebox(0,0){\normalsize $\beta_{2}$}}
 \put(-28,1){\makebox(0,0){\normalsize $\beta_{n}$}}
 \put(170,35){\makebox(0,0){\normalsize $\alpha_{1}$}}
 \put(5,25){\makebox(0,0){\normalsize $\beta_{n-1}$}}
 \put(114,-30){\makebox(0,0){\normalsize $\gamma$}}
 \put(86,-30){\makebox(0,0){\normalsize $\delta$}}
\end{picture}
\end{center}
and the ideal $I\lo n\po$ in the path algebra $K\Delta\lo n\po$
generated by the elements $\alpha_{n-1}\gamma$ and
$\alpha_{n-1}\delta$, $\psi_{n}$ is a canonical automorphism of
the repetitive algebra $\widehat{C}\lo n\po$ such that
$\psi_{n}^{2}=\nu_{\widehat{C}\lo n\po}$, and $\rho_{n}$ is a
rigid automorphism of $\widehat{C}\lo n\po$ with
$\rho_{n}^{2}=id_{\widehat{C}\lo n\po}$. Then $C\lo n\po$ is a
tubular extension of the path algebra of the unique convex
subquiver of $\Delta\lo n\po$ of Euclidean type
$\widetilde{\mathbbm{D}}_{n+4}$ of tubular type $\lo 2,2,2n\!+1
\!\po$, and so $C\lo n\po$ is a representation-infinite tilted
algebra of Euclidean type $\widetilde{\mathbbm{D}}_{2n+3}$ (see
\cite[(4.9)]{bib:Rin}). Hence the required shape of
$\Gamma_{\Gamma^{*}\lo n\po}^{s}$ follows from \cite[Section
4]{bib:AsNeSk} and \cite[Section 2]{bib:Sko}.
\end{proof}

Our next aim is to describe the shape of the stable
Auslander-Reiten quivers of the algebras $\Lambda\lo
p,q,k,s,\theta\po$ and the structure of modules lying on the mouth
of their large stable tubes.

Let $p$, $q$, $k$, $s$ be integers such that $p,q\geqslant 0$,
$k\geqslant 2$, $1\leqslant s\leqslant k-1$, $\gcd\lo s,k\po=1$,
and $\gcd\lo s+2,k\po=1$. Since $\gcd\lo s+2,k\po=1$, there is
exactly one integer $m_{k,s}$ such that $1\leqslant
m_{k,s}\leqslant k-1$, and
\begin{displaymath}
m_{k,s}\lo s+2\po+1\equiv 0 \lo \md k\po.
\end{displaymath}
By definition, $\Lambda\lo p,q,k,s,\theta\po$ is the bound quiver
algebra
\begin{displaymath}
\Omega^{\lo 1\po}\!\lo T_{p,q,k},\sigma_{s},\theta\po=KQ\lo
T_{p,q,k},\sigma_{s}\po/\overline{I}^{\lo 1\po}\!\lo
T_{p,q,k},\sigma_{s},\theta\po,
\end{displaymath}
where the vertices of the quiver $Q\lo T_{p,q,k},\sigma_{s}\po$
are the edges of the Brauer graph $T_{p,q,k}$, and there is an
arrow $i\longrightarrow j$ in $Q\lo T_{p,q,k},\sigma_{s}\po$ if
$j$ is the direct successor of $i$ in the partial order $p\lo
T_{p,q,k},\sigma_{s}\po$, as defined in Section \ref{morita}. We
note that $\Lambda\lo p,q,k,s,\theta\po$ is a special biserial
algebra \cite{bib:SkWa}. Moreover, since we compose the arrows in
$Q\lo T_{p,q,k},\sigma_{s}\po$ from left to right, the category
$\md \Lambda\lo p,q,k,s,\theta\po$ of finite dimensional (left)
$\Lambda\lo p,q,k,s,\theta\po$-modules is identified with the
category of finite dimensional representations (over $K$) of the
bound quiver opposite to the bound quiver $\lo Q\lo
T_{p,q,k},\sigma_{s}\po, \overline{I}^{\lo 1\po}\!\lo
T_{p,q,k},\sigma_{s},\theta\po\po$. For each edge $i$ of
$T_{p,q,k}$, we denote by $S\lo i\po$ the simple $\Lambda\lo
p,q,k,s,\theta\po$-module at $i$, and by $P\lo i\po$ the
projective cover of $S\lo i\po$ in $\md\Lambda\lo
p,q,k,s,\theta\po$. The edges of $T_{p,q,k}$ of the form $i^{r}$,
$1\leqslant i\leqslant k$, $1\leqslant r\leqslant p$, are said to
be outer edges of $T_{p,q,k}$, while the edges of the form
$i_{l}$, $1\leqslant i\leqslant k$, $1\leqslant l\leqslant q$, are
said to be the inner edges of $T_{p,q,k}$. We call the simple
modules $S\lo i^{r}\po$ the outer simple modules and the simple
modules $S\lo i_{l}\po$ the inner simple modules. We describe now
the structure of indecomposable projective $\Lambda\lo
p,q,k,s,\theta\po$-modules. Since $\Lambda\lo p,q,k,s,\theta\po$
is a special biserial algebra, every indecomposable projective
$\Lambda\lo p,q,k,s,\theta\po$-module $P$ is either serial (has a
unique composition series) or biserial ($\rad P/\soc P$ is a
direct sum of two serial modules). Observe also, that for each
vertex $i$ of $Q\lo T_{p,q,k},\sigma_{s}\po$, the vertex
$\sigma_{s}\lo i\po$ is the target of maximal paths in $Q\lo
T_{p,q,k},\sigma_{s}\po$ with source $i$ which are not in
$\overline{I}^{\lo 1\po}\!\lo T_{p,q,k},\sigma_{s},\theta\po$, and
hence $S\lo i\po$ is the socle of the projective module
$P\lo\sigma_{s}\lo i\po\po$. This also shows that the action of
$\sigma_{s}$ on the vertices of $Q\lo T_{p,q,k},\sigma_{s}\po$ is
exactly the Nakayama permutation of the algebra $\Lambda\lo
p,q,k,s,\theta\po$. For each outer edge $i^{r}$ of $T_{p,q,k}$,
the path
\begin{displaymath}
i^{r}\rightarrow \ldots \rightarrow i^{p}\rightarrow i\! +\!
s\rightarrow \lo i\! +\! s\po_{1}\rightarrow \ldots \rightarrow
\lo i\! +\! s\po_{q}\rightarrow i\! +\! s\! -\! 1\rightarrow \lo
i\! +\! s\po^{1}\rightarrow \ldots \rightarrow \lo i\!
+\!s\po^{r}=\sigma_{s}\lo i^{r}\po
\end{displaymath}
is the unique path in $Q\lo T_{p,q,k},\sigma_{s}\po$ from $i^{r}$
to $\sigma_{s}\lo i^{r}\po$ of length at least $2$ which is not in
$\overline{I}^{\lo 1\po}\!\lo T_{p,q,k},\sigma_{s},\theta\po$, and
consequently $P\lo\lo i\! +\! s\po^{r}\po$ is the serial module
with socle $S\lo i^{r}\po$, and the simple composition factors
given by the simple modules at the vertices of this path.
Similarly, for each inner edge $i_{l}$ of $T_{p,q,k}$, the path
\begin{displaymath}
i_{l}\rightarrow \ldots \rightarrow i_{q}\rightarrow i\! -\!
1\rightarrow i^{1}\rightarrow \ldots \rightarrow i^{p}\rightarrow
i\! +\! s\rightarrow \lo i\! +\! s\po_{1}\rightarrow \ldots
\rightarrow \lo i\! +\!s\po_{l}=\sigma_{s}\lo i_{l}\po
\end{displaymath}
is the unique path in $Q\lo T_{p,q,k},\sigma_{s}\po$ from $i_{l}$
to $\sigma_{s}\lo i_{l}\po$ of length at least $2$ which is not in
$\overline{I}^{\lo 1\po}\!\lo T_{p,q,k},\sigma_{s},\theta\po$, and
consequently $P\lo\lo i\! +\! s\po_{l}\po$ is the serial module
with socle $S\lo i_{l}\po$, and the simple composition factors
given by the simple modules at the vertices of this path. Finally,
for an edge $i$ of the cycle of $T_{p,q,k}$, we have two paths in
$Q\lo T_{p,q,k},\sigma_{s}\po$ from $i$ to $\sigma_{s}\lo i\po=i\!
+\! s$ which are not in $\overline{I}^{\lo 1\po}\!\lo
T_{p,q,k},\sigma_{s}\po$, and they are of the forms
\begin{displaymath}
i\rightarrow \lo i\! +\! 1\po^{1}\rightarrow \ldots \rightarrow
\lo i\! +\! 1\po^{p}\rightarrow i\! +\! 1\! +\! s\rightarrow \lo
i\! +\! 1\! +\! s\po_{1}\rightarrow \ldots \rightarrow \lo i\! +\!
1\! +\! s\po_{q} \rightarrow i\! +\! s,
\end{displaymath}
\begin{displaymath}
i\rightarrow i_{1}\rightarrow \ldots \rightarrow i_{q}\rightarrow
i\! -\! 1\rightarrow i^{1}\rightarrow \ldots \rightarrow i^{p}
\rightarrow i\! +\! s.
\end{displaymath}
Hence, $P\lo i\! +\! s\po$ is the biserial module with socle $S\lo
i\po$, and $\rad P\lo i\! +\! s\po/\soc P\lo i\! +\! s\po$ is the
direct sum of two serial modules: one with the top $S\lo\lo i\!
+\! 1\! +\! s\po_{q}\po$, the socle $S\lo\lo i\! +\! 1\po^{1}\po$,
and the simple composition factors formed by the simple modules at
the vertices of the first path, and the second with the top $S\lo
i^{p}\po$, the socle $S\lo i_{1}\po$, and the simple composition
factors formed by the simple modules at the vertices of the second
path. We also note that, for $p=0$, any path of the form
$i\rightarrow \lo i\! +\! 1\po^{1}\rightarrow \ldots \rightarrow
\lo i\! +\! 1\po^{p}\rightarrow i\! +\! 1\! +\! s$ is reduced to
the arrow $i\rightarrow i\! +\! 1\! +\! s$. Similarly, for $q=0$,
any path of the form $i\rightarrow i_{1}\rightarrow \ldots
\rightarrow i_{q}\rightarrow i\! -\! 1$ is reduced to the arrow
$i\rightarrow i\! -\! 1$. For $i\in \lk 1,2,\ldots,k\pk$, we
denote by $U\lo i\po$ the factor module of $P\lo i\po$ by the
unique serial submodule with the top $S\lo\lo i\! -\!
s\po^{p}\po$, if $p\geqslant 1$, and $S\lo i\! -\! s\! -\! 1\po$
if $p=0$. Similarly, for $i\in \lk 1,2,\ldots,k\pk$, we denote by
$V\lo i\po$ the factor module of $P\lo i\po$ by the unique serial
submodule with the top $S\lo\lo i\! +\! 1\po_{q}\po$, if
$q\geqslant 1$, and $S\lo i\! +\! 1\po$ if $q=0$.

The following proposition describes the stable Auslander-Reiten
quivers of the algebras $\Lambda\lo p,q,k,s,\theta\po$.
\begin{proposition} \label{structar}
$\Gamma_{\Lambda\lo p,q,k,s,\theta\po}^{s}$ consists of one
Euclidean component $\mathbbm{Z}\widetilde{\mathbbm{A}}_{2\lo
p+q+1\po k-1}$, a stable tube ${\cal T}_{p,q,k,s}^{out}$ of rank
$\lo 2p+1\po k$, a stable tube ${\cal T}_{p,q,k,s}^{inn}$ of rank
$\lo 2q+1\po k$, and a family of stable tubes of rank $1$ indexed
by $K\setminus\lk 0\pk$. Moreover,
\begin{arablist}
\item the mouth of the tube ${\cal T}_{p,q,k,s}^{out}$ consists of
the modules $U\lo i\po$, $S\lo i^{r}\po$, and $P\lo i^{r}\po/\soc
P\lo i^{r}\po$, $1\leqslant i\leqslant k$, $1\leqslant r\leqslant
p$;
\item the mouth of the tube ${\cal T}_{p,q,k,s}^{inn}$ consists of
the modules $V\lo i\po$, $S\lo i_{l}\po$, and $P\lo i_{l}\po/\soc
P\lo i_{l}\po$, $1\leqslant i\leqslant k$, $1\leqslant l\leqslant
q$.
\end{arablist}
\end{proposition}
\begin{proof}
Let $\Lambda=\Lambda\lo p,q,k,s,\theta\po$, $Q=Q\lo
T_{p,q,k},\sigma_{s}\po$ and $\overline{I}=\overline{I}^{\lo
1\po}\!\lo T_{p,q,k},\sigma_{s},\theta\po$. We know that $\Lambda$
is a special biserial one-parametric selfinjective algebra (see
Theorem \ref{thm1par}). Then it follows from \cite[Theorem
2.1]{bib:ErSk} that $\Gamma_{\Lambda}^{s}$ consists of one
Euclidean component $\mathbbm{Z}\widetilde{\mathbbm{A}}_{r}$ and a
family of stable tubes indexed by $\mathbbm{P}_{1}\lo
K\po=K\cup\lk\infty\pk$, containing a family of stable tubes of
rank $1$ indexed (say) by $K\setminus\lk 0\pk$. Moreover, if
$\Gamma_{\Lambda}^{s}$ admits two stable tubes of ranks $m$ and
$n$ greater than $1$ then $r=m\! +\! n\! -\! 1$. Therefore, since
$k\geqslant 2$, it remains to prove that $\Gamma_{\Lambda}^{s}$
admits the required stable tubes ${\cal T}_{p,q,k,s}^{out}$ of
rank $\lo 2p+1\po k$ and ${\cal T}_{p,q,k,s}^{inn}$ of rank $\lo
2q+1\po k$. We will apply the known formulae for the
Auslander-Reiten sequences of string modules over special biserial
algebras (see \cite[Chapter II]{bib:Erd}). We also note that the
ideal $\overline{I}$ is not an admissible ideal of the path
algebra $KQ$. Indeed, for any outer edge $i^{r}$ of $T_{p,q,k}$ we
have $C\lo i^{r},p\lo T_{p,q,k},\sigma_{s},v_{i}\po\po-C\lo
i^{r},p\lo T_{p,q,k},\sigma_{s},w\lo
i\po^{r}\po\po\in\overline{I}$, where $v\lo i\po$ is the vertex of
$i^{r}$ on the unique cycle of $T_{p,q,k}$ and $w\lo i\po^{r}$ the
second (outer) vertex of $i^{r}$, and $C\lo i^{r},p\lo
T_{p,q,k},\sigma_{s},w\lo i\po^{r}\po\po$ is the arrow
$i^{r}\rightarrow \sigma_{s}\lo i^{r}\po$. We call the arrow
$i^{r}\rightarrow \sigma_{s}\lo i^{r}\po$ a superfluous arrow of
the bound quiver $\lo Q,\overline{I}\po$. Similarly, for any inner
edge $i_{l}$ of $T_{p,q,k}$ we have $C\lo i_{l},p\lo
T_{p,q,k},\sigma_{s},v_{i}\po\po-C\lo i_{l},p\lo
T_{p,q,k},\sigma_{s},w\lo i\po_{l}\po\po\in\overline{I}$, where
$v\lo i\po$ is the vertex of $i_{l}$ on the unique cycle of
$T_{p,q,k}$ and $w\lo i\po_{l}$ the second (inner) vertex of
$i_{l}$, and $C\lo i_{l},p\lo T_{p,q,k},\sigma_{s},w\lo
i\po_{l}\po\po$ is the arrow $i_{l}\rightarrow \sigma_{s}\lo
i_{l}\po$. We call the arrow $i_{l}\rightarrow \sigma_{s}\lo
i_{l}\po$ a superfluous arrow of the bound quiver $\lo
Q,\overline{I}\po$.
\begin{arablist}
\item \textit{The tube ${\cal T}_{p,q,k,s}^{out}$}. Consider first
the case $p=0$. Observe that, for each $i\in\lk 1,2,\ldots,k\pk$,
the quiver $Q$ contains a string of the form
\begin{displaymath}
i\leftarrow \lo i\! +\! 1\po_{q}\leftarrow\ldots\leftarrow\lo i\!
+\! 1\po_{1}\leftarrow i\! +\! 1 \rightarrow i\! +\! s\! +\! 2
\leftarrow\lo i\! +\! s\! +\! 3\po_{q}\leftarrow
\ldots\leftarrow\lo i\! +\! s\! +\! 3\po_{1},
\end{displaymath}
and denote by $L\lo i\po$ the associated string $\Lambda$-module.
Then we have in $\md\Lambda$ an Auslander-Reiten sequence of the
form
\begin{displaymath}
0\rightarrow U\lo i\po\rightarrow L\lo i\po\rightarrow U\lo i\!
+\! s\! +\! 2 \po\rightarrow 0,
\end{displaymath}
and hence $U\lo i\po=\tau_{\Lambda}U\lo i\! +\! s\! +\! 2\po$.
Since $\gcd\lo s\! +\! 2,k\po=1$, we infer that the modules $U\lo
1\po, U\lo 2\po, \ldots, U\lo k\po$ form a periodic orbit of a
component ${\cal T}_{0,q,k,s}^{out}$ of $\Gamma_{\Lambda}^{s}$.
Clearly, ${\cal T}_{0,q,k,s}^{out}$ is a stable tube of rank $k$
and $U\lo 1\po, U\lo 2\po, \ldots, U\lo k\po$ lie on the mouth of
${\cal T}_{0,q,k,s}^{out}$, because the modules $L\lo 1\po, L\lo
2\po, \ldots, L\lo k\po$ are indecomposable.

\hspace{1cm}Assume now $p\geqslant 1$. Fix $i\in\lk
1,2,\ldots,k\pk$. Note that the quiver $Q$ contains a string of
the form
\begin{displaymath}
i^{p}\rightarrow i\! +\! s\leftarrow\lo i\! +\! s\! +\!
1\po_{q}\leftarrow\ldots\leftarrow\lo i\! +\! s\! +\!
1\po_{1}\leftarrow i\! +\! s\! +\! 1 \leftarrow \lo i\! +\!
1\po^{p}\leftarrow \ldots\leftarrow\lo i\! +\! 1\po^{1},
\end{displaymath}
and denote by $L\lo i\po$ the associated string $\Lambda$-module.
Then we have in $\md\Lambda$ an Auslander-Reiten sequence of the
form
\begin{displaymath}
0\rightarrow S\lo i^{p}\po\rightarrow L\lo i\po\rightarrow U\lo
i\! +\! s\po\rightarrow 0,
\end{displaymath}
and hence $S\lo i^{p}\po=\tau_{\Lambda}U\lo i\! +\! s\po$.
Moreover, we have $U\lo i\! +\! s\po=\rad P\lo\lo i\! +\! s\! +\!
1\po^{1}\po$. Further, $Q$ contains also a string of the form
\begin{displaymath}
i^{1}\leftarrow i\! -\! 1\rightarrow\lo i\! -\!
1\po_{1}\rightarrow\ldots\rightarrow\lo i\! -\!
1\po_{q}\rightarrow i\! -\! 2 \rightarrow \lo i\! -\!
1\po^{1}\rightarrow \ldots\rightarrow\lo i\! -\! 1\po^{p},
\end{displaymath}
and denote by $M\lo i\po$ the associated string $\Lambda$-module.
Then we have in $\md\Lambda$ an Auslander-Reiten sequence of the
form
\begin{displaymath}
0\rightarrow P\lo\lo i\! -\! 1\po^{p}\po/S\lo\lo i\! -\! 1\! -\! s
\po^{p}\po\rightarrow M\lo i\po\rightarrow S\lo
i^{1}\po\rightarrow 0,
\end{displaymath}
and hence $P\lo\lo i\! -\! 1\po^{p}\po/S\lo\lo i\! -\! 1\! -\! s
\po^{p}\po=\tau_{\Lambda}S\lo i^{1}\po$. Observe that $S\lo\lo i\!
-\! 1\! -\! s \po^{p}\po=\soc P\lo\lo i\! -\! 1\po^{p}\po$.
Finally, if $p\geqslant 2$, we have in $Q$ a path
\begin{displaymath}
i^{1}\rightarrow i^{2}\rightarrow \ldots\rightarrow i^{p},
\end{displaymath}
and, for $r\in\lk 1,\ldots,p\! -\! 1\pk$, $i^{r}\rightarrow
i^{r+1}$ is the unique nonsuperfluous arrow in $Q$ starting at
$i^{r}$ (respectively, ending in $i^{r+1}$). Then we have in
$\md\Lambda$ Auslander-Reiten sequences of the forms
\begin{displaymath}
0\rightarrow S\lo i^{r}\po\rightarrow L\lo i^{r}\po\rightarrow
S\lo i^{r+1}\po\rightarrow 0,
\end{displaymath}
$1\leqslant r\leqslant p\! -\! 1$, where $L\lo i^{r}\po$ is the
serial $\Lambda$-module given by the string $i^{r}\rightarrow
i^{r+1}$, and consequently $S\lo i^{r}\po=\tau_{\Lambda}S\lo
i^{r+1}\po$. Observe also that $\Omega_{\Lambda} S\lo
i^{r}\po=\rad P\lo i^{r}\po$. Therefore, $S\lo
i^{r}\po=\tau_{\Lambda}S\lo i^{r+1}\po$ implies $\rad P\lo
i^{r}\po=\tau_{\Lambda}\rad P\lo i^{r+1}\po$, for $r\in\lk
1,\ldots, p\! -\! 1\pk$, $p\geqslant 2$. Combining the
calculations above and invoking again our assumption $\gcd\lo s\!
+\! 2,k\po=1$, we conclude that the modules $U\lo i\po$, $S\lo
i^{r}\po$, $P\lo i^{r}\po/\soc P\lo i^{r}\po$, $1\leqslant
i\leqslant k$, $1\leqslant r\leqslant p$, form the mouth of a
stable tube ${\cal T}_{p,q,k,s}^{out}$ of $\Gamma_{\Lambda}^{s}$
of rank $\lo 2p\! +\! 1\po k=k\! +\! pk\! +\! pk$.
\item \textit{The tube ${\cal T}_{p,q,k,s}^{inn}$}. Consider first
the case $q=0$. For $i\in\lk 1,2,\ldots,k\pk$, the quiver $Q$
contains a string of the form
\begin{displaymath}
i\leftarrow \lo i\! -\! s\po^{p}\leftarrow\ldots\leftarrow\lo i\!
-\! s\po^{1}\leftarrow i\! -\! s\! -\! 1\rightarrow i\! -\! s\!
-\! 2 \leftarrow\lo i\! -\! 2s\! -\! 2\po^{p}\leftarrow
\ldots\leftarrow\lo i\! -\! 2s\! -\! 2\po^{1},
\end{displaymath}
and denote by $R\lo i\po$ the associated string $\Lambda$-module.
Then we have in $\md\Lambda$ an Auslander-Reiten sequence of the
form
\begin{displaymath}
0\rightarrow V\lo i\po\rightarrow R\lo i\po\rightarrow V\lo i\!
-\! s\! - \!2 \po\rightarrow 0,
\end{displaymath}
and hence $V\lo i\po=\tau_{\Lambda}V\lo i\! -\! s\! -\! 2\po$.
Since $\gcd\lo s\! +\! 2,k\po=1$, we infer as above that the
modules $V\lo 1\po, V\lo 2\po, \ldots, V\lo k\po$ form the mouth
of a stable tube ${\cal T}_{p,0,k,s}^{inn}$ of
$\Gamma_{\Lambda}^{s}$ of rank $k$.

\hspace{1cm}Assume now $q\geqslant 1$. Fix $i\in\lk
1,2,\ldots,k\pk$. Observe that $Q$ contains a string of the form
\begin{displaymath}
i_{p}\rightarrow i\! -\! 1\leftarrow\lo i\! -\! 1\! -\!
s\po^{p}\leftarrow\ldots\leftarrow\lo i\! -\! 1\! -\!
s\po^{1}\leftarrow i\! -\! 2\! -\! s \leftarrow \lo i\! -\! 1\!
-\! s\po_{q}\leftarrow \ldots\leftarrow\lo i\! -\! 1\! -\! s
\po_{1},
\end{displaymath}
and denote by $R\lo i\po$ the associated string $\Lambda$-module.
Then we have in $\md\Lambda$ an Auslander-Reiten sequence of the
form
\begin{displaymath}
0\rightarrow S\lo i_{q}\po\rightarrow R\lo i\po\rightarrow V\lo
i\! -\! 1\po\rightarrow 0,
\end{displaymath}
and hence $S\lo i_{q}\po=\tau_{\Lambda}V\lo i\! -\! 1\po$.
Moreover, we have $V\lo i\! -\! 1\po=\rad P\lo\lo i\! -\!
1\po_{1}\po$. Further, $Q$ contains also a string of the form
\begin{displaymath}
i_{1}\leftarrow i\rightarrow\lo i\! +\!
1\po^{1}\rightarrow\ldots\rightarrow\lo i\! +\!
1\po^{p}\rightarrow i\! +\! 1\! +\! s\rightarrow \lo i\! +\! 1\!
+\! s\po_{1}\rightarrow \ldots\rightarrow\lo i\! +\! 1\! +\!
s\po_{q},
\end{displaymath}
and denote by $N\lo i\po$ the associated string $\Lambda$-module.
Then we have in $\md\Lambda$ an Auslander-Reiten sequence of the
form
\begin{displaymath}
0\rightarrow P\lo\lo i\! +\! 1\! +\! s\po_{q}\po/S\lo\lo i\! +\!
1\po_{q}\po\rightarrow N\lo i\po\rightarrow S\lo
i_{1}\po\rightarrow 0,
\end{displaymath}
and hence $P\lo\lo i\! +\! 1\! +\! s\po_{q}\po/S\lo\lo i\! +\! 1
\po_q{p}\po=\tau_{\Lambda}S\lo i_{1}\po$. Observe also that
$S\lo\lo i\! +\! 1\po_{q}\po=\soc P\lo\lo i\! +\! 1\! +\!
s\po_{q}\po$. Finally, if $q\geqslant 2$, we have in $Q$ a path
\begin{displaymath}
i_{1}\rightarrow i_{2}\rightarrow \ldots\rightarrow i_{q},
\end{displaymath}
and, for $l\in\lk 1,\ldots,q\! -\! 1\pk$, $i_{l}\rightarrow
i_{l+1}$ is the unique nonsuperfluous arrow of $Q$ starting at
$i_{l}$ (respectively, ending in $i_{l+1}$). Then we have in
$\md\Lambda$ Auslander-Reiten sequences of the forms
\begin{displaymath}
0\rightarrow S\lo i_{l}\po\rightarrow W\lo i_{l}\po\rightarrow
S\lo i_{l+1}\po\rightarrow 0,
\end{displaymath}
$1\leqslant l\leqslant q\! -\! 1$, where $W\lo i_{l}\po$ is the
serial $\Lambda$-module given by the string $i_{l}\rightarrow
i_{l+1}$, and consequently $S\lo i_{l}\po=\tau_{\Lambda}S\lo
i_{l+1}\po$. Since $\Omega_{\Lambda} S\lo i_{l}\po=\rad P\lo
i_{l}\po$, $S\lo i_{l}\po=\tau_{\Lambda}S\lo i_{l+1}\po$ implies
$\rad P\lo i_{l}\po=\tau_{\Lambda}\rad P\lo i^{l+1}\po$. Moreover,
observe that $\tau_{\Lambda}^{2q+1}S\lo i_{1}\po=S\lo\lo i\! +\!
s\po_{1}\po$ for any $i\in\lk 1,2,\ldots,k\pk$. Invoking now our
assumption $\gcd\lo s,k\po=1$, we conclude that the modules $V\lo
i\po$, $S\lo i_{l}\po$, $P\lo i_{l}\po/\soc P\lo i_{l}\po$,
$1\leqslant i\leqslant k$, $1\leqslant l\leqslant q$, form the
mouth of a stable tube ${\cal T}_{p,q,k,s}^{inn}$ of
$\Gamma_{\Lambda}^{s}$ of rank $\lo 2q\! +\! 1\po k=k\! +\! qk\!
+\! qk$.
\end{arablist}
\end{proof}

We call ${\cal T}_{p,q,k,s}^{out}$ the \textit{outer stable tube}
of $\Gamma_{\Lambda\lo p,q,k,s,\theta\po}^{s}$ and ${\cal
T}_{p,q,k,s}^{inn}$ the \textit{inner stable tube} of
$\Gamma_{\Lambda\lo p,q,k,s,\theta\po}^{s}$.

The following proposition describes the action of the syzygy
operator on the outer and inner stable tubes of
$\Gamma_{\Lambda\lo p,q,k,s,\theta\po}^{s}$.

\begin{proposition}\label{action}
\begin{arablist}
\item $\Omega_{\Lambda\lo p,q,k,s,\theta\po}$ fixes the outer stable tube
${\cal T}_{p,q,k,s}^{out}$ of $\Gamma_{\Lambda\lo
p,q,k,s,\theta\po}^{s}$, and $\Omega_{\Lambda\lo
p,q,k,s,\theta\po}=\tau_{\Lambda\lo p,q,k,s,\theta\po}^{m_{k,s}\lo
2p+1\po+p+1}$ on ${\cal T}_{p,q,k,s}^{out}$.
\item $\Omega_{\Lambda\lo p,q,k,s,\theta\po}$ fixes the inner stable tube
${\cal T}_{p,q,k,s}^{inn}$ of $\Gamma_{\Lambda\lo
p,q,k,s,\theta\po}^{s}$, and $\Omega_{\Lambda\lo
p,q,k,s,\theta\po}=\tau_{\Lambda\lo p,q,k,s,\theta\po}^{\lo
k-m_{k,s}\po\lo 2q+1\po-q}$ on ${\cal T}_{p,q,k,s}^{inn}$.
\end{arablist}
\end{proposition}

\begin{proof}
Let $\Lambda=\Lambda\lo p,q,k,s,\theta\po$, ${\cal T}^{out}={\cal
T}_{p,q,k,s}^{out}$, ${\cal T}^{inn}={\cal T}_{p,q,k,s}^{inn}$,
and $m=m_{k,s}$. We will apply the formulae on the actions of
$\tau_{\Lambda}$ on ${\cal T}^{out}$ and ${\cal T}^{inn}$
established in the proof of Proposition \ref{structar}.
\begin{arablist}
\item It is enough to prove that
\begin{displaymath}
\Omega_{\Lambda}U\lo 1\po=\tau_{\Lambda}^{m\lo 2p+1\po+p+1}U\lo
1\po.
\end{displaymath}

Assume $p=0$. Observe that $\Omega_{\Lambda}^{-}U\lo 1\po=U\lo s\!
+\! 2\po$. Moreover, we have $\tau_{\Lambda}^{-m-1}U\lo 1\po=U\lo
1\! +\!\lo m\! +\! 1\po\lo s\! +\! 2\po\po=U\lo s\! +\! 3\! +\!
m\lo s\! +\! 2\po\po= U\lo s\! +\! 2\po$, since $m\lo s\! +\!
2\po\! +\! 1\equiv 0 \lo \md k\po$ implies $s\! +\! 3\! +\! m\lo
s\! +\! 2\po\equiv s\! +\! 2 \lo \md k\po$. Hence,
$\Omega_{\Lambda}^{-}U\lo 1\po=\tau_{\Lambda}^{-\lo m+1\po}U\lo
1\po$, or equivalently $\Omega_{\Lambda}U\lo
1\po=\tau_{\Lambda}^{m+1}U\lo 1\po$.

Assume $p\geqslant 1$. In this case, we have $U\lo 1\po=\rad P\lo
2^{1}\po$, and then $\Omega_{\Lambda}^{-}U\lo 1\po=S\lo 2^{1}\po$.
We have $\tau_{\Lambda}^{-p-1}U\lo 1\po=\tau_{\Lambda}^{-p-1}\rad
P\lo 2^{1}\po=\tau_{\Lambda}^{-}\tau_{\Lambda}^{-p}\rad P\lo
2^{1}\po=\tau_{\Lambda}^{-}P\lo 2^{p}\po/\soc P\lo 2^{p}\po=S\lo
3^{1}\po$ and $\tau_{\Lambda}^{-m\lo 2p+1\po}S\lo 3^{1}\po=S\lo\lo
3\! +\! m\lo s\! +\! 2\po\po^{1}\po=S\lo 2^{1}\po$, because $m\lo
s\! +\! 2\po\! +\! 1\equiv 0 \lo \md k\po$ implies $3\! +\! m\lo
s\! +\! 2\po\equiv 2 \lo \md k\po$. Therefore, we obtain
$\Omega_{\Lambda}^{-}U\lo 1\po=\tau_{\Lambda}^{-m\lo
2p+1\po-p-1}U\lo 1\po$, or equivalently $\Omega_{\Lambda}U\lo
1\po=\tau_{\Lambda}^{m\lo 2p+1\po+p+1}U\lo 1\po$.
\item It is enough to prove that
\begin{displaymath}
\Omega_{\Lambda}V\lo 1\po=\tau_{\Lambda}^{\lo k-m\po\lo
2q+1\po-q}V\lo 1\po.
\end{displaymath}

Assume $q=0$. Observe that $\Omega_{\Lambda}V\lo 1\po=V\lo 2\po$
and $\tau_{\Lambda}^{k-m}V\lo 1\po=V\lo 1\! +\!\lo k\! -\! m\po\lo
s\! +\! 2\po\po$. Then $\Omega_{\Lambda}V\lo
1\po=\tau_{\Lambda}^{k-m}V\lo 1\po$, because $m\lo s\! +\! 2\po\!
+\! 1\equiv 0 \lo \md k\po$ implies $1\! +\!\lo k\! -\! m\po\lo
s\! +\! 2\po\equiv 2 \lo \md k\po$.

Assume $q\geqslant 1$. In this case, we have $V\lo 1\po=\rad P\lo
1_{1}\po=\Omega_{\Lambda}S\lo 1_{1}\po$. Hence
$\Omega_{\Lambda}V\lo 1\po=\Omega_{\Lambda}^{2}S\lo 1_{1}\po$.
Thus $\Omega_{\Lambda}V\lo 1\po=\tau_{\Lambda}^{\lo k-m\po\lo
2q+1\po-q}V\lo 1\po$ is equivalent to $V\lo
1\po=\tau_{\Lambda}^{\lo k-m\po\lo 2q+1\po-q}S\lo 1_{1}\po$. We
have $\tau_{\Lambda}^{\lo k-m\po\lo 2q+1\po-q}S\lo 1_{1}\po=S\lo
\lo 1\! +\! \lo k\! -\! m\po\lo s\! +\! 2\po\po_{1}\po=S\lo
2_{1}\po$, because $m\lo s\! +\! 2\po\! +\! 1\equiv 0 \lo \md
k\po$ implies $1\! +\!\lo k\! -\! m\po\lo s\! +\! 2\po\equiv 2 \lo
\md k\po$. Therefore, we obtain $\tau_{\Lambda}^{\lo k-m\po\lo
2q+1\po-q}S\lo 1_{1}\po=\tau_{\Lambda}^{-q}S\lo
2_{1}\po=\tau_{\Lambda}^{-}S\lo 2_{q}\po=V\lo 1\po$, as required.
\end{arablist}
\end{proof}


\section{\hspace{-2.5ex}. Weakly symmetric algebras} \label{wsalg}


For the proof of Theorem \ref{mainthm} we need also a family of
one-parametric weakly symmetric algebras of Euclidean type
$\widetilde{\mathbbm{A}}_{m}$.

Let $T$ be a Brauer graph with exactly one cycle ${\cal R}_{k}$,
having an odd number $k$ of edges. Take the trivial rotation
$\sigma_{s}=\sigma_{0}$ (with $s=0$) of $T$. Observe that $\gcd\lo
s\! +\! 2,k\po=\gcd\lo 2,k\po=1$. Moreover, for each vertex $v$ of
$T$, the cyclic order $p\lo T,\sigma_{0},v\po$, as defined in
Section \ref{morita}, is the cyclic order of edges of the Brauer
graph $T$ adjacent to $v$. Then, for $\theta\in K\backslash\lk
0\pk$, consider the algebra $\Omega^{\lo 1\po}\!\lo
T,\theta\po=\Omega^{\lo 1\po}\!\lo T,\sigma_{0},\theta\po$. For
$\theta=1$, the algebra $\Omega^{\lo 1\po}\!\lo T,1\po$ is just
the one-parametric special biserial algebra $\Lambda^{\prime}\!\lo
T\po$ considered in \cite{bib:BoHoSk}, \cite{bib:BoSk1} and
\cite{bib:BoSk2}. In general, $\Omega^{\lo 1\po}\!\lo T,\theta\po$
is a one-parametric weakly symmetric special biserial algebra
which is socle equivalent to $\Omega^{\lo 1\po}\!\lo
T,1\po=\Lambda^{\prime}\!\lo T\po$.

For $p,q\geqslant 1$ and $k=1$, consider the Brauer graph
$T_{p,q}=T_{p,q,1}$ of the form
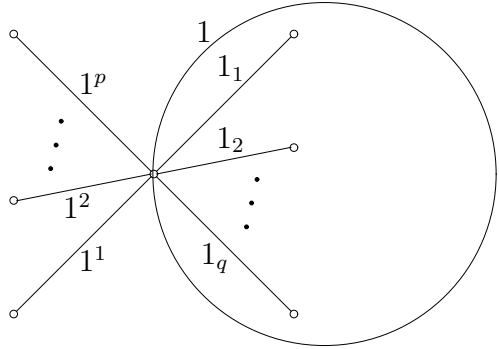
\begin{center}
\unitlength=1.0pt
\begin{picture}(191.00,131.00)(20.00,-6.00)
\put(81.00,59.00){\circle{3.00}} \put(134.00,69.00){\circle{3.00}}
\put(28.00,112.00){\circle{3.00}} \put(28.00,6.00){\circle{3.00}}
\put(134.00,112.00){\circle{3.00}}
\put(134.00,6.00){\circle{3.00}} \put(28.00,49.00){\circle{3.00}}
\put(82.00,59.00){\line(5,1){51.00}}
\put(80.00,60.00){\line(-1,1){51.00}}
\put(82.00,60.00){\line(1,1){51.00}}
\put(80.00,58.00){\line(-1,-1){51.00}}
\put(82.00,58.00){\line(1,-1){51.00}}
\put(80.00,59.00){\line(-5,-1){51.00}}
\put(100.00,113.00){\makebox(0,0){\normalsize $1$}}
\put(110.00,72.00){\makebox(0,0){\normalsize $1_{2}$}}
\put(52.00,47.00){\makebox(0,0){\normalsize $1^{2}$}}
\put(110.00,99.00){\makebox(0,0){\normalsize $1_{1}$}}
\put(58.00,93.00){\makebox(0,0){\normalsize $1^{p}$}}
\put(58.00,27.00){\makebox(0,0){\normalsize $1^{1}$}}
\put(104.00,27.00){\makebox(0,0){\normalsize $1_{q}$}}
\multiput(120.00,57.00)(-2,-9){3}{\circle*{1.5}}
\multiput(42.00,61.00)(2,9){3}{\circle*{1.5}}
\put(145.50,59.00){\circle{130}}
\end{picture}
\end{center}
Further, for $\theta\in K\backslash\lk 0\pk$, consider the
algebras $\Lambda\lo p,q,\theta\po=\Lambda\lo p,q,1,0,\theta\po$.
Then $\Lambda\lo p,q,\theta\po$ is a one-parametric weakly
symmetric special biserial algebra given by the quiver
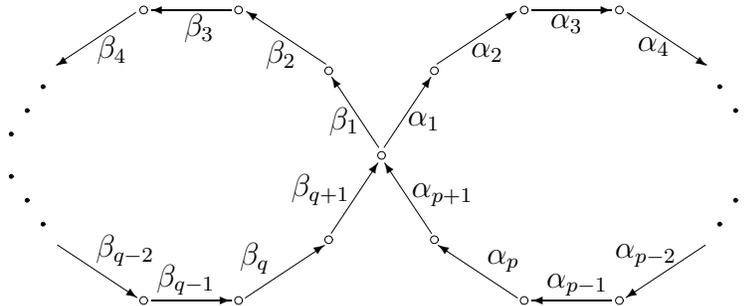
\begin{center}
\unitlength=1.0pt
\begin{picture}(301,126)(-48,-78)
 \put(99,-12){\vector(-2,3){17.8}}
 \put(101,-12){\vector(2,3){17.8}}
 \put(79,19.5){\vector(-3,2){30}}
 \put(121,19.5){\vector(3,2){30}}
 \put(43,40){\vector(-1,0){30}}
 \put(157,40){\vector(1,0){30}}
 \put(7,39){\vector(-3,-2){30}}
 \put(193,39){\vector(3,-2){30}}

 \put(81,-44.5){\vector(2,3){17.8}}
 \put(119,-44.5){\vector(-2,3){17.8}}
 \put(48.5,-69){\vector(3,2){30}}
 \put(151.5,-69){\vector(-3,2){30}}
 \put(13,-70){\vector(1,0){30}}
 \put(187,-70){\vector(-1,0){30}}
 \put(-22.5,-49){\vector(3,-2){30}}
 \put(222.5,-49){\vector(-3,-2){30}}

 \put(80,17){\circle{3}}
 \put(120,17){\circle{3}}
 \put(46,40){\circle{3}}
 \put(154,40){\circle{3}}
 \put(10,40){\circle{3}}
 \put(190,40){\circle{3}}

 \put(80,-47){\circle{3}}
 \put(120,-47){\circle{3}}
 \put(46,-70){\circle{3}}
 \put(154,-70){\circle{3}}
 \put(10,-70){\circle{3}}
 \put(190,-70){\circle{3}}

 \put(100,-15){\circle{3}}

 \multiput(-28,11)(-6,-9){3}{\circle*{1.5}}
 \multiput(-40,-23)(6,-9){3}{\circle*{1.5}}
 \multiput(228,11)(6,-9){3}{\circle*{1.5}}
 \multiput(240,-23)(-6,-9){3}{\circle*{1.5}}

 \put(86,-2){\makebox(0,0){\normalsize $\beta_{1}$}}
 \put(116,-2){\makebox(0,0){\normalsize $\alpha_{1}$}}
 \put(62,23){\makebox(0,0){\normalsize $\beta_{2}$}}
 \put(140,24){\makebox(0,0){\normalsize $\alpha_{2}$}}
 \put(31,33){\makebox(0,0){\normalsize $\beta_{3}$}}
 \put(170,35){\makebox(0,0){\normalsize $\alpha_{3}$}}
 \put(-2,25){\makebox(0,0){\normalsize $\beta_{4}$}}
 \put(203,26){\makebox(0,0){\normalsize $\alpha_{4}$}}

 \put(77,-28){\makebox(0,0){\normalsize $\beta_{q+1}$}}
 \put(123,-30){\makebox(0,0){\normalsize $\alpha_{p+1}$}}
 \put(52,-54){\makebox(0,0){\normalsize $\beta_{q}$}}
 \put(146,-55){\makebox(0,0){\normalsize $\alpha_{p}$}}
 \put(26,-62){\makebox(0,0){\normalsize $\beta_{q-1}$}}
 \put(174,-64){\makebox(0,0){\normalsize $\alpha_{p-1}$}}
 \put(3,-51){\makebox(0,0){\normalsize $\beta_{q-2}$}}
 \put(200,-53){\makebox(0,0){\normalsize $\alpha_{p-2}$}}
\end{picture}
\end{center}
and the relations:
\begin{displaymath}
\begin{array}{l}
\alpha_{1}\alpha_{2}\ldots\alpha_{p}\alpha_{p+1}\beta_{1}\beta_{2}\ldots\beta_{q}\beta_{q+1}=
 \theta\beta_{1}\beta_{2}\ldots\beta_{q}\beta_{q+1}\alpha_{1}\alpha_{2}\ldots\alpha_{p}\alpha_{p+1},\\
 \alpha_{p+1}\alpha_{1}=0,\; \beta_{q+1}\beta_{1}=0,\\
 \alpha_{i}\alpha_{i+1}\ldots\alpha_{p}\alpha_{p+1}\beta_{1}\ldots\beta_{q}\beta_{q+1}\alpha_{1}\ldots\alpha_{i-1}\alpha_{i}=0,\;
 2\leqslant i \leqslant p\!+\!1,\\
 \beta_{j}\beta_{j+1}\ldots\beta_{q}\beta_{q+1}\alpha_{1}\ldots\alpha_{p}\alpha_{p+1}\beta_{1}\ldots\beta_{j-1}\beta_{j}=0,\;
 2\leqslant j \leqslant q\!+\!1.
\end{array}
\end{displaymath}
We also note that for $p=q=0$, $\Lambda\lo 0,0,\theta\po$ is the
local $4$-dimensional algebra
$A\lo\theta\po=K\ld\alpha,\beta\pd/\lo\alpha^{2},\beta^{2},\alpha\beta-\theta\beta\alpha\po$,
and for $p\! +\! q\geqslant 1$ and $\theta=1$, $\Lambda\lo
p,q,1\po$ is the symmetric algebra $A\lo p\! +\! 1,q\! +\! 1\po$
considered in \cite{bib:BoHoSk}.

The following proposition is then a direct consequence of
arguments applied in the proofs of Proposition \ref{redstep1} and
Theorem \ref{normalform}.
\begin{tw}\label{prop81}
Let $T$ be a Brauer graph with exactly one cycle, having an odd
number $k$ of edges. Denote by $p$ the total number of outer edges
of $T$ and by $q$ the total number of inner edges of $T$. Then,
for any $\theta\in K\backslash\lk 0\pk$, the algebra $\Omega^{\lo
1\po}\! \lo T,\theta\po$ is derived equivalent to $\Lambda\lo
p^{\prime},q^{\prime},\theta\po$, where $p^{\prime}=p\! +\!
\frac{k-1}{2}$ and $q^{\prime}=q\! +\! \frac{k-1}{2}$.
\end{tw}

We also note that, in the symmetric case $\theta=1$, the above
proposition is related with \cite[Theorem 2.3]{bib:Kau}.

As we noticed above, for $p,q\geqslant 1$ and $\theta\in
K\backslash\lk 0\pk$, the algebras $\Lambda\lo p,q,\theta\po$ and
$\Lambda\lo p,q,1\po$ are socle equivalent (the factor algebras by
the socles are isomorphic). Then the following two propositions
follow directly from \cite[Propositions 5.3 and 5.4]{bib:BoHoSk}
and their proofs.
\begin{proposition}\label{prop82}
For $p,q\geqslant 0$ and $\theta\in K\backslash\lk 0\pk$, the
stable Auslander-Reiten quiver $\Gamma_{\Lambda\lo
p,q,\theta\po}^{s}$ consists of one Euclidean component
$\mathbbm{Z}\widetilde{\mathbbm{A}}_{2\lo p+q\po+1}$, an outer
stable tube ${\cal T}_{p,q,\theta}^{out}$ of rank $2p\! +\! 1$, an
inner stable tube ${\cal T}_{p,q,\theta}^{inn}$ of rank $2q\! +\!
1$, and a family of stable tubes of rank $1$.
\end{proposition}
\begin{proposition}\label{prop83}
For $p,q\geqslant 0$ and $\theta\in K\backslash\lk 0\pk$, we have
\begin{arablist}
\item $\Omega_{\Lambda\lo p,q,\theta\po}$ fixes the outer stable
tube ${\cal T}_{p,q,\theta}^{out}$ and $\Omega_{\Lambda\lo
p,q,\theta\po}=\tau_{\Lambda\lo p,q,\theta\po}^{p+1}$ on ${\cal
T}_{p,q,\theta}^{out}$.
\item $\Omega_{\Lambda\lo p,q,\theta\po}$ fixes the inner stable
tube ${\cal T}_{p,q,\theta}^{inn}$ and $\Omega_{\Lambda\lo
p,q,\theta\po}=\tau_{\Lambda\lo p,q,\theta\po}^{q+1}$ on ${\cal
T}_{p,q,\theta}^{inn}$.
\end{arablist}
\end{proposition}

For the proof of Theorem \ref{mainthm}, we need the following
result.
\begin{proposition}\label{prop84}
Let $p,p^{\prime},q,q^{\prime}\geqslant 0$, $k\geqslant 2$ and $s$
be integers with $1\leqslant s\leqslant k\! -\! 1$, $\gcd\lo s\!
+\! 2,k\po=1$, $\gcd\lo s,k\po=1$, and $\theta,\theta^{\prime}\in
K\backslash\lk 0\pk$. Then the algebras $\Lambda\lo
p,q,k,s,\theta\po$ and $\Lambda\lo
p^{\prime},q^{\prime},\theta\po$ are not stably equivalent.
\end{proposition}
\begin{proof}
Suppose that the algebras $\Lambda=\Lambda\lo p,q,k,s,\theta\po$
and $\Lambda^{\prime}=\Lambda\lo p^{\prime},q^{\prime},\theta\po$
are stably equivalent, and let
$F:\ulmod\Lambda\longrightarrow\ulmod\Lambda^{\prime}$ be a functor
inducing a stable
equivalence. Since $\Lambda$ and $\Lambda^{\prime}$ are
selfinjective algebras of Loewy length at least $3$, applying
\cite[(X.1.9), (X.1.12)]{bib:AuReSm}, we have
$F\tau_{\Lambda}=\tau_{\Lambda^{\prime}}F$ and
$F\Omega_{\Lambda}=\Omega_{\Lambda^{\prime}}F$. In particular, the
stable Auslander-Reiten quivers $\Gamma_{\Lambda}^{s}$ and
$\Gamma_{\Lambda^{\prime}}^{s}$ are isomorphic translation
quivers. Further, since $k\geqslant 2$, it follows from
Proposition \ref{structar} that the stable tubes ${\cal
T}_{\Lambda}^{out}$ and ${\cal T}_{\Lambda}^{inn}$ have ranks at
least $2$. Then, applying Proposition \ref{prop82}, we conclude
that $\Gamma_{\Lambda^{\prime}}^{s}$ has two stable tubes of ranks
at least $2$, and hence the stable tubes ${\cal
T}_{\Lambda^{\prime}}^{out}$ and ${\cal
T}_{\Lambda^{\prime}}^{inn}$ are unique stable tubes of ranks at
least $2$ in $\Gamma_{\Lambda^{\prime}}^{s}$. In particular, we
have $p^{\prime},q^{\prime}\geqslant 1$ and
\begin{displaymath}
\lk F\lo {\cal T}_{\Lambda}^{out}\po, F\lo {\cal
T}_{\Lambda}^{inn}\po\pk=\lk {\cal T}_{\Lambda^{\prime}}^{out},
{\cal T}_{\Lambda^{\prime}}^{inn}\pk.
\end{displaymath}
If $F\lo {\cal T}_{\Lambda}^{out}\po={\cal
T}_{\Lambda^{\prime}}^{out}$ then, applying Propositions
\ref{structar}, \ref{action}, \ref{prop82} and \ref{prop83}, we
get the equalities
\begin{displaymath}
\lo 2p \!+\! 1\po k=2p^{\prime} \!+\! 1 \mbox{ and } m\lo 2p \!+\!
1\po\!+\! p\! +\! 1=p^{\prime}\! +\! 1,
\end{displaymath}
where $m=m_{k,s}$ is such that $m\lo s\! +\! 2\po\! +\! 1\equiv
0\lo\md k\po$. Similarly, if $F\lo {\cal
T}_{\Lambda}^{out}\po={\cal T}_{\Lambda^{\prime}}^{inn}$, we get
the equalities
\begin{displaymath}
\lo 2p \!+\! 1\po k=2q^{\prime} \!+\! 1 \mbox{ and } m\lo 2p \!+\!
1\po\!+\! p\! +\! 1=q^{\prime}\! +\! 1.
\end{displaymath}
The both cases force $k=2m\! +\! 1$. On the other hand, $m\lo s\!
+\! 2\po\! +\! 1\equiv 0\lo\md k\po$ is equivalent to $\lo 2m\!
+\! 1\po\! +\! sm\equiv 0\lo \md k\po$. Hence $k=2m\! +\! 1$
implies $2m\! +\! 1 | sm$, a contradiction, because $1\leqslant
s\leqslant k\! -\! 1=2m$ and $\gcd\lo 2m\! +\! 1,s\po=\gcd\lo
k,s\po=1$.
\end{proof}


\section{\hspace{-2.5ex}. Proof of Theorem \ref{mainthm}} \label{proofthm1}


The aim of this section is to complete the proof of Theorem
\ref{mainthm}. The implication $\lo 1\po\Longrightarrow\lo 2\po$
follows from Theorem \ref{normalform} and Proposition
\ref{normaldn}. The implication $\lo 2\po\Longrightarrow\lo 3\po$
is the direct consequence of the following general result proved
by J. Rickard.

\begin{proposition}\label{prop71}
Let $A$ and $\Lambda$ be derived equivalent selfinjective
algebras. Then $A$ and $\Lambda$ are stably equivalent.
\end{proposition}
\begin{proof}
See \cite[Corollary 2.2]{bib:Ric2}.
\end{proof}

We will need also the following general result.

\begin{proposition} \label{prop72}
Let $A$ and $\Lambda$ be stably equivalent selfinjective algebras,
and assume that $\Lambda$ is symmetric of Loewy length at least
$3$. Then $A$ is weakly symmetric.
\end{proposition}
\begin{proof}
See \cite[Proposition 5.2]{bib:BoHoSk}.
\end{proof}

We are now in position to prove the implication $\lo
3\po\Longrightarrow\lo 1\po$ of Theorem \ref{mainthm}. Let $A$ be
a standard selfinjective algebra which is stably equivalent to an
algebra of the form $\Lambda\lo p,q,k,s,\theta\po$ or
$\Gamma^{*}\lo n\po$. Then the stable Auslander-Reiten quiver
$\Gamma_{A}^{s}$ of $A$ is isomorphic either to
$\Gamma_{\Lambda\lo p,q,k,s,\theta\po}^{s}$ or to
$\Gamma_{\Gamma^{*}\lo n\po}$. In particular, $A$ is a
one-parametric selfinjective algebra of an Euclidean type
$\widetilde{\mathbbm{A}}_{2\lo p+q+1\po k-1}$ or
$\widetilde{\mathbbm{D}}_{2n+3}$ (see \cite{bib:Sko} and
Propositions \ref{structard} and \ref{structar}). Moreover, since
the algebras $\Lambda\lo p,q,k,s,\theta\po$ and $\Gamma^{*}\!\lo
n\po$ are not weakly symmetric, applying Proposition \ref{prop72},
we conclude that $A$ is not symmetric. Then applying the main
results of \cite{bib:BoSk1}, \cite{bib:BoSk2} and \cite{bib:BoSk3}
we obtain that $A$ is isomorphic to an algebra of one of the forms
$\Omega^{\lo 1\po}\!\lo T,\sigma_{s},\theta\po$, $\Omega^{\lo
1\po}\!\lo T,\theta\po$, or $\Omega^{\lo 2\po}\!\lo
T,v_{1},v_{2}\po$. Suppose $A$ is isomorphic to an algebra
$\Omega^{\lo 1\po}\!\lo T,\theta\po$. Then it follows from Theorem
\ref{prop81} and Proposition \ref{prop71} that $A$ is stably
equivalent to an algebra of the form $\Lambda\lo
p^{\prime},q^{\prime},\theta^{\prime}\po$, for same
$p^{\prime},q^{\prime}\geqslant 0$ and $\theta^{\prime}\in
K\backslash\lk 0\pk$. Moreover, by our assumption $\lo 3\po$, $A$
is stably equivalent to an algebra $\Lambda\lo p,q,k,s,\theta\po$
or $\Gamma^{*}\!\lo n\po$. Since the stable Auslander-Reiten
quivers $\Gamma_{\Lambda\lo
p^{\prime},q^{\prime},\theta^{\prime}\po}^{s}$ and
$\Gamma_{\Gamma^{*}\!\lo n\po}$ have, by Propositions
\ref{structard} and \ref{prop82}, the unique Euclidean components
of different types, we conclude that $A$ is stably equivalent
simultaneously to an algebra $\Lambda\lo p,q,k,s,\theta\po$ and to
an algebra $\Lambda\lo p^{\prime},q^{\prime},\theta^{\prime}\po$.
But this contradicts Proposition \ref{prop84}. Therefore, the
implication $\lo 3\po\Longrightarrow\lo 1\po$ holds. This finishes
the proof of Theorem \ref{mainthm}.


\section{\hspace{-2.5ex}. Stable equivalences} \label{stableequiv}


It follows from Proposition \ref{structard} that two algebras
$\Gamma^{*}\lo m\po$ and $\Gamma^{*}\lo n\po$ are stably
equivalent if and only if $m=n$, or equivalently $\Gamma^{*}\lo
m\po$ and $\Gamma^{*}\lo n\po$ are isomorphic. Clearly, two
algebras $\Lambda\lo p,q,k,s,\theta\po$ and $\Gamma^{*}\lo n\po$
are never stably equivalent, because their stable Auslander-Reiten
quivers are not isomorphic. The following proposition gives
necessary conditions for two algebras $\Lambda\lo
p,q,k,s,\theta\po$ and $\Lambda\lo
p^{\prime},q^{\prime},k^{\prime},s^{\prime},\theta^{\prime}\po$ to
be stably equivalent.

\begin{proposition} \label{propiso}
Let $\Lambda\lo p,q,k,s,\theta\po$ and $\Lambda\lo
p^{\prime},q^{\prime},k^{\prime},s^{\prime},\theta^{\prime}\po$ be
stably equivalent algebras. Then one of two cases holds:
\begin{arablist}
\item $p=p^{\prime}$, $q=q^{\prime}$, $k=k^{\prime}$,
$s=s^{\prime}$;
\item $p=q^{\prime}$, $q=p^{\prime}$, $k=k^{\prime}$, $s\neq
s^{\prime}$, $m_{k,s}+m_{k^{\prime},s^{\prime}}=k-1$.
\end{arablist}
\end{proposition}

\begin{proof}
We abbreviate $\Lambda=\Lambda\lo p,q,k,s,\theta\po$ and
$\Lambda^{\prime}=\Lambda\lo
p^{\prime},q^{\prime},k^{\prime},s^{\prime},\theta^{\prime}\po$.
We have
\begin{displaymath}
p,q\geqslant 0,\mbox{ } k\geqslant 2,\mbox{ } 1\leqslant
s\leqslant k-1,\mbox{ } \gcd\lo s,k\po=1,\mbox{ } \gcd\lo
s+2,k\po=1,
\end{displaymath}
\begin{displaymath}
p^{\prime},q^{\prime}\geqslant 0,\mbox{ } k^{\prime}\geqslant
2,\mbox{ } 1\leqslant s^{\prime}\leqslant k^{\prime}-1,\mbox{ }
\gcd\lo s^{\prime},k^{\prime}\po=1,\mbox{ } \gcd\lo
s^{\prime}+2,k^{\prime}\po=1.
\end{displaymath}
We also note that since $k\geqslant 2$ and $k^{\prime}\geqslant
2$, $\Lambda$ and $\Lambda^{\prime}$ are of Loewy length at least
$3$, and then for any stable equivalence
$F:\ulmod\Lambda\longrightarrow \ulmod\Lambda^{\prime}$ we have
$F\tau_{\Lambda}=\tau_{\Lambda^{\prime}}F$ and
$F\Omega_{\Lambda}=\Omega_{\Lambda^{\prime}}F$ (see \cite[(X.1.9),
(X.1.12)]{bib:AuReSm}). Recall also that $\Gamma_{\Lambda}^{s}$
has exactly two large stable tubes: the outer tube ${\cal
T}_{p,q,k,s}^{out}$ of rank $\lo 2p\! +\! 1\po k$ and the inner
tube ${\cal T}_{p,q,k,s}^{inn}$ of rank $\lo 2q\! +\! 1\po k$.
Similary, $\Gamma_{\Lambda^{\prime}}^{s}$ has exactly two large
stable tubes: the outer tube ${\cal
T}_{p^{\prime},q^{\prime},k^{\prime},s^{\prime}}^{out}$ of rank
$\lo 2p^{\prime}\! +\!1\po k^{\prime}$ and the inner tube ${\cal
T}_{p^{\prime},q^{\prime},k^{\prime},s^{\prime}}^{inn}$ of rank
$\lo 2q^{\prime}\! +\!1\po k^{\prime}$. Since the stable
Auslander-Reiten quivers of $\Lambda$ and $\Lambda^{\prime}$ are
isomorphic, we obtain
\begin{displaymath}
\lk\lo 2p\! +\! 1\po k,\lo 2q\! +\! 1\po k\pk=\lk \lo
2p^{\prime}\! +\! 1\po k^{\prime},\lo 2q^{\prime}\! +\! 1\po
k^{\prime}\pk.
\end{displaymath}
Hence, we get
\begin{displaymath}
\lo 2p\! +\! 1\po k\! +\!\lo 2q\! +\! 1\po k=\lo 2p^{\prime}\! +\!
1\po k^{\prime}\! +\!\lo 2q^{\prime}\! +\! 1\po k^{\prime},
\end{displaymath}
and consequently the equality
\begin{displaymath}
k\lo p\! +\! q\! +\! 1\po= k^{\prime}\lo p^{\prime}\! +\!
q^{\prime}\! +\! 1\po.
\end{displaymath}

Let $m=m_{k,s}$ and $m^{\prime}=m_{k^{\prime},s^{\prime}}$. We have
the congruences
\begin{displaymath}
m\lo s\! +\! 2\po\! +\! 1\equiv 0 \lo \md k\po \mbox{ and }
m^{\prime}\lo s^{\prime}\! +\! 2\po\! +\! 1\equiv 0 \lo \md
k^{\prime}\po,
\end{displaymath}
or equivalently
\begin{displaymath}
\lo 2m\! +\! 1\po\! +\! sm\equiv 0 \lo \md k\po \mbox{ and } \lo
2m^{\prime}\! +\! 1\po\! +\! s^{\prime}m^{\prime}\equiv 0 \lo \md
k^{\prime}\po.
\end{displaymath}
Hence $\gcd\lo m,k\po=1$ and $\gcd\lo m^{\prime},k^{\prime}\po=1$.
Then $\gcd\lo s,k\po=1$ and $\gcd\lo s^{\prime},k^{\prime}\po=1$
imply $\gcd\lo 2m\! +\! 1,k\po=1$ and $\gcd\lo 2m^{\prime}\! +\!
1,k^{\prime}\po=1$.

We have two cases to consider. Assume first that ${\cal
T}_{p,q,k,s}^{out}={\cal
T}_{p^{\prime},q^{\prime},k^{\prime},s^{\prime}}^{out}$ and ${\cal
T}_{p,q,k,s}^{inn}={\cal
T}_{p^{\prime},q^{\prime},k^{\prime},s^{\prime}}^{inn}$. Comparing
the ranks of large stable tubes of $\Gamma_{\Lambda}^{s}$ and
$\Gamma_{\Lambda^{\prime}}^{s}$, we then have
\begin{displaymath}
\lo 2p\! +\! 1\po k=\lo 2p^{\prime}\! +\! 1\po k^{\prime} \mbox{
and } \lo 2q\! +\! 1\po k=\lo 2q^{\prime}\! +\! 1\po k^{\prime}.
\end{displaymath}
Further, comparing the actions of $\Omega_{\Lambda}$ and
$\Omega_{\Lambda^{\prime}}$ on the large stable tubes, we obtain
from Proposition \ref{action} the equalities
\begin{displaymath}
m\lo 2p\! +\! 1\po\! +\! p\! +\! 1=m^{\prime}\lo 2p^{\prime}\! +\!
1\po\! +\! p^{\prime}\! +\! 1 \mbox{ and } \lo k\! -\! m\po\lo
2q\! +\! 1\po\! -\! q=\lo k^{\prime}\! -\! m^{\prime}\po\lo
2q^{\prime}\! +\! 1\po\! -\! q^{\prime}.
\end{displaymath}
Moreover, the second equality and $\lo 2q\! +\! 1\po k=\lo
2q^{\prime}\! +\! 1\po k^{\prime}$ imply $m\lo 2q\! +\! 1\po\! +\!
q=m^{\prime}\lo 2q^{\prime}\! +\! 1\po\! +\! q^{\prime}$. Summing
up both equalities we obtain $\lo 2m\! +\! 1\po\lo p\! +\! q\! +\!
1\po=\lo 2m^{\prime}\! +\! 1\po\lo p^{\prime}\! +\! q^{\prime}\!
+\! 1\po$. Finally, invoking the equality $k\lo p\! +\! q\! +\!
1\po=k^{\prime}\lo p^{\prime}\! +\! q^{\prime}\! +\! 1\po$, we
obtain
\begin{displaymath}
\lo 2m\! +\! 1\po k^{\prime}=\lo 2m^{\prime}\! +\! 1\po k.
\end{displaymath}
But then $\gcd\lo 2m\! +\! 1,k\po=1$ and $\gcd\lo 2m^{\prime}\!
+\! 1,k^{\prime}\po=1$ imply $k | k^{\prime}$ and $k^{\prime} |
k$. Hence $k=k^{\prime}$ and consequently $p=p^{\prime}$,
$q=q^{\prime}$, $m=m^{\prime}$.

Assume that ${\cal T}_{p,q,k,s}^{out}={\cal
T}_{p^{\prime},q^{\prime},k^{\prime},s^{\prime}}^{inn}$ and ${\cal
T}_{p,q,k,s}^{inn}={\cal
T}_{p^{\prime},q^{\prime},k^{\prime},s^{\prime}}^{out}$. Then
comparing the ranks of these tubes, we obtain
\begin{displaymath}
\lo 2p\! +\! 1\po k=\lo 2q^{\prime}\! +\! 1\po k^{\prime} \mbox{
and } \lo 2q\! +\! 1\po k=\lo 2p^{\prime}\! +\! 1\po k^{\prime}.
\end{displaymath}
Further, comparing the actions of $\Omega_{\Lambda}$ and
$\Omega_{\Lambda^{\prime}}$ on these tubes, we obtain from
Proposition \ref{action} the equalities
\begin{displaymath}
m\lo 2p\! +\! 1\po\! +\! p\! +\! 1=\lo k^{\prime}\! -\!
m^{\prime}\po\lo 2q^{\prime}\! +\! 1\po\! -\! q^{\prime} \mbox{
and } \lo k\! -\! m\po\lo 2q\! +\! 1\po\! -\! q=m^{\prime}\lo
2p^{\prime}\! +\! 1\po\! +\! p^{\prime}\! +\! 1.
\end{displaymath}
Hence, we obtain
\begin{eqnarray*}
m^{\prime}\lo 2q^{\prime}\! +\! 1\po\! +\! q^{\prime} & = &
k^{\prime}\lo 2q^{\prime}\! +\! 1\po\! -\! m\lo 2p\! +\! 1\po\!
-\! p\! -\! 1 \\
 & = & k\lo 2p\! +\! 1\po\! -\! m\lo 2p\! +\! 1\po\!
-\! p\! -\! 1 \\
& = & \lo k\! -\! m\po\lo 2p\! +\! 1\po\! -\! p\! -\! 1.
\end{eqnarray*}
Summing up, we get the equality
\begin{displaymath}
\lo 2m^{\prime}\! +\! 1\po\lo p^{\prime}\! +\! q^{\prime}\! +\!
1\po=\lo 2\lo k\! -\! m\po\! -\! 1\po\lo p\! +\! q\! +\! 1\po.
\end{displaymath}
Invoking again the equality $k\lo p\! +\! q\! +\!
1\po=k^{\prime}\lo p^{\prime}\! +\! q^{\prime}\! +\! 1\po$, we
obtain
\begin{displaymath}
\lo 2m^{\prime}\! +\! 1\po k=\lo 2\lo k\! -\! m\po\! -\! 1\po
k^{\prime}.
\end{displaymath}
We have $\gcd\lo 2m^{\prime}\! +\! 1,k^{\prime}\po=1$ and $\gcd\lo
2m\! +\! 1,k\po=1$. Moreover, $\gcd\lo 2m\! +\! 1,k\po=1$ implies
$\gcd\lo 2\lo k\! -\! m\po\! -\! 1,k\po=1$. Therefore, $k |
k^{\prime}$ and $k^{\prime} | k$, and so $k=k^{\prime}$. Then
$p=q^{\prime}$, $q=p^{\prime}$. Finally, $ 2m^{\prime}\! +\! 1=
2\lo k\! -\! m\po\! -\! 1$, and consequently $m\! +\!
m^{\prime}=k\! -\! 1$.
\end{proof}

We note that the condition (2) in the above proposition forces $k$
to be odd. Indeed, if $k$ is even then $m_{k,s}$ and
$m_{k^{\prime},s^{\prime}}$ have different parity, and this
contradicts the congruence $m_{k,s}\lo s\! +\! 2\po\! +\!
m_{k^{\prime},s^{\prime}}\lo s^{\prime}\! +\! 2\po\! +\! 2\equiv 0
\lo\md k\po$, because $s\! +\! 2$ and $s^{\prime}\! +\! 2$ are odd
by $\gcd\lo s,k\po=1$ and $\gcd\lo
s^{\prime},k^{\prime}\po=\gcd\lo s^{\prime},k\po=1$.

We also have the following fact.
\begin{proposition}\label{prop102}
Two algebras $\Lambda\lo p,q,k,s,\theta\po$ and $\Lambda\lo
p^{\prime},q^{\prime},k^{\prime},s^{\prime},\theta^{\prime}\po$
are isomorphic if and only if one of the following cases holds:
\begin{arablist}
\item $p=p^{\prime}$, $q=q^{\prime}$, $k=k^{\prime}$,
$s=s^{\prime}$, $\theta^{\prime}=\theta$;
\item $p=q^{\prime}$, $q=p^{\prime}$, $k=k^{\prime}$, $s\neq
s^{\prime}$, $m_{k,s}+m_{k^{\prime},s^{\prime}}=k-1$,
$\theta^{\prime}=\theta^{-1}$.
\end{arablist}
\end{proposition}
\begin{proof}
Let
\begin{displaymath}
\Lambda\lo p,q,k,s,\theta\po=KQ\lo
T_{p,q,k},\sigma_{s}\po/\overline{I}^{\lo 1\po}\!\lo
T_{p,q,k},\sigma_{s},\theta\po
\end{displaymath}
and
\begin{displaymath}
\Lambda\lo
p^{\prime},q^{\prime},k^{\prime},s^{\prime},\theta^{\prime}\po=KQ\lo
T_{p^{\prime},q^{\prime},k^{\prime}},\sigma_{s^{\prime}}\po/\overline{I}^{\lo
1\po}\!\lo
T_{p^{\prime},q^{\prime},k^{\prime}},\sigma_{s^{\prime}},\theta^{\prime}\po.
\end{displaymath}
We denote by $i_{j}$, for $j=1,2,\ldots q$, $i=1,2,\ldots,k$, the
vertices of the quiver $Q\lo T_{p,q,k},\sigma_{s}\po$
(respectively, $Q\lo
T_{p^{\prime},q^{\prime},k^{\prime}},\sigma_{s^{\prime}}\po$)
corresponding to the inner edges of the graph $T_{p,q,k}$
(respectively, $T_{p^{\prime},q^{\prime},k^{\prime}}$) and by
$i^{j}$ , for $j=1,2,\ldots p$, $i=1,2,\ldots,k$, the vertices of
the quiver $Q\lo T_{p,q,k},\sigma_{s}\po$ (respectively, $Q\lo
T_{p^{\prime},q^{\prime},k^{\prime}},\sigma_{s^{\prime}}\po$)
corresponding to the outer edges of the graph $T_{p,q,k}$
(respectively, $T_{p^{\prime},q^{\prime},k^{\prime}}$). Moreover,
we denote by $\alpha_{i_{j}}:i_{j}\longrightarrow i_{j+1}$, for
$j=1,2,\ldots q\! -\! 1$, $\alpha_{i_{0}}:i\longrightarrow i_{1}$,
$\alpha_{i_{q}}:i_{q}\longrightarrow i\! -\! 1$,
$\beta_{i^{j}}:i^{j}\longrightarrow i^{j+1}$, for $j=1,2,\ldots
p\! -\! 1$, $\beta_{i^{0}}:i\! -\! 1\longrightarrow i^{1}$,
$\beta_{i^{p}}:i^{p}\longrightarrow i\! +\! s$, for
$i=1,2,\ldots,k$, the arrows of the quiver $Q\lo
T_{p,q,k},\sigma_{s}\po$ (respectively, $Q\lo
T_{p^{\prime},q^{\prime},k^{\prime}},\sigma_{s^{\prime}}\po$). The
ideal $\overline{I}^{\lo 1\po}\!\lo
T_{p,q,k},\sigma_{s},\theta\po$ contains the following essential
commutativity generators:
\begin{displaymath}
\alpha_{i_{0}}\alpha_{i_{1}}\ldots\alpha_{i_{q}}\beta_{i^{0}}\beta_{i^{1}}\ldots\beta_{i^{p}}-
\beta_{\lo i\! +\! 1\po ^{0}}\beta_{\lo i\! +\!
1\po^{1}}\ldots\beta_{\lo i\! +\! 1\po^{p}}\alpha_{\lo i\! +\! 1\!
+\! s\po_{0}}\alpha_{\lo i\! +\! 1\! +\!
s\po_{1}}\ldots\alpha_{\lo i\! +\! 1\! +\! s\po_{q}},
\end{displaymath}
for $i=2,3,\ldots k$, and
\begin{displaymath}
\alpha_{1_{0}}\alpha_{1_{1}}\ldots\alpha_{1_{q}}\beta_{1^{0}}\beta_{1^{1}}\ldots\beta_{1^{p}}-
\theta\beta_{2^{0}}\beta_{2^{1}}\ldots\beta_{2^{p}}\alpha_{\lo s\!
+\! 2\po_{0}}\alpha_{\lo s\! +\! 2\po_{1}}\ldots\alpha_{\lo s\!
+\! 2\po_{q}}.
\end{displaymath}
Further, the ideal
$\overline{I}^{\lo 1\po}\!\lo
T_{p^{\prime},q^{\prime},k^{\prime}},\sigma_{s^{\prime}},\theta^{\prime}\po$
contains the following essential commutativity generators:
\begin{displaymath}
\alpha_{i_{0}}\alpha_{i_{1}}\ldots\alpha_{i_{q}}\beta_{i^{0}}\beta_{i^{1}}\ldots\beta_{i^{p}}-
\beta_{\lo i\! +\! 1\po ^{0}}\beta_{\lo i\! +\!
1\po^{1}}\ldots\beta_{\lo i\! +\! 1\po^{p}}\alpha_{\lo i\! +\! 1\!
+\! s^{\prime}\po_{0}}\alpha_{\lo i\! +\! 1\! +\!
s^{\prime}\po_{1}}\ldots\alpha_{\lo i\! +\! 1\! +\!
s^{\prime}\po_{q}},
\end{displaymath}
for $i=2,3,\ldots k$, and
\begin{displaymath}
\alpha_{1_{0}}\alpha_{1_{1}}\ldots\alpha_{1_{q}}\beta_{1^{0}}\beta_{1^{1}}\ldots\beta_{1^{p}}-
\theta^{\prime}\beta_{2^{0}}\beta_{2^{1}}\ldots\beta_{2^{p}}\alpha_{\lo
s^{\prime}\! +\! 2\po_{0}}\alpha_{\lo s^{\prime}\! +\!
2\po_{1}}\ldots\alpha_{\lo s^{\prime}\! +\! 2\po_{q}}.
\end{displaymath}
Assume that the algebras $\Lambda\lo p,q,k,s,\theta\po$ and
$\Lambda\lo
p^{\prime},q^{\prime},k^{\prime},s^{\prime},\theta^{\prime}\po$
are isomorphic and let \linebreak $f:\Lambda\lo
p,q,k,s,\theta\po\longrightarrow\Lambda\lo
p^{\prime},q^{\prime},k^{\prime},s^{\prime},\theta^{\prime}\po$ be
a $K$-algebra isomorphism. In particular, the algebras $\Lambda\lo
p,q,k,s,\theta\po$ and $\Lambda\lo
p^{\prime},q^{\prime},k^{\prime},s^{\prime},\theta^{\prime}\po$
are stably equivalent. Then it follows from Proposition
\ref{propiso} that one of the following cases holds:
\begin{alphlist}
\item $p=p^{\prime}$, $q=q^{\prime}$, $k=k^{\prime}$,
$s=s^{\prime}$;
\item $p=q^{\prime}$, $q=p^{\prime}$, $k=k^{\prime}$, $s\neq
s^{\prime}$, $m_{k,s}+m_{k^{\prime},s^{\prime}}=k-1$.
\end{alphlist}

Assume (a) holds. Then we may assume (without loss of generality)
that $f\lo\alpha_{i_{j}}\po=\alpha_{i_{j}}$, for $j=0,1,\ldots,q$,
$i=1,s\! +\! 2$, and $f\lo\beta_{i^{j}}\po=\beta_{i^{j}}$, for
$j=0,1,\ldots,p$, $i=1,2$. Moreover, we have the equalities:
\begin{equation*}
\begin{split}
&
\theta^{\prime}\beta_{2^{0}}\beta_{2^{1}}\ldots\beta_{2^{p}}\alpha_{\lo
s+2\po_{0}}\alpha_{\lo s+2\po_{1}}\ldots\alpha_{\lo s+2\po_{q}}=
\\ &
= \alpha_{1_{0}}\alpha_{1_{1}}\ldots\alpha_{1_{q}}\beta_{1^{0}}\beta_{1^{1}}\ldots\beta_{1^{p}} \\
& = f\lo\alpha_{1_{0}}\po f\lo\alpha_{1_{1}}\po\ldots
f\lo\alpha_{1_{q}}\po f\lo\beta_{1^{0}}\po
f\lo\beta_{1^{1}}\po\ldots f\lo\beta_{1^{p}}\po \\
& = f\lo\alpha_{1_{0}}\alpha_{1_{1}}\ldots\alpha_{1_{q}}\beta_{1^{0}}\beta_{1^{1}}\ldots\beta_{1^{p}}\po \\
& = f\lo
\theta\beta_{2^{0}}\beta_{2^{1}}\ldots\beta_{2^{p}}\alpha_{\lo
s+2\po_{0}}\alpha_{\lo s+2\po_{1}}\ldots\alpha_{\lo s+2\po_{q}}\po
\\
& = \theta f\lo\beta_{2^{0}}\po f\lo\beta_{2^{1}}\po\ldots
f\lo\beta_{2^{p}}\po f\lo\alpha_{\lo s+2\po_{0}}\po
f\lo\alpha_{\lo s+2\po_{1}}\po\ldots f\lo\alpha_{\lo
s+2\po_{q}}\po \\
& = \theta\beta_{2^{0}}\beta_{2^{1}}\ldots\beta_{2^{p}}\alpha_{\lo
s+2\po_{0}}\alpha_{\lo s+2\po_{1}}\ldots\alpha_{\lo s+2\po_{q}}.
\end{split}
\end{equation*}
This implies $\theta^{\prime}=\theta$.\newline

Assume (b) holds. Then we may assume (without loss of generality)
that $f\lo\alpha_{1_{j}}\po=\beta_{2^{j}}$, $f\lo\alpha_{\lo s +
2\po_{j}}\po=\beta_{1^{j}}$, for $j=0,1,\ldots q$, and
$f\lo\beta_{2^{j}}\po=\alpha_{1_{j}}$,
$f\lo\beta_{1^{j}}\po=\alpha_{\lo s^{\prime} + 2\po_{j}}$, for
$j=0,1,\ldots, p$. Moreover, we have the equalities:
\begin{equation*}
\begin{split}
&
\theta^{\prime}\beta_{2^{0}}\beta_{2^{1}}\ldots\beta_{2^{p}}\alpha_{\lo
s^{\prime}+2\po_{0}}\alpha_{\lo
s^{\prime}+2\po_{1}}\ldots\alpha_{\lo s^{\prime}+2\po_{q}}=
\\ &
=
\alpha_{1_{0}}\alpha_{1_{1}}\ldots\alpha_{1_{q}}\beta_{1^{0}}\beta_{1^{1}}\ldots\beta_{1^{p}}
\\ & = f\lo\beta_{2^{0}}\po f\lo\beta_{2^{1}}\po\ldots
f\lo\beta_{2^{p}}\po f\lo\alpha_{\lo s+2\po_{0}}\po
f\lo\alpha_{\lo s+2\po_{1}}\po\ldots f\lo\alpha_{\lo
s+2\po_{q}}\po
\\ & =
f\lo\beta_{2^{0}}\beta_{2^{1}}\ldots\beta_{2^{p}}\alpha_{\lo
s+2\po_{0}}\alpha_{\lo s+2\po_{1}}\ldots\alpha_{\lo s+2\po_{q}}\po
\\ & = f\lo
\theta^{-1}\alpha_{1_{0}}\alpha_{1_{1}}\ldots\alpha_{1_{q}}\beta_{1^{0}}\beta_{1^{1}}\ldots\beta_{1^{p}}\po
\\ & = \theta^{-1} f\lo\alpha_{1_{0}}\po
f\lo\alpha_{1_{1}}\po\ldots f\lo\alpha_{1_{q}}\po
f\lo\beta_{1^{0}}\po
f\lo\beta_{1^{1}}\po\ldots f\lo\beta_{1^{p}}\po \\
& =
\theta^{-1}\beta_{2^{0}}\beta_{2^{1}}\ldots\beta_{2^{p}}\alpha_{\lo
s^{\prime}+2\po_{0}}\alpha_{\lo
s^{\prime}+2\po_{1}}\ldots\alpha_{\lo s^{\prime}+2\po_{q}}.
\end{split}
\end{equation*}
This implies $\theta^{\prime}=\theta^{-1}$.

If (1) holds then obviously the algebras $\Lambda\lo
p,q,k,s,\theta\po$ and $\Lambda\lo
p^{\prime},q^{\prime},k^{\prime},s^{\prime},\theta^{\prime}\po$
are isomorphic. Assume (2) holds. We define a homomorphism
$f_{k,s,s^{\prime}}:\Lambda\lo p,q,k,s,\theta\po\longrightarrow
\Lambda\lo q,p,k,s^{\prime},\theta^{\prime}\po$ of $K$-algebras as
follows. Since $\gcd\lo s,k\po=1$, for each $i=1,2,\ldots,k$,
there exists exactly one $0\leqslant j\lo i\po\leqslant k\! -\! 1$
such that $i\equiv 1\! +\! j\lo i\po s \lo\md k\po$. Let $t\lo
i\po$ be the unique number such that $1\leqslant t\lo
i\po\leqslant k$ and $t\lo i\po\equiv 1+j\lo i\po s^{\prime}\lo\md
k\po$. We define $f_{k,s,s^{\prime}}\lo e_{i}\po=e_{t\lo i\po}$,
$f_{k,s,s^{\prime}}\lo\alpha_{i_{j}}\po=\beta_{\lo t\lo 1\po\! +\!
1\po^{j}}$, for $j=0,1,\ldots,q$, and
$f_{k,s,s^{\prime}}\lo\beta_{i^{j}}\po=\alpha_{t\lo i\! -\!
1\po}^{j}$, for $j=0,1,\ldots,p$. Let $i\in\lk 1,2,\ldots,k\pk$.
In order to prove that, for any element $\omega$ of the ideal
$\overline{I}^{\lo 1\po}\!\lo T_{p,q,k,\sigma_{s}},\theta\po$,
$f_{k,s,s^{\prime}}\lo\omega\po$ is an element of the ideal
$\overline{I}^{\lo 1\po}\!\lo
T_{p^{\prime},q^{\prime},k^{\prime}},\sigma_{s^{\prime}},\theta^{\prime}\po$,
we need the following four congruences:
\begin{romalist}
\item $t\lo i\! -\! 1\po\equiv t\lo i\po\! +\! s^{\prime}\! +\! 1 \lo\md k\po$,
\item $t\lo i\! +\! s\po\equiv t\lo i\! +\! s\! +\! 1\po\! +\! s^{\prime}\! +\! 1 \lo\md k\po$,
\item $t\lo i\! +\! s\! +\! 1\po\equiv t\lo i\po\! -\! 1 \lo\md
k\po$,
\item $t\lo i\! +\! s\po\equiv t\lo i\! -\! 1\po\! -\! 1 \lo\md
k\po$.
\end{romalist}
Observe that it is enough prove the congruences (i) and (iii).

Let $m=m_{k,s}$ and $m^{\prime}=m_{k,s^{\prime}}$. We first prove
that
\begin{displaymath}
ss^{\prime}+s+s^{\prime}\equiv 0 \lo \md k\po.
\end{displaymath}
The congruences
\begin{displaymath}
m\lo s +\! 2\po\! +\! 1\equiv 0\lo \md k\po \mbox{ , }
m^{\prime}\lo s^{\prime} +\! 2\po\! +\! 1\equiv 0\lo \md k\po
\end{displaymath}
are the congruences
\begin{equation}\label{gwiazdka}
ms\! +\!2m\! +\! 1\equiv 0\lo \md k\po \mbox{ , }
m^{\prime}s^{\prime} +\! 2m^{\prime}\! +\! 1\equiv 0\lo \md k\po.
\end{equation}
The conditions $\gcd\lo s,k\po=1$ and $\gcd\lo s^{\prime},k\po=1$
imply the congruences
\begin{equation}\label{dwiegwiazdki}
mss^{\prime}\! +\!2ms^{\prime}\! +\! s^{\prime}\equiv 0\lo \md
k\po \mbox{ , } m^{\prime}s^{\prime}s +\! 2m^{\prime}s\! +\!
s\equiv 0\lo \md k\po.
\end{equation}
Moreover, $m\! +\! m^{\prime}=k\! -\! 1$ implies $m\! +\!
m^{\prime}\equiv -1\lo \md k\po$, and hence
\begin{displaymath}
m^{\prime}\equiv -m\! -\! 1\lo\md k\po.
\end{displaymath}
From (\ref{gwiazdka}) we then obtain
\begin{displaymath}
ms\! +\!2m\! +\! 1\equiv 0\lo \md k\po \mbox{ , } \lo -m\! -\!
1\po s^{\prime} +\! 2\lo -m\! -\! 1\po\! +\! 1\equiv 0\lo \md
k\po,
\end{displaymath}
and hence
\begin{equation}\label{jedynka}
m\lo s\! -\! s^{\prime}\po\equiv s^{\prime}\lo\md k\po.
\end{equation}
From (\ref{dwiegwiazdki}) we obtain
\begin{displaymath}
mss^{\prime}\! +\!2ms^{\prime}\! +\! s^{\prime}\equiv 0\lo \md
k\po \mbox{ , } \lo -m\! -\! 1\po ss^{\prime} +\! 2\lo -m\! -\!
1\po s\! +\! s\equiv 0\lo \md k\po,
\end{displaymath}
and hence
\begin{equation}\label{dwojka}
ss^{\prime}\! +\! s\! -\! s^{\prime}\! +\! 2m\lo s\! -\!
s^{\prime}\po\equiv 0\lo\md k\po.
\end{equation}
Combining (\ref{jedynka}) and (\ref{dwojka}) we get the required
congruence
\begin{displaymath}
ss^{\prime}+s+s^{\prime}\equiv 0 \lo \md k\po.
\end{displaymath}
Further, $i\equiv 1\! +\! j\lo i\po\lo\md k\po$ and $i\! -\! 1
\equiv 1\! +\! j\lo i\! -\! 1\po\lo\md k\po$ imply
\begin{displaymath}
\lo j\lo i\po\! -\! j\lo i\! -\! 1\po\po s\equiv 1\lo\md k\po.
\end{displaymath}
Since $\gcd\lo s^{\prime},k\po=1$, this is equivalent to
\begin{displaymath}
\lo j\lo i\po\! -\! j\lo i\! -\! 1\po\po ss^{\prime}\equiv
s^{\prime}\lo\md k\po.
\end{displaymath}
Applying $ss^{\prime}+s+s^{\prime}\equiv 0 \lo \md k\po$, we then
get
\begin{displaymath}
\lo j\lo i\po\! -\! j\lo i\! -\! 1\po\po ss^{\prime}\equiv
-ss^{\prime}\! -\! s\lo\md k\po.
\end{displaymath}
Since $\gcd\lo s,k\po=1$, this implies that
\begin{displaymath}
\lo j\lo i\po\! -\! j\lo i\! -\! 1\po\po s^{\prime}\equiv
-s^{\prime}\! -\! 1\lo\md k\po.
\end{displaymath}
Finally, $t\lo i\po\equiv 1\! +\! j\lo i\po s^{\prime}\lo\md k\po$
and $t\lo i\! -\! 1\po\equiv 1\! +\! j\lo i\! -\! 1\po
s^{\prime}\lo\md k\po$ imply
\begin{displaymath}
t\lo i\po\! -\! t\lo i\! -\! 1\po\equiv \lo j\lo i\po\! -\! j\lo
i\! -\! 1\po\po s^{\prime}\lo\md k\po,
\end{displaymath}
and consequently $t\lo i\po\! -\! t\lo i\! -\! 1\po\equiv
-s^{\prime}\! -\! 1\lo\md k\po$. Therefore, we obtain
\begin{displaymath}
t\lo i\! -\! 1\po\equiv t\lo i\po\! +\! s^{\prime}\! +\! 1 \lo\md
k\po,
\end{displaymath}
as required in (i).

Let $i\equiv 1\! +\! j\lo i\po\lo\md k\po$ and $i\! +\! s\! +\! 1
\equiv 1\! +\! j\lo i\! +\! s\!+\! 1\po\lo\md k\po$. This implies
\begin{displaymath}
\lo j\lo i\po\! -\! j\lo i\! +\! s\!+\! 1\po\po s\equiv -s\! -\!
1\lo\md k\po.
\end{displaymath}
Since $\gcd\lo s^{\prime},k\po=1$, this is equivalent to
\begin{displaymath}
\lo j\lo i\po\! -\! j\lo i\! +\! s\!+\! 1\po\po ss^{\prime}\equiv
-ss^{\prime}\! -\! s^{\prime}\lo\md k\po.
\end{displaymath}
Applying $ss^{\prime}+s+s^{\prime}\equiv 0 \lo \md k\po$, we then
get
\begin{displaymath}
\lo j\lo i\po\! -\! j\lo i\! +\! s\!+\! 1\po\po ss^{\prime}\equiv
s\lo\md k\po.
\end{displaymath}
Since $\gcd\lo s,k\po=1$, this implies that
\begin{displaymath}
\lo j\lo i\po\! -\! j\lo i\! +\! s\!+\! 1\po\po s^{\prime}\equiv
1\lo\md k\po.
\end{displaymath}
Finally, $t\lo i\po\equiv 1\! +\! j\lo i\po s^{\prime}\lo\md k\po$
and $t\lo i\! +\! s\!+\! 1\po\equiv 1\! +\! j\lo i\! +\! s\!+\!
1\po s^{\prime}\lo\md k\po$ imply
\begin{displaymath}
t\lo i\po\! -\! t\lo i\! +\! s\!+\! 1\po\equiv \lo j\lo i\po\! -\!
j\lo i\! +\! s\!+\! 1\po\po s^{\prime}\lo\md k\po,
\end{displaymath}
and consequently $t\lo i\po\! -\! t\lo i\! +\! s\!+\! 1\po\equiv
1\lo\md k\po$. Therefore, we obtain
\begin{displaymath}
t\lo i\! +\! s\! +\! 1\po\equiv t\lo i\po\! -\! 1 \lo\md k\po,
\end{displaymath}
as required in (iii).

We note that $f_{s,s^{\prime},k}\circ
f_{s^{\prime},s,k}=\id_{\Lambda\lo
p^{\prime},q^{\prime},k^{\prime},s^{\prime},\theta^{\prime}\po}$
and $f_{s^{\prime},s,k}\circ f_{s,s^{\prime},k}=\id_{\Lambda\lo
p,q,k,s,\theta\po}$. Hence the algebras $\Lambda\lo
p,q,k,s,\theta\po$ and $\Lambda\lo
p^{\prime},q^{\prime},k^{\prime},s^{\prime},\theta^{\prime}\po$
are isomorphic.
\end{proof}

We do not know yet if two algebras of the forms $\Lambda\lo
p,q,k,s,\theta\po$ and $\Lambda\lo
p^{\prime},q^{\prime},k^{\prime},s^{\prime},\theta^{\prime}\po$
are stably equivalent if and only if they are isomorphic.

We ends this section with two examples illustrating the above
considerations.

\begin{exa} \rm
Let $k=8$. Then the set of all $s\in \{1,2,\ldots,7\}$ satisfying
$\gcd\lo s,8\po=1$ and $\gcd\lo s\! +\! 2,8\po=1$ consists of $1$,
$3$, $5$ and $7$. Moreover, we have $m_{8,1}=5$, $m_{8,3}=3$,
$m_{8,5}=1$, $m_{8,7}=7$. Observe that for any $s,s^{\prime}\in\lk
1,3,5,7\pk$ we have $m_{8,s}\! +\! m_{8,s^{\prime}}\neq k\! -\!
1=7$. Therefore, for any $p,q,p^{\prime},q^{\prime}\geqslant 0$,
$s,s^{\prime}\in\lk 1,3,5,7\pk$ with $s\neq s^{\prime}$ and
$\theta\in K\backslash\lk 0\pk$, we have $\Lambda\lo
p,q,8,s,\theta\po\ncong\Lambda\lo
p^{\prime},q^{\prime},8,s^{\prime},\theta\po$.
\end{exa}

\begin{exa} \rm
Let $k=5$, $s=1$ and $s^{\prime}=2$. Then $\gcd\lo s,k\po=1$,
$\gcd\lo s\! +\! 2,k\po=1$, $\gcd\lo s^{\prime},k\po=1$, $\gcd\lo
s^{\prime}\! +\! 2,k\po=1$. Moreover, $m_{5,1}=3$, $m_{5,2}=1$,
and so $m_{5,1}\! +\! m_{5,2}=4=k\! -\! 1$. Hence, for any
$p,q,p^{\prime},q^{\prime}\geqslant 0$ and $\theta\in
K\backslash\lk 0\pk$, we have $\Lambda\lo
p,q,5,1,\theta\po\cong\Lambda\lo
p^{\prime},q^{\prime},5,2,\theta\po$.
\end{exa}


\section{\hspace{-2.5ex}. Proof of Theorem \ref{mainthm2}} \label{secrepdim}


As a consequence of our main theorem we deduce in this section
some results on homological invariants of one-parametric
selfinjective algebras.

Recall that the \textit{representation dimension} of an algebra
$A$ was defined by Auslander \cite{bib:Auslander} as the number
$$\rdim(A):=\inf \{\gldim\End_A(N)~|~N\mbox{ generator-cogenerator}\}.
$$
A finitely generated $A$-module $N$ is a generator-cogenerator for
$A$ if it contains all projective indecomposable and all injective
indecomposable $A$-modules as direct summands. This invariant
gives a way of measuring how far an algebra is from being of
finite representation type. In fact, $\rdim(A)\le 2$ if and only
if $A$ is of finite representation type \cite{bib:Auslander}.
\medskip

We shall prove the following result which contributes further
to the long-term project of determining the representation dimension
for all (selfinjective) algebras of tame representation type.

\begin{tw} \label{repdim}
Let $A$ be a standard one-parametric selfinjective algebra. Then
$A$ has representation dimension 3.
\end{tw}

\begin{proof}
We shall use the main results of \cite{bib:EHIS}. Firstly, if an
algebra $A$ can be embedded into an algebra $B$, preserving
radicals, and if $B$ is of finite representation type, then $A$
has representation dimension at most 3 (\cite[Theorem
1.1]{bib:EHIS}). As an application, we get secondly that every
special biserial algebra has representation dimension at most 3
(\cite[(1.3)]{bib:EHIS}). We also use that for selfinjective
algebras, the representation dimension is invariant under derived
equivalences \cite{bib:Xi}.

For $A$ weakly symmetric, the theorem is proved in our paper
\cite{bib:BoHoSk2}. Therefore, it suffices to prove the theorem
for the normal forms of our derived equivalence classification for
standard one-parametric selfinjective algebras obtained in this
paper.

The normal forms $\Lambda(p,q,k,s,\theta)$ are special biserial,
hence they have representation dimension at most 3. Since they are
not of finite representation type, the representation dimension is
actually 3. The normal forms $\Gamma^*(n)$ are not special
biserial. In order to show that $\rdim\lo \Gamma^*(n)\po=3$, we
shall construct a radical embedding of a certain factor algebra of
$\Gamma^*(n)$ into an algebra of finite representation type. This
proof is similar to the one for algebras of Euclidean type carried
out in \cite{bib:BoHoSk2}. Since this paper might not be easily
available, we include a complete proof here for the convenience of
the reader.

Suitable factor algebras for determining representation dimensions
are suggested by the following useful result
(\cite[(1.2)]{bib:EHIS}): {\em Let $A$ be a basic algebra, and let
$P$ be an indecomposable projective-injective $A$-module. Consider
the factor algebra $B:=A/soc(P)$ modulo the socle of $P$. If
$repdim\lo B\po\le 3$, then also $repdim\lo A\po\le 3$.}

We now consider the algebras $\Gamma^*(n)$. Let $v$ be the unique
vertex of the quiver of $\Gamma^*(n)$ with three incoming and
outgoing arrows, with corresponding projective indecomposable
module $P(v)$. We consider the factor algebra
$\tilde{\Gamma}^*(n):=\Gamma^*(n)/\soc\lo P_v\po$. We get
additional relations $\alpha_1\alpha_2=0$, $\gamma_1\gamma_2=0$
and $(\beta_1\ldots\beta_n)^2=0$. We have to show that
$\tilde{\Gamma}^*(n)$ has representation dimension 3.

In order to construct a radical embedding of $\tilde{\Gamma}^*(n)$
into an algebra of finite representation type we use the technique
of {\em splitting datum} as introduced in \cite{bib:EHIS}. We can
distribute the edges adjacent to $v$ as follows: $E_1\cup
E_2:=\{\beta_n\}\cup \{\alpha_2,\gamma_2\}$ and $S_1\cup
S_2:=\{\beta_1\}\cup \{\alpha_1,\gamma_1\}$. The conditions for a
splitting datum are satisfied: all products in $E_iS_j$ are zero
for $i\neq j$ and all relations of $\tilde{\Gamma}^{*}(n)$ are
monomial. Any splitting datum gives rise to a radical embedding
into an algebra whose quiver is obtained by splitting the vertex
under consideration. In our case, we get a radical embedding of
$\tilde{\Gamma}^{*}(n)$ into the algebra
$\tilde{\Gamma}^{*}(n)^{sp}$ given by the disjoint union of a
cyclic quiver with edges $\beta_1,\ldots,\beta_n$ and relations
$(\beta_1\ldots\beta_n)^2=0$, $\beta_j\beta_{j+1}\ldots
\beta_n\beta_1\ldots\beta_n\beta_1\ldots\beta_{j-1}\beta_j=0$, for
$2\le j\le n$, and the quiver


\begin{center}

\setlength{\unitlength}{2368sp}%
\begingroup\makeatletter\ifx\SetFigFont\undefined
\def\x#1#2#3#4#5#6#7\relax{\def\x{#1#2#3#4#5#6}}%
\expandafter\x\fmtname xxxxxx\relax \def\y{splain}%
\ifx\x\y   
\gdef\SetFigFont#1#2#3{%
  \ifnum #1<17\tiny\else \ifnum #1<20\small\else
  \ifnum #1<24\normalsize\else \ifnum #1<29\large\else
  \ifnum #1<34\Large\else \ifnum #1<41\LARGE\else
     \huge\fi\fi\fi\fi\fi\fi
  \csname #3\endcsname}%
\else
\gdef\SetFigFont#1#2#3{\begingroup
  \count@#1\relax \ifnum 25<\count@\count@25\fi
  \def\x{\endgroup\@setsize\SetFigFont{#2pt}}%
  \expandafter\x
    \csname \romannumeral\the\count@ pt\expandafter\endcsname
    \csname @\romannumeral\the\count@ pt\endcsname
  \csname #3\endcsname}%
\fi
\fi\endgroup
\begin{picture}(3055,615)(2692,-1861)
\thinlines \put(4201,-1561){\circle{100}}
\put(5701,-1561){\circle{100}} \put(2851,-1486){\vector( 1,
0){1200}} \put(3976,-1636){\vector(-1, 0){1125}}
\put(4351,-1486){\vector( 1, 0){1200}}
\put(5476,-1636){\vector(-1, 0){1125}}
\put(2738,-1561){\circle{100}}
\put(3280,-1400){\makebox(0,0)[lb]{\smash{\SetFigFont{10}{12.0}{rm}
$\alpha_2$}}}
\put(3280,-1861){\makebox(0,0)[lb]{\smash{\SetFigFont{10}{12.0}{rm}
$\alpha_1$}}}
\put(4800,-1411){\makebox(0,0)[lb]{\smash{\SetFigFont{10}{12.0}{rm}
$\gamma_1$}}}
\put(4800,-1861){\makebox(0,0)[lb]{\smash{\SetFigFont{10}{12.0}{rm}
$\gamma_2$}}}
\end{picture}
\end{center}


\noindent with relations $\alpha_1\alpha_2=0$,
$\gamma_1\gamma_2=0$, $\alpha_2\gamma_1=0$, $\gamma_2\alpha_1=0$.
Note that this new algebra $\tilde{\Gamma}^{*}(n)^{sp}$ is special
biserial and of finite representation type \cite{bib:SkWa}.

By the radical embedding theorem from \cite{bib:EHIS} we can deduce
that
$\tilde{\Gamma}^*(n)$ has representation dimension 3, as desired.
\end{proof}

\bigskip

For selfinjective algebras, the representation dimension
is intimately related to the
notion of the dimension of a triangulated category as recently
defined by R. Rouquier \cite{Rouquier}. Let us briefly recall
the main definitions.
Let $A$ be a selfinjective algebra with stable module
category $\underline{\md}\,A$ of finitely generated
$A$-modules.
Recall that this is a triangulated category with shift
$\Omega^{-1}$ (see \cite{bib:Hap2}). For any $A$-module $M$ denote
by $\langle M \rangle_i$ the full subcategory of
$\underline{\md}\,A$ with objects the direct summands of the
modules obtained from $M$ by taking $i$-fold extensions of finite
direct sums of shifts. Then following R. Rouquier \cite{Rouquier},
the \textit{dimension} of $\underline{\md}\,A$, denoted
$\dim\underline{\md}\,A$, is the minimal integer $d\ge 0$ such
that there is an $A$-module $M$ with $\underline{\md}\,A=\langle
M\rangle_{d+1}$. In particular, we have $\dim\underline{\md}\,A=0$
if and only if $A$ has finite representation type. As with the
representation dimension, it seems to be very hard to determine
this dimension of the stable module category, even for examples of
fairly small and well-known algebras.

%
%

The representation dimension of a selfinjective algebra is closely
related to the dimension of its stable module category, as
observed by R. Rouquier.

\begin{proposition}
Let $A$ be a selfinjective algebra. Then $\rdim(A)\ge 2+
\dim\underline{\md}\,A.$
\end{proposition}
\begin{proof}
See \cite[Proposition 6.9]{Rouquier}.
\end{proof}

Using this connection we can deduce the dimension of the
stable module category of standard selfinjective
one-parametric algebras from
our Theorem \ref{repdim}.

\begin{corollary} Let $A$ be a standard one-parametric selfinjective
algebra. Then the stable module category $\underline{\md}\,A$ has
dimension 1.
\end{corollary}

It does not seem to be clear how one can determine
the dimensions of stable module categories directly,
without referring to the representation dimension.

The following problem would be interesting to consider.
\medskip

\noindent
{\bf Question 1:} Is $\dim\underline{\md}\,A\le 1$ for all
(selfinjective) algebras of tame representation type?
\medskip

This is a weaker version of the following long standing open
problem.
\medskip

\noindent
{\bf Question 2:} Is $\rdim(A)\le 3$ for all
(selfinjective) algebras of tame representation type?
\medskip

Note that the only known examples of algebras of representation
dimension greater 3 were only recently discovered by Rouquier
\cite{Rouquier}, and they are all of wild representation type.

\bigskip\noindent{\bf Acknowledgements.} The first and the third
named authors gratefully acknowledge support from the Polish
Scientific Grant KBN No. 1 P03A 018 27.


\end{document}